\documentclass[leqno,12pt]{article}
\usepackage{amsmath}
\usepackage{amssymb}
\usepackage{amsthm}
\usepackage{latexsym}

\DeclareMathOperator*{\esup}{ess\,sup}
\DeclareMathOperator{\HS}{HS}
\DeclareMathOperator{\IHS}{IHS}
\DeclareMathOperator{\sgn}{sgn}
\DeclareMathOperator{\supp}{supp}
\DeclareMathOperator{\Tr}{Tr}

\newcommand{\blim}{\nolimits} %`Below limits'.  An integral over a region $R$ will be written as $$\int\blim_R$$$; thus, if \blim is redefined to return \limits, then an integral over a region will appear with the region beneath (rather than beside) the integral sign.
\newcommand{\tblim}{\nolimits} %`Text below limits'.  As above, but in inline mode.
\newcommand{\bblim}{\nolimits} %`Below double limits'.  An integral over two regions $R$ and $S$ will be written as $$\int\bblim_{\substack{R \\ S}}$$; thus, if \bblim is redefined to return \limits, then an integral over two regions will appear with the region beneath (rather than beside) the integral sign.
\newcommand{\inlim}{\nolimits} %`(Text) indexing limits'.  This controls whether an indexing subscript (as in $\max\inlim_x$) appears beside or below the operator it indexes. 
\newcommand{\ic}{\!} %`Integral correction'.  This is used in the context $$2\ic\int$$, where usually the distance between the 2 and the $\int$ would be too large.
\newcommand{\bic}{\!\!\!} %`Between integrals correction'.  This is used for when we have a multiple (i.e., double, triple &c.) integral, when usually too much space is inserted between the two integrals.
\newcommand{\tbic}{\!\!} %`Text between integrals correction'.  As above, but in inline mode.
\newcommand{\tic}{\!} %`Text integral correction'.  This is used in the context $2\tic\int$ (notice that this is an inline equation), as above.
\newcommand{\dic}{\!\!} %`Derivative integral correction'.  This is used for when we are taking the derivative of an expression which begins with an integral.

\newcommand{\FracInt}{{\mathbb I}}
\newcommand{\ksub}{k}%%At one point in the text, $k$ seems to be used to represent two different things.  Changing the symbol returned by this command will change all occurrences that (seemed to me to) mean one thing.
\newcommand{\mult}{\!\cdot\!}
\newcommand{\nm}[2]{||#1||_{H^{#2}}}
\newcommand{\dnm}[2]{||#1||_{\dot{H}^{#2}}}
\newcommand{\anm}[3]{||#1||_{H^{#2}(#3)}}
\newcommand{\adnm}[3]{||#1||_{\dot{H}^{#2}(#3)}}
\newcommand{\snm}[3]{||#1||_{H^{#2}_{#3}}}
\newcommand{\asnm}[4]{||#1||_{H^{#2}_{#4}(#3)}}
\newcommand{\dsnm}[3]{||#1||_{\dot{H}^{#2}_{#3}}}
\newcommand{\adsnm}[4]{||#1||_{\dot{H}^{#2}_{#4}(#3)}}
\newcommand{\pd}[2]{\frac{\partial #1}{\partial #2}}

\newcommand{\pref}[1]{(\ref{#1})}

\newcommand{\nw}[1]{#1} %This command is used so that newly inserted blocks of text can easily be found in the TeX code.
\newcommand{\old}[1]{}  %This swallows its argument, an old text removed in the new draft, without deleting it from the code.
\newcommand{\comment}[1]{}  %This swallows its argument, a comment on the content.  This should be replaced by something like {\forcemath\verbatim\bf #1} (the command can appear in math mode, and/or contain TeX commands) before printing a draft, so that typist comments will be seen by the reader.

\theoremstyle{plain}
\newtheorem{corollary}{Corollary}[section]
\newtheorem{lemma}{Lemma}[section]
\newtheorem{proposition}{Proposition}[section]
\newtheorem{theorem}{Theorem}[section]

\theoremstyle{definition}
\newtheorem{definition}{Definition}[section]

\theoremstyle{remark}
\newtheorem{fact}{Fact}[section]
\newtheorem{remark}{Remark}[section]
\newtheorem{claim}[equation]{Claim}

\newcounter{sublemma}
\newcounter{subtheorem}

\renewcommand{\thelemma}{\arabic{section}.\arabic{lemma}}
\renewcommand{\thesublemma}{\thelemma.\arabic{sublemma}}

\renewcommand{\thetheorem}{\arabic{section}.\arabic{theorem}}
\renewcommand{\thesubtheorem}{\thetheorem.\arabic{subtheorem}}
\newcommand{\newsect}[2]{\section{#1}
\setcounter{equation}{0}
\setcounter{corollary}{0}
\setcounter{definition}{0}
\setcounter{lemma}{0}
\setcounter{proposition}{0}
\setcounter{remark}{0}
\setcounter{theorem}{0}
\label{#2}} %This command starts and labels a new section, and `reinitialises' all the relevant counters.

\begin{document}
\begin{center}
\large\bf THE GENERALIZED KORTEWEG-DE VRIES \\
EQUATION ON THE HALF LINE
\end{center}

\centerline{\sc J. E. Colliander \footnote{\rm J.\ E.\ C.\ was supported in part by an N.\ S.\ F.\ Postdoctoral Research Fellowship and N.\ S.\ F.\ grant DMS 0100595.} and C. E. Kenig \footnote{\rm C.\ E.\ K.\ was supported in part by N.\ S.\ F.\ Grant DMS 9500725.}}

\begin{center}
\begin{minipage}{4in}
\footnotesize {\sc Abstract}.  The initial-boundary value problem for the generalized Korteweg-de Vries equation on a half-line is studied by adapting the initial value techniques developed by Kenig, Ponce and Vega and Bourgain to the initial-boundary setting.  The approach consists of replacing the initial-boundary problem by a forced initial value problem.  The forcing is selected to satisfy the boundary condition by inverting a Riemann-Liouville fractional integral.
\end{minipage}
\end{center}

\newsect{Introduction}{1}

This paper introduces a method to solve initial-boundary value problems for nonlinear dispersive partial differential equations by recasting these problems as initial value problems with an appropriate forcing term.  This reformulation transports the robust theory of initial value problems to the initial-boundary value setting.  The procedure is applied here to solve the initial-boundary value problem for the generalized Korteweg-de Vries equation on the half-line.  We expect that generalizations of the ideas described below will be useful in solving problems in higher dimensions, such as the nonlinear Schr\"odinger equation posed on a spatial domain under Dirichlet boundary conditions.  The methods introduced here may be viewed as a dispersive generalization of the classical Calder\'on Projector method \cite{Calderon} used in elliptic and parabolic problems.

Consider the {\it initial-boundary value problem for the generalized Korteweg-de Vries equation on the right half-line,}
\begin{equation}
\label{1.1}
\begin{cases}
\partial_tu + \partial_x^3u + \frac{1}{k + 1}\partial_xu^{k + 1} = 0, & x > 0, t \in [0, T]\text{ and }k \in {\mathbb N} \\
u(x, 0) = \phi(x),                                                    & x \geq 0     \\
u(0, t) = f(t),                                                       & t \in [0, T]
\end{cases}
\end{equation}
with $\phi$ and $f$ satisfying certain regularity hypotheses, and $T > 0$.  This problem will be solved using a {\it forced initial value problem for the generalized Korteweg-de Vries equation on the line,}
\begin{equation}
\label{1.2}
\begin{cases}
\partial_t\tilde{u} + \partial_x^3\tilde{u} + \frac{1}{k + 1}\partial_x\tilde{u}^{k + 1} = \delta_0(x)g(t), & x \in {\mathbb R}\text{ and }t \in [0, T] \\
\tilde{u}(x, 0) = \tilde{\phi}(x), & x \in {\mathbb R},
\end{cases}
\end{equation}
where $\delta_0$ denotes the Dirac mass at $x = 0$ and $\tilde{\phi}$ is a nice extension of $\phi$.  The {\it boundary forcing function} $g$ is selected to ensure that
\begin{equation}
\label{1.3}
\tilde{u}(0, t) = \tilde{f}(t), \quad t \in [0, T],
\end{equation}
where $\tilde{f}$ is a nice extension of $f$.  We will sometimes refer to \pref{1.3} as a boundary condition, even though $\{x = 0\}$ is not a boundary for the $\tilde{u}$ problem \pref{1.2}.  Upon constructing the solution $\tilde{u}$ of \pref{1.2}, we obtain the solution $u$ of \pref{1.1} by restriction, as
\begin{equation}
\label{1.4}
u = \tilde{u}\big|_{\{x \geq 0\} \times \{t : t \in [0, T]\}}.
\end{equation}
The solution of \pref{1.2} satisfying \pref{1.3} is constructed using the initial value machinery in \cite{5}, \cite{10} and \cite{11} and the inversion of a Riemann-Liouville fractional integration operator.  First, an explicit solution of the linearization of \pref{1.1}, using the linearization of \pref{1.2} satisfying \pref{1.3}, is constructed.  This reveals the natural regularity relationships among $\phi$, $f$ and $g$.  We verify that the space-time norms (used in \cite{5} and \cite{10}) of the linear solution $\tilde{u}$ are bounded by related norms of $\tilde{\phi}$ and $\tilde{f}$.  Next, we consider \pref{1.2} with the nonlinear term replaced by a given function $h$ of space-time.  The solution of the resulting inhomogeneous analogue of \pref{1.2} satisfying \pref{1.3} is explicitly constructed in terms of $\tilde{\phi}$, $\tilde{f}$ and $h$, and the natural inhomogeneous estimates from \cite{5} and \cite{10} are established.  Finally, the nonlinear term is treated as an inhomogeneity and a multilinear estimate exploiting time localization and/or scaling proves the required contraction estimate.

The linearization of \pref{1.2}, obtained by removing the nonlinear term $\frac{1}{k + 1}\partial_xu^{k + 1}$, is solved using Duhamel's formula,
\begin{equation}
\label{1.5}
u(x, t) = S(t)\tilde{\phi}(x) + \int_0^t S(t - t')\delta_0(x)g(t')dt',
\end{equation}
where $S$ is the linear solution operator given by
\begin{equation}
\label{1.6}
S(t)\phi(x) = \int e^{i(x\xi + t\xi^3)}\Hat{\Tilde{\phi}}(\xi)d\xi = \int t^{-1/3}A\biggl(\frac{x - x'}{t^{1/3}}\biggr)\tilde{\phi}(x')dx',
\end{equation}
where $A$ is the Airy function.  We describe how to select the boundary forcing $g$ in this context and reveal the natural regularity relationships among $\tilde{\phi}$, $\tilde{f}$ and $g$.  The condition \pref{1.3} determines the boundary forcing function if we can solve
\begin{equation}
\label{1.7}
\int_0^t S(t - t')\delta_0(x)g(t')dt'\Bigr|_{\{x = 0\}} = \tilde{f}(t) - S(t)\tilde{\phi}(x)\Bigr|_{\{x = 0\}}
\end{equation}
for $g$.  The right-hand side of \pref{1.7} suggests that the regularity properties of $\tilde{f}$ should be the same as those expressed by $S(t)\tilde{\phi}\bigl|_{\{x = 0\}}$.  For $\tilde{\phi} \in \dot{H}^s({\mathbb R}_x) = \bigl\{f : |\xi|^s\hat{f}(\xi) \in L^2({\mathbb R}_\xi)\bigr\}$, the change of variables argument used to prove the local smoothing property shows that $S(t)\tilde{\phi}\bigl|_{\{x = 0\}} \in \dot{H}^{(s + 1)/3}({\mathbb R}_t)$.  The convolution representation of $S(t)$ allows the left-hand side of \pref{1.7} to be reexpressed as $$\int_0^t \frac{1}{(t - t')^{1/3}}A\biggl(\frac{x}{(t - t')^{1/3}}\biggr)g(t')dt'.$$  Since the Airy function $A$ is continuous and $A(0) \neq 0$, we can evaluate at $x = 0$ and rewrite \pref{1.7} as $$A(0)\ic\int_0^t \frac{1}{(t - t')^{1/3}}g(t')dt' = \tilde{f}_1(t),$$ where we have introduced the notation $\tilde{f}_1(t) = \tilde{f}(t) - S(t)\tilde{\phi}\bigl|_{\{x = 0\}} \in \dot{H}^{(s + 1)/3}_0({\mathbb R}_t)$.  Recalling the theory in \cite{14} (see also \ref{3} below) of the Riemann-Liouville fractional integral
\begin{equation}
\label{1.8}
\FracInt_\alpha(h)(t) = \frac{1}{\Gamma(\alpha)}\ic\int_0^t (t - s)^{\alpha - 1}h(s)ds,
\end{equation}
we reexpress \pref{1.7} as
\begin{equation}
\label{1.9}
\FracInt_{2/3}(g) = \frac{1}{A(0)\Gamma(2/3)}\tilde{f}_1.
\end{equation}
This equation is solved for $g$ by applying $\FracInt_{-2/3}$.  Moreover, the operator $\FracInt_{\alpha}$ is smoothing of order $-\alpha$, so the boundary forcing function $g$ is seen to have $2/3$ fewer derivatives in $L^2$ than $\tilde{f}_1$.  These remarks reveal that, for $\tilde{\phi} \in H^s({\mathbb R}_x)$, it is natural to assume that $\tilde{f} \in H^{(s + 1)/3}({\mathbb R}_t)$ and to look for $g \in H^{(s - 1)/3}({\mathbb R}_t)$.  Certain technical distinctions among the spaces $H^\sigma({\mathbb R}_t)$ and $\dot{H}^\sigma({\mathbb R}_t)$ will be clarified below using time localization.

We apply this basic idea -- {\it apply boundary forcing to ensure that the $\{x = 0\}$ condition $u(0, t) = f(t)$ holds with the forcing $g$ selected by inverting a certain fractional integral} -- in the nonlinear setting using a contraction mapping argument to obtain the following results for the initial-boundary value problem \pref{1.1}.  Let ${\mathbb R}_y^+$ denote the open right half-line $\{y : y > 0\}$.

\begin{theorem}
\label{Th1.1}
Consider \pref{1.1} in the cases $k \in {\mathbb N}$ and set $s_1 = 0$, $s_2 = 1/4$, $s_3 = 1/12$ and $s_\ksub%%Is this really supposed to be the same $k$ as above?
 = \frac{1}{2} - \frac{2}{\ksub}$ (for $\ksub \geq 4$).  For any $\phi \in H^{s_\ksub}({\mathbb R}_x^+)$ and $f \in H^{(s_\ksub + 1)/3}({\mathbb R}_t^+)$, there exists a $T = T(\phi, f) > 0$ and a solution $u \in C\bigl([0, T]; H^{s_\ksub}({\mathbb R}_x^+)\bigr)$ of the initial-boundary value problem \pref{1.1}.  When $\ksub < 4$, $T = T(||\phi||_{H^{s_\ksub}}, ||f||_{H^{(s_\ksub + 1)/3}})$.  The map $(\phi, f) \mapsto u$ taking the initial and boundary data to the solution is Lipschitz-continuous from $H^{s_\ksub}({\mathbb R}_x^+) \times H^{(s_\ksub + 1)/3}({\mathbb R}_t^+)$ to $C\bigl([0, T]; H^{s_\ksub}({\mathbb R}_x^+)\bigr)$.
\end{theorem}

The cases where $\ksub \geq 2$ are proved using the selected boundary forcing procedure outlined above in the function space framework introduced by Kenig, Ponce and Vega in \cite{10} to treat the corresponding initial value problems.  The proof in the case $\ksub = 1$ uses the forcing procedure and a contraction mapping argument in a \nw{modification of a} function space introduced by Bougain in \cite{5}.  The need for the modified spaces is explained in Remark \ref{Rk5.1'}.  An optimization of the approach in \cite{5}, carried out in \cite{11}, established local well-posedness of the initial value problem for the standard (i.e., $\ksub = 1$) KdV equation in $H^s({\mathbb R})$ for $s > -3/4$.  An improvement in the $k = 1$ case of Theorem \ref{Th1.1} to $s_1 > -3/4$ is likely to be true, but we have chosen not to carry it out in the present paper{\footnote{This improvement was completed by Justin Holmer while this paper was under review.}}.  In particular, the selected forcing procedure extends the theory of the initial value problem for the generalized KdV equations to the initial-boundary value setting.

Our local (in time) results can be combined with integration by parts to yield global (in time) results.  For instance, we have:

\begin{theorem}
\label{Th1.2}
{\ }
\begin{enumerate}
\item\label{Th1.2(i)} When $\ksub = 1$, $\phi \in L^2({\mathbb R}^+)$ and $f \in H^{7/12}({\mathbb R}^+)$, the results of Theorem \ref{Th1.1} extend to any interval $[0, T]$.
\item\label{Th1.2(ii)}  When $\ksub = 2$, $\phi \in H^1({\mathbb R}^+)$, $f \in H^{11/12}({\mathbb R}^+)$ and $\phi(0) = f(0)$, we have well-posedness  of \pref{1.1} in $H^1 \times H^{2/3}$ (in the sense of Theorem \ref{Th1.1}) on any interval $[0, T]$, provided that $||f||_{L^2({\mathbb R}^+)}$ is sufficiently small.
\item\label{Th1.2(iii)}  When $\ksub = 3$, the result in \pref{Th1.2(ii)} holds for $f \in H^{5/4}({\mathbb R}^+)$, provided that $||f||_{L^2({\mathbb R}^+)}$ is sufficiently small.
\item\label{Th1.2(iv)}  For $\ksub \geq 4$, the result in \pref{Th1.2(ii)} holds for $f \in H^{11/12}({\mathbb R}^+)$, provided that $||f||_{L^2({\mathbb R}^+)}$ and $||\phi||_{L^2({\mathbb R}^+)}$ are sufficiently small.
\end{enumerate}
\end{theorem}

\nw{(See \ref{7}.)  Note that $7/12 > 1/3$; $11/12 > 2/3$; and $5/4 > 2/3$.  This is because ``conservation laws'' are not ``exact'' in the quarter-plane setting.  (Again, see \ref{7} and \cite{3}.)  \comment{Should this be part of the theorem?}}

The boundary forcing method is flexible.  In particular, Theorem \ref{Th1.1} applies to the initial-boundary value problem obtained by replacing $\partial_x^3$ by $\partial_x^3 + c\partial_x$ for some $c \in {\mathbb R}$ in \pref{1.1}.  \old{The first-order transport term $c\partial_xu$ may be treated by noting that, for the local (in time) theory (and $T < \infty$), the multipliers $e^{it\xi^3}$ and $e^{i(t\xi^3 + ct\xi)}$ satisfy the same estimates.}  \nw{See Remarks \ref{Rk6.4}\comment{What Remark 6.4?} and \ref{Rk7.5}.  (Transport terms arise naturally in physical models.)}  The method also constructs solutions to the analogous initial value problem on the left half-line $\{x < 0\}$.  (Uniqueness in the linear problem in the left half-line %, and for transport terms like $c\partial_x$ for $c < 0$, are 
is not established here, and indeed does not hold.  See \cite{F}, \cite{FS}.  Also, the global (in time) results of Theorem \ref{Th1.2} do not apply to the left half-line%, or to transport terms like $c\partial_x$ for $c < 0$
.)  We plan to apply this approach to the generalized KdV problem on a finite interval with Dirichlet boundary conditions in a future publication.  Recent work of  Colin and Ghidaglia on the problem on a finite interval has just appeared \cite{CG2001}. Initial-boundary value problems for other dispersive equations in one spatial dimension may also be studied using the method and the estimates obtained in \cite{12%%Is this the right reference?  It's hard to read.
}.  We plan to extend these ideas to treat higher-dimensional initial-boundary value problems, such as NLS on a domain $\Omega \subseteq {\mathbb R}^2$ with Dirichlet boundary conditions (see \cite{6} and the references therein).  We believe that this paper is the first application of a general method for treating initial-boundary value problems for nonlinear dispersive PDEs.

The initial-boundary value problem \pref{1.1} (with a first-order transport term) physically models the evolution of small-amplitude shallow water long waves propagating in a channel with forcing applied at the left end (see \cite{9}).  Bona and Winther (in \cite{2} and \cite{4}), employing energy methods, treated the right half-line problem \pref{1.1} with a transport term in the case $\ksub%%Or just $k$?
 = 1$ with more regular initial and boundary data ($(\phi, f) \in H^s({\mathbb R}_x) \times H^{1 + [s/3]}({\mathbb R}_t^+)$).  Polynomial generalizations of the standard KdV equation were studied by Bona and Luo in \cite{3} under higher regularity assumptions on the data.  A recent preprint of Bona, Sun and Zhang (\cite{1}) uses a Laplace transform technique in the framework of function spaces developed in \cite{10} to treat the $\ksub = 1$ version of the right half-line problem \pref{1.1} with transport term for initial data $\phi \in H^s({\mathbb R}^+)$ and boundary data $f \in H^{(s + 1)/3}({\mathbb R}^+)$ provided $s > 3/4$. The half-line problem for KdV (and generalizations) has also been considered by Faminski{\u\i} \cite{F88}, \cite{F96}, \cite{F}.

We outline the rest of the paper.  The next two sections present background material required for the proof of Theorem \ref{Th1.1}.  \ref{2} develops the theory of $L^2$-based Sobolev spaces on half-lines and intervals.  \ref{3} recalls the Riemann-Liouville fractional integral and establishes its mapping properties between $L^2$-based Sobolev spaces.  \ref{4} provides estimates in space-time mixed-norm spaces for the linear solution operator and the inhomogeneous and boundary forcing Duhamel terms\old{ and contains a uniqueness result}.  \ref{5} establishes similar space-time estimates in Bourgain spaces.  In \ref{6}, we construct and estimate solutions of homogeneous and inhomogeneous linear analogues of \pref{1.2}; \ref{6} also contains a uniqueness result.  The nonlinear problems are addressed in \ref{7}.

\newsect{Sobolev Spaces on Half-Lines and Intervals}{2}

This section develops the $L^2$-based Sobolev spaces on half-lines and intervals required for our treatment of certain initial-boundary value problems in later sections.

\begin{definition}
\label{Df2.1}
For $s \geq 0$, we write $$H^s({\mathbb R}^+) = \bigl\{f = F\bigl|_{{\mathbb R}^+} : F \in H^s({\mathbb R})\bigr\}$$ and $\anm{f}{s}{{\mathbb R}^+} = \inf \anm{F}{s}{\mathbb R}$.  Similarly, for $s \geq 0$, $$\dot{H}^s({\mathbb R}^+) = \bigl\{f = F\bigl|_{{\mathbb R}^+} : F \in \dot{H}^s({\mathbb R})\bigr\}$$ and $\adnm{f}{s}{{\mathbb R}^+} = \inf \adnm{F}{s}{\mathbb R}$.  For $-\infty < s < +\infty$, $$H^s_0({\mathbb R}^+) = \bigl\{f \in H^s({\mathbb R}) : \supp f \subseteq [0, \infty)\bigr\}$$ and $$\dot{H}^s_0({\mathbb R}^+) = \bigl\{f \in \dot{H}^s({\mathbb R}) : \supp f \subseteq [0, \infty)\bigr\}.$$
\end{definition}

Recall that $H^s({\mathbb R})$ is the set of distributions satisfying $(1 + |\xi|)^s\hat{f}(\xi) \in L_\xi^2$, where $\hat{f}$ denotes the Fourier transform in $x$.  The space $\dot{H}^s({\mathbb R})$ is the homogeneous analogue consisting of distributions satisfying $|\xi|^s\hat{f}(\xi) \in L_\xi^2$.

\begin{remark}
\label{Rk2.1}
For $k = 0, 1, 2, \ldots$, define $$W^{2, k}({\mathbb R}^+) = \bigl\{f : D_x^\alpha%%Or do we want this bizarre $(d/dx)^\alpha$ thing?
f \in L^2({\mathbb R}^+) = H^0({\mathbb R}^+)\text{ for }|\alpha| \leq k\bigr\}.$$  There is a bounded extension operator $E_k : W^{2, k}({\mathbb R}^+) \rightarrow W^{2, k}({\mathbb R})$.  Hence $H^k({\mathbb R}^+) = W^{2, k}({\mathbb R}^+)$.  (In fact, the operator $E_k$ can be defined independently of $k$, so that we have a single operator $E : W^{2, k}({\mathbb R}^+) \rightarrow W^{2, k}({\mathbb R})$ for each $k$.)  Note that this implies that $\anm{f}{k}{{\mathbb R}^+} \sim ||f||_{W^{2, k}({\mathbb R}^+)}$ for $k = 0, 1, 2, \ldots$.  Consequently, the spaces $H^s({\mathbb R}^+)$ interpolate by the complex method $\bigl[L^2({\mathbb R}^+), W^{2, k}({\mathbb R}^+)\bigr]_\theta = H^{\theta k}({\mathbb R}^+)$ for $0 \leq \theta \leq 1$.  By reiteration, similar results hold for $H^s({\mathbb R}^+)$, and we see that $H^s({\mathbb R}^+)$ is a complex interpolation scale for $s \geq 0$.  (See, for instance, \cite{8}.)
\end{remark}

\begin{remark}
\label{Rk2.2}
The space $C_0^\infty({\mathbb R}^+)$ is dense in $H^s_0({\mathbb R}^+)$ for $-\infty < s < +\infty$.
\end{remark}

\begin{definition}
\label{Df2.2}
For $0 < \alpha < +\infty$, $H^{-\alpha}({\mathbb R}^+)$ is the space of linear functionals on $C_0^\infty({\mathbb R}^+)$, with the norm $$||g||_{H^{-\alpha}({\mathbb R}^+)} = \sup \{|g(f)| : f \in C_0^\infty({\mathbb R}^+)\text{ and }||f||_{H^\alpha({\mathbb R})} \leq 1\}.$$  Let ${\cal S}({\mathbb R}^+) = \{f = F\bigl|_{{\mathbb R}^+} : F \in {\cal S}({\mathbb R})\}$ and define, for $s > 0$, $$\dot{H}^{-s}({\mathbb R}^+) = \bigl\{g \in {\cal S}'%%It's "'" here and "*" below.  Which should we use?
({\mathbb R}) : \sup |g(f)| < +\infty\bigr\},$$ where the $\sup$ in the definition above is taken over $f \in C_0^\infty({\mathbb R}^+)$ satisfying $\adnm{f}{s}{{\mathbb R}^+} \leq 1$.
\end{definition}

\begin{proposition}
\label{Pn2.1}
For $\alpha \geq 0$, $H^{-\alpha}({\mathbb R}^+)$ is the dual space of $H^\alpha_0({\mathbb R}^+)$, and $H^{-\alpha}_0({\mathbb R}^+)$ is the dual space of $H^\alpha({\mathbb R}^+)$.  Also, ${\cal S}({\mathbb R}^+)$ is dense in $H^\alpha({\mathbb R}^+)$ for $-\infty < \alpha < +\infty$.
\end{proposition}

\begin{proof}{}
The first assertion follows from Remark \ref{Rk2.2}.  Now suppose that $f \in H^{-\alpha}_0({\mathbb R}^+)$ and define a linear functional on $H^\alpha({\mathbb R}^+)$ by $u \mapsto f(u_1)$, where $u_1 \in H^\alpha({\mathbb R})$ satisfies $u_1\bigl|_{{\mathbb R}^+} = u$.  Such a functional is well-defined because, if $u_2$ is another extension of $u$, then $u_1 - u_2 \in H^\alpha_0({\mathbb R}^-)$, so $u_1 - u_2$ is the limit in the $H^\alpha$ norm of functions in $C_0^\infty({\mathbb R}^-)$ by Remark \ref{Rk2.2}, and $f$ kills such functions.  Thus, we have a natural mapping $H^{-\alpha}_0({\mathbb R}^+) \rightarrow \bigl(H^\alpha({\mathbb R}^+)\bigr)^*$.  This map has an inverse given by restriction, as follows:  If $f \in \bigl(H^\alpha({\mathbb R}^+)\bigr)^*$, then we define $\tilde{f} \in H^{-\alpha}_0({\mathbb R}^+)$ by $\tilde{f}(F) = f\bigl(F\bigl|_{{\mathbb R}^+}\bigr)$.  Clearly, $\tilde{f} \in H^{-\alpha}_0({\mathbb R}^+)$; and $f \mapsto \tilde{f}$ is the inverse of the previous map.

The fact that ${\cal S}({\mathbb R}^+)$ is dense in $H^\alpha({\mathbb R}^+)$ for $\alpha \geq 0$ follows from the density of ${\cal S}({\mathbb R})$ in $H^\alpha({\mathbb R})$.  It is easy to see that $H^\alpha_0({\mathbb R}^+)$ is reflexive, since it is a closed subspace of $H^\alpha({\mathbb R})$ for $-\infty < \alpha < +\infty$.  Hence %%You don't need this.  The dual of any space is reflexive.
$H^{-\alpha}({\mathbb R}^+)$ is reflexive for $\alpha > 0$.  If ${\cal S}({\mathbb R}^+)$ were not dense, by the Hahn-Banach theorem there would exist $f \neq 0$ in $\bigl(H^{-\alpha}({\mathbb R}^+)\bigr)^* = H^\alpha_0({\mathbb R}^+)$ such that $\int f\phi dx = 0$ for all $\phi \in {\cal S}({\mathbb R}^+)$, which is a contradiction.
\end{proof}

\begin{corollary}
\label{Co2.1}
For $\alpha \geq 0$, $H^{-\alpha}_0({\mathbb R}^+)$ is a complex interpolation scale.
\end{corollary}

\begin{proof}{}
Remark \ref{Rk2.1} showed that $H^\alpha({\mathbb R}^+)$ is a complex interpolation scale, and Proposition \ref{Pn2.1} established that $H^{-\alpha}_0({\mathbb R}^+)$ is the dual space of $H^\alpha({\mathbb R}^+)$.  The dual space of a complex interpolation scale is a complex interpolation scale.
\end{proof}

\begin{proposition}
\label{Pn2.2}
For $\alpha \geq 0$, $H^\alpha_0({\mathbb R}^+)$ is a complex interpolation scale.
\end{proposition}

\begin{proof}{}
See Proposition 2.11 in \cite{8}.
\end{proof}

Note that, if $F \in H^\alpha({\mathbb R})$ for some $\alpha > 1/2$, then $F$ is uniformly continuous, and hence $F(0)$ is well-defined.  If $f \in H^\alpha({\mathbb R}^+)$ and $F\bigl|_{{\mathbb R}^+} = f$, then we set $\Tr(f) = F(0)$.  For $f \in {\cal S}({\mathbb R}^+)$, we have $\Tr(f) = f(0)$.

\begin{proposition}
\label{Pn2.3}
If $\alpha > 1/2$, then $\Tr$ is a well-defined bounded linear operator from $H^\alpha({\mathbb R}^+)$ to $\mathbb R$.
\end{proposition}

\begin{proposition}
\label{Pn2.4}
If $1/2 < \alpha < 3/2$, then $$H^\alpha_0({\mathbb R}^+) = \{f \in H^\alpha({\mathbb R}^+) : \Tr(f) = 0\}.$$
\end{proposition}

\begin{proof}{}
Since $C_0^\infty({\mathbb R}^+)$ is dense in $H^\alpha_0({\mathbb R}^+)$, it is clear that, if $f \in H^\alpha_0({\mathbb R}^+)$, then $\Tr(f) = 0$.  For the converse direction, we use two lemmas from \cite{8}:

\begin{lemma}
\label{Lm2.1}
{\bf (3.7 from \cite{8%%It was 8 above -- do we want 1 now?
}.)}  If $1/2 < \alpha < 3/2$, then $$\int_0^\infty |f(x) - f(0)|^2\frac{dx}{x^{2\alpha}} \leq C\anm{f}{\alpha}{\mathbb R}^2.$$
\end{lemma}

\begin{lemma}
\label{Lm2.2}
{\bf (3.8 from \cite{8}.)}  If $1/2 < \alpha < 3/2$, then $$\anm{\chi_{(0, +\infty)}f}{\alpha}{\mathbb R} \leq C\Biggl[\anm{f}{\alpha}{\mathbb R} + \biggl(\int_0^\infty |f(x)|^2\frac{dx}{x^{2\alpha}}\biggr)^{1/2}\Biggr].$$
\end{lemma}

We now complete the proof of the proposition by showing that, if $f \in H^\alpha({\mathbb R}^+)$ satisfies $\Tr(f) = 0$, then $\anm{\chi_{(0, +\infty)}f}{\alpha}{\mathbb R} \leq C\anm{f}{\alpha}{{\mathbb R}^+}$.  This is a direct consequence of the lemmas above.
\end{proof}

\begin{lemma}
\label{Lm2.3}
{\bf (3.5 from \cite{8}.)}  For $0 \leq \alpha < 1/2$, $$\anm{\chi_{(0, +\infty)}f}{\alpha}{\mathbb R} \leq C\anm{f}{\alpha}{\mathbb R}.$$
\end{lemma}

\begin{proposition}
\label{Pn2.5}
Let $\mu \in C^\infty({\mathbb R})$ be such that $\partial_x^j\mu$ is bounded.  Then, for $f \in H^\alpha_0({\mathbb R}^+)$, with $-\infty < \alpha < +\infty$, $$\asnm{\mu f}{\alpha}{{\mathbb R}^+}{0} \leq C\asnm{f}{\alpha}{{\mathbb R}^+}{0}.$$
\end{proposition}

\begin{proof}{}
We need only prove this for $f \in C_0^\infty({\mathbb R}^+)$.  First, consider the case where $0 \leq \alpha < +\infty$.  By complex interpolation, it suffices to treat the case where $\alpha = k$, where $k$ is an integer, which follows by the Leibniz rule since, for $f \in C_0^\infty({\mathbb R}^+)$, $\asnm{f}{k}{{\mathbb R}^+}{0} = ||f||_{W^{2, k}({\mathbb R}^+)}$.  Note that, by taking $f \in {\cal S}({\mathbb R}^+)$, we see that the same result applies to $H^\alpha({\mathbb R}^+)$ for $\alpha \geq 0$ since $\anm{f}{k}{{\mathbb R}^+} \sim ||f||_{W^{2, k}({\mathbb R}^+)}$.  Therefore, to treat the case $\alpha < 0$, we use duality and Proposition \ref{Pn2.1}.
\end{proof}

\begin{proposition}
\label{Pn2.6}
For $0 \leq \alpha < 3/2$ and $\alpha \neq 1/2$, we have that $$\asnm{f}{\alpha}{{\mathbb R}^+}{0} \sim \anm{f}{\alpha}{{\mathbb R}^+}$$ for $f \in H^\alpha_0({\mathbb R}^+)$.
\end{proposition}

\begin{proof}{}
By density of $C_0^\infty({\mathbb R}^+)$ in $H^\alpha_0({\mathbb R}^+)$, it is enough to do this for $f \in C_0^\infty({\mathbb R}^+)$.  Then we clearly have $\anm{f}{\alpha}{{\mathbb R}^+} \leq \asnm{f}{\alpha}{{\mathbb R}^+}{0}$.  For the other inequality, assume first that $0 \leq \alpha < 1/2$.  Let $F \in H^\alpha({\mathbb R})$ satisfy $F\bigl|_{{\mathbb R}^+} = f$ and $\anm{F}{\alpha}{\mathbb R} \leq 2\anm{f}{\alpha}{{\mathbb R}^+}$.  Then, by Lemma \ref{Lm2.3}, $$\anm{f}{\alpha}{\mathbb R} = \anm{\chi_{(0, +\infty)}F}{\alpha}{\mathbb R} \leq C\anm{F}{\alpha}{\mathbb R} \leq 2C\anm{f}{\alpha}{{\mathbb R}^+}.$$  For $1/2 < \alpha \leq 3/2$, a similar proof works by noting that the inequality $$\anm{\chi_{(0, +\infty)}F}{\alpha}{\mathbb R} \leq C\anm{F}{\alpha}{\mathbb R}$$ holds for all $F \in H^\alpha({\mathbb R})$ satisfying $F(0) = 0$.
\end{proof}

\begin{proposition}
\label{Pn2.7}
For $0 < \alpha < 1/2$ and $f \in H^{-\alpha}_0({\mathbb R}^+)$, we have $$\asnm{f}{-\alpha}{{\mathbb R}^+}{0} \sim \anm{f}{-\alpha}{{\mathbb R}^+}.$$
\end{proposition}

\begin{proof}{}
Let $h \in C_0^\infty({\mathbb R}^+)$ satisfy $\anm{h}{\alpha}{\mathbb R} \leq 1$.  Consider now that $$|f(h)| \leq \anm{f}{-\alpha}{\mathbb R}\anm{h}{\alpha}{\mathbb R} \leq \anm{f}{-\alpha}{\mathbb R} = \asnm{f}{-\alpha}{{\mathbb R}^+}{0}.$$  Thus, $\anm{f}{-\alpha}{\mathbb R} \leq \asnm{f}{-\alpha}{{\mathbb R}^+}{0}$.  We need to show that $\asnm{f}{-\alpha}{{\mathbb R}^+}{0} \leq C\anm{f}{-\alpha}{{\mathbb R}^+}$.  Again, it suffices to assume that $f \in C_0^\infty({\mathbb R}^+)$.  Let now $g \in C_0^\infty({\mathbb R})$ satisfy $\anm{g}{\alpha}{\mathbb R} \leq 1$.  Then $$\asnm{f}{-\alpha}{{\mathbb R}^+}{0} = \anm{f}{-\alpha}{\mathbb R} = \sup_g |f(g)| = \sup_g |f(\chi_{(0, +\infty)}g)|.$$  Now, $\anm{\chi_{(0, +\infty)}g}{\alpha}{\mathbb R} \leq C\anm{g}{\alpha}{\mathbb R} \leq 1$, and so $\chi_{(0, +\infty)}g \in H^\alpha_0({\mathbb R}^+)$ satisfies $\asnm{g}{\alpha}{\mathbb R}{0}\linebreak \leq C$.  By Proposition \ref{Pn2.1}, $$\sup_g |f(g)| \leq \anm{f}{-\alpha}{{\mathbb R}^+}\anm{\chi_{(0, +\infty)}g}{\alpha}{\mathbb R},$$ and the proposition follows.
\end{proof}

\begin{proposition}
\label{Pn2.8}
Suppose that $f \in H^\alpha_0({\mathbb R}^+)$ for some $-1/2 < \alpha < 1/2$ and that $\supp f \subseteq [0, 1]$.  Then $$\adsnm{f}{\alpha}{{\mathbb R}^+}{0} \sim \asnm{f}{\alpha}{{\mathbb R}^+}{0}.$$
\end{proposition}

\begin{proof}{}
It is enough to show the result for $f \in C_0^\infty\bigl((0, 1)\bigr)$.  Note that, for $\alpha \geq 0$, $\adnm{f}{\alpha}{\mathbb R} \leq \anm{f}{\alpha}{\mathbb R}$, and, for $\alpha < 0$, $\anm{f}{\alpha}{\mathbb R} \leq \adnm{f}{\alpha}{\mathbb R}$.  Suppose that $0 \leq \alpha < 1/2$.  We need to show that $\anm{f}{\alpha}{\mathbb R} \leq C\adnm{f}{\alpha}{\mathbb R}$.  Since $$\anm{f}{\alpha}{\mathbb R} \sim ||f||_{L^2({\mathbb R})} + ||D^\alpha f||_{L^2({\mathbb R})} = ||f||_{L^2({\mathbb R})} + \adnm{f}{\alpha}{\mathbb R},$$ we only need to show that $||f||_{L^2({\mathbb R})} \leq C||D^\alpha f||_{L^2({\mathbb R})}$.  We know, however, that $f = {\cal I}\comment{Is this what we want?}_\alpha D^\alpha f$, where ${\cal I}_\alpha$ is the standard fractional integral operator; so, by the fractional integration theorem, $$||f||_{L^q({\mathbb R})} \leq C||D^\alpha f||_{L^2({\mathbb R})}, \quad \frac{1}{q} = \frac{1}{2} - \alpha.$$  Since $f$ has compact support, H\"older's inequality gives the desired bound.  Next, we wish to show that $$\adnm{f}{-\alpha}{\mathbb R} \leq C\anm{f}{-\alpha}{\mathbb R}, \quad 0 < \alpha < \frac{1}{2}.$$  Choose $\mu \in C_0^\infty({\mathbb R})$ satisfying $\supp \mu \subseteq (-2, 2)$ and $\mu \equiv 1$ on $[0, 1]$, so that $g = \mu f = f$.  We will show that $|\hat{f}(\tau)| \leq C\anm{f}{-\alpha}{\mathbb R}$ for $|\tau| \leq 1$.  In fact,
\begin{multline*}
|\hat{f}(\tau)| = |\hat{\mu} * \hat{f}(\tau)| = \Bigl|\int \hat{\mu}(\tau - \eta)\hat{f}(\eta)d\eta\Bigr| \\
= \biggl|\int \hat{\mu}(\tau - \eta)(1 + |\eta|)^\alpha\frac{\hat{f}(\eta)}{(1 + |\eta|)^\alpha}d\eta\biggr| \leq C\nm{f}{-\alpha}
\end{multline*}
for $|\tau| \leq 1$, since $\hat{\mu} \in {\cal S}({\mathbb R})$.
\end{proof}

\begin{remark}
\label{Rk2.3}
The conclusion of Proposition \ref{Pn2.8} also holds for $0 \leq \alpha \leq 1$.  We need to show that, for $f \in H^\alpha_0\bigl((0, 1)\bigr)$, $||f||_{L^2} \leq C||D^\alpha f||_{L^2}$ for $1/2 \leq \alpha \leq 1$.  The result clearly holds for $\alpha = 1$.  For $1/2 < \alpha < 1$, note that $|f(x) - f(y)| \leq C|x - y|^{\alpha - 1/2}$, where $C \leq C||D^\alpha f||_{L^2}$%%What in the world does this mean?  Are those $C$s the same?
.  In fact, if $$u(x, t) = \int e^{ix\xi}e^{-t|\xi|}\hat{f}(\xi)d\xi$$ is the harmonic extension of $f$, then $\bigl|\pd{u}{t}(x, t)\bigr| \leq \frac{C}{t^{1/2}t^{1 - \alpha}}$ by a simple use of Cauchy-Schwarz\comment{There is a check in the text; what does it mean?}, and hence the result follows from the proof of Proposition 7, Chapter V in \cite{15}.  (For $\alpha = 1/2$, we can instead obtain a BMO estimate for $f$.)  Taking $y = 0$, we see that $|f(x)| \leq C||D^\alpha f||_{L^2}$, as desired.
\end{remark}

\begin{remark}
\label{Rk2.4}
With $f$ as in Proposition \ref{Pn2.8} and $-1/2 < \alpha \leq 1$ with $\alpha \neq 1/2$, we have $\adsnm{f}{\alpha}{{\mathbb R}^+}{0} \simeq \anm{f}{\alpha}{{\mathbb R}^+}$.  This is because of Propositions \ref{Pn2.6}, \ref{Pn2.7} and \ref{Pn2.8} and Remark \ref{Rk2.3}.
\end{remark}

%%Are we supposed to have anything about properties of the Airy function?

\newsect{The Riemann-Liouville Fractional Integral}{3}

Our method for constructing solutions of certain initial-boundary value problems below exploits properties of a fractional integration operator whose properties are described in this section.

\begin{definition}
\label{Df3.1}
For $h \in C_0^\infty({\mathbb R}^+)$ and $\Re \alpha > 0$, let $$\FracInt_\alpha(h)(t) = \frac{1}{\Gamma(\alpha)}\ic\int_0^t (t - s)^{\alpha - 1}h(s)ds.$$
\end{definition}

The operator $\FracInt_\alpha$ is the {\it Riemann-Liouville fractional integration operator} (of order $\alpha$).  It is shown in \cite{14} that $\FracInt_\alpha(h)(t)$ extends to an analytic function of $\alpha \in {\mathbb C}$ with the following properties:
\begin{gather*}
\FracInt_{-k}(h)(t) = \frac{d^kh}{dt^k}(t), \quad k \in {\mathbb N}, \\
\FracInt_1(h)(t) = \int_0^t h(s)ds, \\
\FracInt_0(h) = h \\
\intertext{ and }
\FracInt_\alpha \FracInt_\beta(h) = \FracInt_{\alpha + \beta}(h).
\end{gather*}

The rest of this section establishes mapping properties of the family $\{\FracInt_\alpha\}$ on the $H^s_0({\mathbb R}^+)$ spaces.

\begin{proposition}
\label{Pn3.1}
The estimate $$||\FracInt_{i\gamma}(h)||_{L^2({\mathbb R}^+)} \leq C||h||_{L^2({\mathbb R}^+)}$$ holds for all $h \in C_0^\infty({\mathbb R}^+)$ and $\gamma \in {\mathbb R}$.
\end{proposition}

We postpone the proof of Proposition \ref{Pn3.1} until later in this section.

\begin{corollary}
\label{Co3.1}
For $0 \leq \Re \alpha \leq 1$, we have that $$||\FracInt_{-\alpha}(h)||_{L^2({\mathbb R}^+)} \leq C\asnm{h}{\alpha}{{\mathbb R}^+}{0}.$$
\end{corollary}

\begin{proof}{}
Take $h \in C_0^\infty({\mathbb R}^+)$.  Since $\FracInt_{-1}(h) = \frac{dh}{dt}$, $$\FracInt_{-1 + i\gamma}(h) = \FracInt_{i\gamma}\FracInt_{-1}(h) = \FracInt_{i\gamma}\biggl(\frac{dh}{dt}\biggr)$$ and, since $\FracInt_0(h) = h$, the result follows from complex interpolation using Propositions \ref{Pn2.2} and \ref{Pn3.1}.
\end{proof}

\begin{proposition}
\label{Pn3.2}
For $0 < r \leq \alpha$ and $0 \leq \alpha \leq 1$, we have $$\asnm{\FracInt_{-\alpha}(h)}{r - \alpha}{{\mathbb R}^+}{0} \leq C\asnm{h}{r}{{\mathbb R}^+}{0}.$$
\end{proposition}

\begin{proof}{}
We will prove the inequality by complex interpolation.  Note that $r - \alpha \leq 0 < r$, so both sides are complex interpolation scales.  Thus, it suffices to prove this when $\Re \alpha = r$ and when $\Re \alpha = 1$.  When $\Re \alpha = r$, this is Corollary \ref{Co3.1}.  When $\Re \alpha = 1$, write $\alpha = 1 - i\gamma$ with $\gamma \in {\mathbb R}$, so that $\FracInt_{-\alpha}(h) = \FracInt_{i\gamma}\bigl(\frac{dh}{dt}\bigr)$.  We need to prove the following facts:

\begin{fact}
\label{Fc3.1}
For $0 \leq r \leq 1$, $\frac{d}{dt} : H^r_0({\mathbb R}^+) \rightarrow H^{r - 1}_0({\mathbb R}^+)$.
\end{fact}

\begin{fact}
\label{Fc3.2}
For $|\alpha| \leq 1$ and $\gamma \in {\mathbb R}$, $\FracInt_{i\gamma} : H^\alpha_0({\mathbb R}^+) \rightarrow H^\alpha_0({\mathbb R}^+)$.
\end{fact}

If Facts \ref{Fc3.1} and \ref{Fc3.2} hold, then the proof is finished.  First we prove Fact \ref{Fc3.1}.  If $r = 1$, then the fact is clear.  Next, we need to check the case $r = 0$.  Note that $r - 1 \leq 0 \leq r$, so both sides are complex interpolation scales.  We need to show that, if $h \in L^2({\mathbb R}^+)$, then $\frac{dh}{dt} \in H^{-1}_0({\mathbb R}^+)$.  By density, it suffices to consider the case $h \in C_0^\infty({\mathbb R}^+)$.  Then, however, $\supp \frac{dh}{dt} \subseteq (0, +\infty)$, so we only need to verify that $\bigl|\bigl|\frac{dh}{dt}\bigr|\bigr|_{H^{-1}({\mathbb R})} \leq C||h||_{L^2({\mathbb R})}$, which is obvious.

We turn our attention to proving Fact \ref{Fc3.2}.  For $\Re \alpha = 0$, this is Proposition \ref{Pn3.1}.  For $\Re \alpha = 1$, we must verify that $\FracInt_{i\gamma} : H^1_0({\mathbb R}^+) \rightarrow H^1_0({\mathbb R}^+)$.  It suffices to consider the case $h \in C_0^\infty({\mathbb R}^+)$.  We have that
\begin{multline*}
\anm{\FracInt_{i\gamma}(h)}{1}{{\mathbb R}^+} \sim ||\FracInt_{i\gamma}(h)||_{L^2({\mathbb R}^+)} + \biggl|\biggl|\frac{d}{dt}\FracInt_{i\gamma}(h)\biggr|\biggr|_{L^2({\mathbb R}^+)} \\
= ||\FracInt_{i\gamma}(h)||_{L^2({\mathbb R}^+)} + \biggl|\biggl|\FracInt_{i\gamma}\frac{dh}{dt}\biggr|\biggr|_{L^2({\mathbb R}^+)} \leq ||h||_{L^2({\mathbb R}^+)} + \biggl|\biggl|\frac{dh}{dt}\biggr|\biggr|_{L^2({\mathbb R}^+)}.
\end{multline*}
Therefore, $\FracInt_{i\gamma}(h) \in H^1({\mathbb R}^+)$.  To show that $\FracInt_{i\gamma}(h) \in H^1_0({\mathbb R}^+)$, in light of Proposition \ref{Pn2.4}, all we need to show is that $\FracInt_{i\gamma}(h)(0) = 0$.

Let $M^\alpha(h)(t) = t^{-\alpha}\FracInt_\alpha(h)(t)$.  Then $$M^\alpha(h)(t) = t^{-\alpha}\frac{1}{\Gamma(\alpha)}\ic\int_0^t (t - s)^{\alpha - 1}h(s)ds = \frac{1}{\Gamma(\alpha)}\ic\int_0^1 (1 - s)^{\alpha - 1}h(st)ds.$$  We break the last expression into two pieces by writing $$M^\alpha(h)(t) = \frac{1}{\Gamma(\alpha)}\ic\int_0^{1/2} (1 - s)^{\alpha - 1}h(st)ds + \frac{1}{\Gamma(\alpha)}\ic\int_{1/2}^1 (1 - s)^{\alpha - 1}h(st)ds =: I + II.$$  The first piece is bounded and the second piece may be rewritten as $$II = \frac{1}{\Gamma(\alpha + 1)}\ic\int_{1/2}^1 (1 - s)^\alpha\frac{d}{ds}h(st)ds + \frac{1}{\Gamma(\alpha + 1)}\biggl(\frac{1}{2}\biggr)^\alpha h\biggl(\frac{1}{2}t\biggr).$$  In this form, it is clear that $$\lim_{\Re \alpha \rightarrow 0^+} \lim_{t \rightarrow 0^+} t^\alpha M^\alpha(h)(t) = 0$$ and so $\FracInt_{i\gamma}(h)(0) = 0$.  Thus $\FracInt_{i\gamma} : H^1_0({\mathbb R}^+) \rightarrow H^1_0({\mathbb R}^+)$ and Fact \ref{Fc3.2} holds for $0 \leq \Re \alpha \leq 1$.

We next need to check that $\FracInt_{i\gamma} : H^{-1}_0({\mathbb R}^+) \rightarrow H^{-1}_0({\mathbb R}^+)$.  Since $H^{-1}_0({\mathbb R}^+)$ is dual to $H^1({\mathbb R}^+)$, we need only verify that, for $f \in H^1({\mathbb R}^+)$ and $h \in H^{-1}_0({\mathbb R}^+)$, $$\frac{1}{\Gamma(i\gamma)}\ic\int_0^\infty \Bigl(\int_0^t (t - s)^{i\gamma - 1}h(s)ds\Bigr)f(t)dt = \frac{1}{\Gamma(i\gamma)}\ic\int_0^\infty h(s)\Bigl(\int_s^\infty f(t)(t - s)^{i\gamma - 1}dt\Bigr)ds$$ is bounded -- i.e., we need to show that, for $f \in H^1({\mathbb R}^+)$, $\Gamma(i\gamma)^{-1}\int_s^\infty f(t)(t - s)^{i\gamma - 1}dt \in H^1({\mathbb R}^+)$.

The $H^1({\mathbb R}^+)$ norm is comparable to the sum of the $L^2({\mathbb R}^+)$ norm and the $\dot{H}^1({\mathbb R}^+)$ norm.  Proposition \ref{Pn3.1} will imply that the $L^2$ contribution is appropriately bounded.  The following integration by parts argument essentially reduces our consideration of the remaining contribution to Proposition \ref{Pn3.1} as well.  Replace $i\gamma$ by an $\alpha \in {\mathbb C}$ with $\Re \alpha > 0$.  The object under examination satisfies $$\frac{1}{\Gamma(\alpha)}\ic\int_s^\infty f(t)(t - s)^{\alpha - 1}dt = \int_s^\infty f(t)\frac{1}{\alpha\Gamma(\alpha)}\frac{d}{dt}(t - s)^\alpha dt.$$  Integrating by parts, we find that this equals $-\bigl(\alpha\Gamma(\alpha)\bigr)^{-1}\int_s^\infty f'(t)(t - s)^\alpha dt$.  The boundary term vanishes since $\Re \alpha > 0$.  Therefore, we observe that $$\frac{d}{dt}\biggl(\frac{1}{\Gamma(\alpha)}\ic\int_s^\infty f(t)(t - s)^{\alpha - 1}dt\biggr) = \frac{1}{\Gamma(\alpha)}\ic\int_s^\infty f'(t)(t - s)^{\alpha - 1}dt,$$ where we have again used the fact that $\Re \alpha > 0$.  Taking the $L_s^2$ norm of both sides, we obtain what we want, provided that
\begin{equation}
\label{3.1}
\Bigl|\Bigl|\lim_{\Re \alpha \rightarrow 0^+} \int_s^\infty g(t)(t - s)^{\alpha - 1}dt\Bigr|\Bigr|_{L_s^2} \leq C||g||_{L_s^2}.
\end{equation}
Note that, by a Fourier transform computation (see 6.12 of Chapter II in \cite{15}), $\Gamma(i\gamma)^{-1}\int_0^\infty f(t)(t - s)^{i\gamma - 1}dt$ is $L^2$-bounded, so \pref{3.1} follows from Proposition \ref{Pn3.1}.
\end{proof}

\begin{proposition}
\label{Pn3.3}
For $0 \leq \theta \leq \alpha$, $0 \leq \alpha \leq 1$ and $\mu \in C_0^\infty({\mathbb R})$, $$\asnm{\mu \FracInt_\alpha(h)}{\alpha - \theta}{{\mathbb R}^+}{0} \leq C\asnm{h}{-\theta}{{\mathbb R}^+}{0}.$$
\end{proposition}

\begin{proof}
Note that $-\theta \leq 0 \leq \alpha - \theta$, so we have complex interpolation scales on both sides of the inequality.  We prove the result in the cases $\Re \alpha = \theta$ and $\Re \alpha = 1$; the result then follows by complex interpolation.

First we consider the case where $\Re \alpha = \theta$.  We wish then to show that $$||\mu \FracInt_{\theta + i\gamma}(h)||_{L^2({\mathbb R}^+)} \leq C\asnm{h}{-\theta}{{\mathbb R}^+}{0}.$$  When $\theta = 0$, this is Proposition \ref{Pn3.1}.  When $\theta = 1$ and $\gamma = 0$, we see by duality that we need to show that $\int_s^\infty \mu(t)f(t)dt \in H^1({\mathbb R}^+)$ for $f \in L^2$.  Since $\mu$ has compact support, the integral is zero for $s$ large, and hence the integral is in $L^2({\mathbb R}^+)$.  Its derivative clearly is in $L^2$, and the proof of this case follows.  When $\gamma \neq 0$ but $\theta = 1$, we have that $\FracInt_{1 + i\gamma}(h) = \FracInt_1\FracInt_{i\gamma}(h)$; by Fact \ref{Fc3.2} above, $\asnm{\FracInt_{i\gamma}(h)}{-1}{{\mathbb R}^+}{0} \sim \asnm{h}{-1}{{\mathbb R}^+}{0}$, so we are done with this case.

Next we consider the case where $\Re \alpha = 1$.  Here we want to verify that $$\asnm{\mu \FracInt_{1 + i\gamma}(h)}{1 - \theta}{{\mathbb R}^+}{0} \leq C\asnm{h}{-\theta}{{\mathbb R}^+}{0}.$$  We checked the case $\theta = 1$ above, so we need only to check the case $\theta = 0$.  By writing $\FracInt_{1 + i\gamma}(h) = \FracInt_1\FracInt_{i\gamma}(h)$, we see that Proposition \ref{Pn3.1} shows that we only need to check $\gamma = 0$.  Thus, we need to check that, if $h \in L^2({\mathbb R}^+)$, then $\mu(t)\int_0^t h(s)ds \in H^1_0({\mathbb R}^+)$.  The quantity $\mu(t)\int_0^t h(s)ds$ is zero at $t = 0$, so, by Proposition \ref{Pn2.4}, it is enough to check that it belongs to $H^1({\mathbb R}^+)$.  The $L^2$ contribution is bounded since $\mu$ has compact support.  The derivative is in $L^2$ since $\mu$ has compact support, $\mu$ and $\mu'$ are bounded and $h \in L^2$.  This completes the proof of Proposition \ref{Pn3.3}, modulo the proof of Proposition \ref{Pn3.1}.
\end{proof}

\begin{proof}[Proof of Proposition \ref{Pn3.1}]
Our task is to estimate $\FracInt_{i\gamma}(h)(t) = \linebreak\Gamma(i\gamma)^{-1}\int_0^t (t - s)^{i\gamma - 1}h(s)ds$ in $L^2({\mathbb R}^+)$.  We interpret this singular expression as the limit as $\alpha \rightarrow 0^+$ of $\FracInt_{\alpha + i\gamma}(h)(t)$.  We change variables by writing $s = tu$ (so $ds = t\,du$) to see that $$\FracInt_{i\gamma}(h)(t) = \frac{t^{i\gamma}}{\Gamma(i\gamma)}\int_0^1 (1 - u)^{i\gamma - 1}h(tu)du.$$  We change variables again by writing $t = e^{-y}$ and $u = e^x$ for $-\infty < x < 0$ (so that $du = e^xdx$) to see that
\begin{align*}
\FracInt_{i\gamma}(h)(e^{-y}) & = \frac{e^{-iy\gamma}}{\Gamma(i\gamma)}\int_{-\infty}^0 (1 - e^x)^{i\gamma - 1}h(e^{x - y})e^xdx \\
                       & = \frac{e^{-iy\gamma}}{\Gamma(i\gamma)}e^{y/2}\ic\int_{-\infty}^0 (1 - e^x)^{i\gamma - 1}e^{x/2}h(e^{x - y})e^{(x - y)/2}dx.
\end{align*}
Suppose that $A$ and $B$ are two functions satisfying $A(x) = B(e^{\pm x})e^{\pm x/2}$.  Then a simple calculation shows that $B \in L^2({\mathbb R}^+)$ if and only if $A \in L^2({\mathbb R})$.  Therefore, $\FracInt_{i\gamma}(h)(t) \in L^2({\mathbb R}_t^+)$ if and only if $\FracInt_{i\gamma}(h)(e^{-y})e^{-y/2} \in L^2({\mathbb R}_y)$, which in turn is true if and only if $$\frac{1}{\Gamma(i\gamma)}\int_{-\infty}^0 (1 - e^x)^{i\gamma - 1}e^{x/2}h(e^{x - y})e^{(x - y)/2}dx \in L^2({\mathbb R}_y).$$  The same simple calculation suggests writing $h(e^x)e^x = f(x)$, so that what we wish to show becomes $$\biggl|\biggl|\frac{1}{\Gamma(i\gamma)}\int_{-\infty}^0 (1 - e^x)^{i\gamma - 1}e^{x/2}f(x - y)dx\biggr|\biggr|_{L^2({\mathbb R}_y)} \leq C||f||_{L^2({\mathbb R})}.$$  Let $$k(x) = \begin{cases}
\Gamma(i\gamma)^{-1}(1 - e^x)^{i\gamma - 1}e^{x/2} & x < 0     \\
0                                                  & x \geq 0.
\end{cases}$$  It suffices to show that $|\hat{k}(\xi)| \leq C$ -- but
\begin{multline*}
\hat{k}(\xi) = \int_{-\infty}^0 \frac{(1 - e^x)^{i\gamma - 1}}{\Gamma(i\gamma)}e^{x/2}e^{ix\xi}dx \\
= \int_{-\infty}^{-1/2} \frac{(1 - e^x)^{i\gamma - 1}}{\Gamma(i\gamma)}e^{x/2}e^{ix\xi}dx + \int_{-1/2}^0 \frac{(1 - e^x)^{i\gamma - 1}}{\Gamma(i\gamma)}e^{x/2}e^{ix\xi}dx \\
=: I(\xi) + II(\xi).
\end{multline*}
The integrand in $I(\xi)$ is integrable, so this contribution is bounded.  The remaining piece is rewritten as $$II(\xi) = \int_{-1/2}^0 e^{ix\xi}\biggl\{\frac{(1 - e^x)^{i\gamma - 1}}{\Gamma(i\gamma)}\biggr\}(e^{x/2} - 1)dx + \int_{-1/2}^0 e^{ix\xi}\frac{(1 - e^x)^{i\gamma - 1}}{\Gamma(i\gamma)}dx.$$  Since $|e^{x/2} - 1| \leq C|x|$ and $|(1 - e^x)^{i\gamma - 1}| \sim |x|^{-1}$, the first piece above has a bounded integrand.  For the second piece of $II(\xi)$, write $$1 - e^x = -x + (1 + x - e^x).$$  Then we can take the Taylor expansion of the function $g(y) = y^{i\gamma - 1}$ and write $$(1 - e^x)^{i\gamma - 1} = (-x)^{i\gamma - 1} + O\biggl(\frac{1}{|x|^2}\biggr)(1 + x - e^x).$$  Since $1 + x - e^x$ vanishes to the second order at $x = 0$, we can ignore the error term.  Finally, $$\biggl|\frac{1}{\Gamma(i\gamma)}\int_{-1/2}^0 e^{ix\xi}(-x)^{i\gamma - 1}dx\biggr| \leq C$$ by a classical calculation (see 6.12 of Chapter II in \cite{15}).  This completes the proof of Proposition \ref{Pn3.1} and our discussion of the properties of the Riemann-Liouville fractional integral.
\end{proof}

\newsect{Some estimates for the group and its associated Duhamel terms, in mixed norm spaces}{4}

\nw{This section begins by recalling the Airy function and its relationship to the linearized KdV equation.  Next, we prove estimates on the three terms arising in the Duhamel representation of solutions of the linearization of \pref{1.2}.  In particular, we establish bounds for the linear solution {\em group}, the {\em Duhamel forcing term} containing the function $g$ and the {\em inhomogeneous Duhamel term}.}

\subsection{The Airy function and KdV}

\begin{definition}
\label{Df4.1}
Define the Airy function $A$ by the property $$\hat{A}(\xi) = e^{i\xi^3}.$$  Similary, define $A_\alpha$ by $$\hat{A}_\alpha(\xi) = |\xi|^\alpha e^{i\xi^3}.$$
\end{definition}

\begin{proposition}
\label{Pn4.1}
The Airy function $A$ is bounded and continuous on $\mathbb R$.  The value $A(0) = C_A \neq 0$.
\end{proposition}

\begin{proof}{}
Lemma 2.7 in \cite{12} shows that $A_\alpha$ is bounded for $0 \leq \alpha \leq \frac{1}{2}$.  This proves that $A$ is bounded, and Sobolev inequalities imply that $A$ is continuous.  Since $$A(x) = \int e^{ix\xi}e^{i\xi^3}d\xi,$$ we know that $$C_A = A(0) = \lim_{\varepsilon \rightarrow 0} \int\blim_{|\xi| < 1/\varepsilon} e^{i\xi^3}d\xi = \lim_{\varepsilon \rightarrow 0} \int\blim_{|\eta| < \varepsilon^{-1/3}} e^{i\eta}\frac{1}{3}\eta^{-2/3}d\eta = c{\cal F}^{-1}(\eta^{-2/3})(1) \neq 0.$$
\end{proof}

\begin{remark}
\label{Rk4.1}
The solution of the initial value problem $$\begin{cases}
\partial_tw + \partial_x^3w = 0 &                    \\
w(x, 0) = \phi(x),              & x \in {\mathbb R},
\end{cases}$$\comment{You want \label{group} here.  Why?} is given as a convolution of the initial data and a time-rescaled Airy function.  Indeed, the solution is
\begin{equation}
\label{4.(0.5)}
w(x, t) = S(t)\phi(x) = \int e^{i(x\xi + t\xi^3)}\hat{\phi}(\xi)d\xi = \int t^{-1/3}A\left(\frac{x - x'}{t^{1/3}}\right)\phi(x')dx'.
\end{equation}
\end{remark}

\subsection{Group estimates}

\begin{lemma}
\label{Lm4.1}
The function $w$ defined above satisfies the following estimates:
\begin{gather}
\label{4.1}\sup_t \anm{w(-, t)}{s}{\mathbb R} \leq C\anm{\phi}{s}{\mathbb R}, \quad -\infty < s < +\infty, \\
\label{4.2}\sup_{\gamma \in {\mathbb R}} ||D_x^{i\gamma}D_x^{s + 1}w||_{L_x^\infty L_t^2} \leq C\adnm{\phi}{s}{\mathbb R}, \quad -\infty < s < +\infty, \\
\label{4.(2.5)}||D_x^s\partial_xw||_{L_x^\infty L_t^2} \leq C\adnm{\phi}{s}{\mathbb R}, \quad -\infty < s < +\infty, \\
\label{4.3}||w||_{L_x^\infty L_T%%Here and several other places the small $t$ suddenly seems to become $T$.  Should this be?
^2} \leq C_T\anm{\phi}{-1}{\mathbb R} \\
\intertext{and}
\label{4.4}||w||_{L_x^\infty\dot{H}^{(s + 1)/3%%Or s + 1/3?
}({\mathbb R}_t)} \leq C\adnm{\phi}{s}{\mathbb R}, \quad -\infty < s < +\infty.
\end{gather}
If $\Psi \in C_0^\infty({\mathbb R})$, then
\begin{equation}
\label{4.(4.5)}
||\Psi(t)w(x, t)||_{C((-\infty, +\infty); H^{(s + 1)/3%%Or s + 1/3?
}({\mathbb R}_t))} \leq C\anm{\phi}{s}{\mathbb R}, \quad -1 \leq s \leq 1
\end{equation}
(here $C$ depends on $\Psi$ and its support).  If $1 \geq s > 1/2$, then $w(0, -)$ is continuous, and $w(0, 0) = \phi(0)$.

Also
\begin{equation}
\label{4.5}
\sup_{\gamma \in {\mathbb R}%%I couldn't read this in the text -- is that what this should be?
}\, e^{-C|\gamma|}||D_x^{i\gamma}D_x^{s - 1/4}w||_{L_x^4L_t^\infty} \leq C\adnm{\phi}{s}{\mathbb R}, \quad -\infty < s < +\infty,
\end{equation}
\begin{equation}
\label{4.(5.25)}
||D_x^sw||_{L_x^5L_t^{10}} \leq C||D_x^s\phi||_{L^2({\mathbb R})},
\end{equation}
\nw{\begin{equation}
\label{4.(5.5)}
||w||_{L_x^2L_T^\infty} \leq C\anm{\phi}{s}{\mathbb R}, \quad s > 3/4,\ \ T < 1
\end{equation}
and
\begin{equation}
\label{4.(5.75)}
||\partial_xw||_{L_t^4L_x^\infty} \leq C\anm{\phi}{3/4}{\mathbb R}.
\end{equation}}
\end{lemma}

\begin{proof}{}
\pref{4.1} follows from the group property of $S$.  \pref{4.2} follows from Theorem 3.5(i) in \cite{10}, and its proof.  To establish \pref{4.3}, write $\phi = \phi_1 + \phi_2$, where $\supp \hat{\phi}_1 \subseteq \{\xi : |\xi| \leq 1\}$.  Note that $\adnm{\phi_2}{-1}{\mathbb R} \leq \anm{\phi}{-1}{\mathbb R}$, and hence the estimate for $\phi_2$ follows from \pref{4.2}.  For the other term, we have
\begin{multline*}
||S(t)\phi_1||_{L_x^\infty L_T^2} \leq C\bigl\{||S(t)\phi_1||_{L_x^2L_T^2} + ||\partial_xS(t)\phi_1||_{L_x^2L_T^2}\bigr\} \\
\leq CT^{1/2}||\phi_1||_{L^2} + CT^{1/2}||\partial_x\phi_1||_{L^2} \leq CT^{1/2}\nm{\phi}{-1}.
\end{multline*}
To establish \pref{4.4}, note that $$w(x, t) = \int e^{ix\xi}e^{it\xi^3}\hat{\phi}(\xi)d\xi = \frac{1}{3}\ic\int e^{ix\eta^{1/3}}e^{it\eta}\hat{\phi}(\eta^{1/3})\frac{d\eta}{\eta^{2/3}},$$ so that, for a fixed $x$, $$\adnm{w(x, -)}{(s + 1)/3%%Or s + 1/3?
}{{\mathbb R}_t}^2 = \frac{1}{9}\ic\int |\eta|^{2(s + 1)/3%%Or 2s + 1/3?
}|\hat{\phi}(\eta^{1/3})|^2\frac{d\eta}{\eta^{4/3}} = \frac{1}{3}\nm{\phi}{s}^2.$$  Next, note that the decomposition $\phi = \phi_1 + \phi_2$, as above, corresponding to $w = w_1 + w_2$, together with the previous estimate, yields $||D_t^{(s + 1)/3%%Or s + 1/3?
}w_2(x, -)||_{L^2} \leq C\nm{\phi}{s}$.  Moreover, $\nm{\phi_1}{-1} \leq \nm{\phi}{s}$, and hence $||w_1(x, -)||_{L_T^2} \leq C\nm{\phi}{s}$.  Also, $\partial_tw_1(x, t) = -S(t)\partial_x^3\phi_1$ and $\nm{\partial_x^3\phi_1}{-1} \leq C\nm{\phi}{s}$, so that $\snm{\Psi w_1(x, -)}{1}{t} \leq C\nm{\phi}{s}$ and hence $\snm{\Psi w_1(x, -)}{(s + 1)/3%%Or s + 1/3?
}{t} \leq C\nm{\phi}{s}$.  Next, if we use the Leibniz rule (Theorem A.12 in \cite{10}) and the bound $|D^\alpha\Psi| \leq \int |\xi|^\alpha|\hat{\Psi}(\xi)|d\xi$, we see that, for a fixed $x$, we have
\begin{multline*}
\bigl|\bigl|D_t^{(s + 1)/3}\bigl(\Psi(t)w_2(x, t)\bigr) - \Psi(t)D_t^{(s + 1)/3%%Or s + 1/3?
}w_2(x, t) - D_t^{(s + 1)/3%%Or s + 1/3?
}\Psi(t)w_2(x, t)\bigr|\bigr|_{L_t^2} \\
\leq C||\Psi||_\infty||D_t^{(s + 1)/3%%Or s + 1/3?
}w_2(x, t)||_{L^2}.
\end{multline*}
Using the bound above for $D_t^{(s + 1)/3%%Or s + 1/3?
}w_2$ and these bounds, we see that $$\bigl|\bigl|D_t^{(s + 1)/3}\bigl(\Psi(t)w_2(x, t)\bigr)\bigr|\bigr|_{L^2} \leq C\nm{\phi}{s}$$ and $||\Psi(t)w_2(x, t)||_{L^2} \leq C\nm{\phi}{s}$.  To obtain the continuity into $H^{(s + 1)/3%%Or s + 1/3?
}$, note that, for $\phi \in H^\infty$, all $t$ and $x$ derivatives of $w$ are continuous; this, together with the previous bound, gives the desired continuity.  The fact that, for $s > 1/2$, $w(0, -)$ is continuous follows from Sobolev embedding and the fact that $(s + 1)/3%%Or s + 1/3?
> 1/2$.  That $w(0, 0) = \phi(0)$ for $s > 1/2$ is an elementary computation.  \pref{4.5} is (3.9) (and its proof) and (3.11)\nw{; \pref{4.(5.5)} is (3.38); and \pref{4.(5.75)} follows from the proof of (3.36), where the latter references are to} \cite{10}.
\end{proof}

\nw{\begin{remark}
\label{Rk4.(1.5)}
The proof of \pref{4.(4.5)} given above uses properties of the cut-off function $\Psi$.  We quantify the dependence of the constant $C$ on $\Psi$.  The proof given above relies on the boundedness of $||\Psi||_{L^\infty}$ and $||D_t^{(s + 1)/3}\Psi||_{L^3}$.  The $L^3$ norm emerges from the last term in the Leibniz rule using the fact that $w_2 \in L^6$ from fractional integration.  By interpolation, $$||D_t^{(s + 1)/3}\Psi||_{L^3} \leq \sup_\gamma ||D^{i\gamma}D^{(s + 1)/3 + 1/6}\Psi||_{L^2}\cdot\sup_\gamma ||D^{i\gamma}\Psi||_{L^\infty}.$$  Since $\frac{s + 1}{3} + \frac{1}{6} < 1$, Sobolev's inequality shows that $$C \leq c\{||\Psi||_{L^2} + ||\Psi'||_{L^2}\}.$$
\end{remark}}

\subsection{Duhamel forcing term estimates}

We now turn our attention to corresponding Duhamel terms of the form
\begin{equation}
\label{4.6}
\int_0^t S(t - t')\delta_0(x)g(t')dt' = w(x, t),
\end{equation}
for suitable $g$.

\begin{lemma}
\label{Lm4.2}
For $h \in C_0^\infty({\mathbb R})$, define $$w(x, t) = \int_0^t S(t - t')\delta_0(x)h(t')dt' = \int_0^t A\biggl(\frac{x}{(t - t')^{1/3}}\biggr)\frac{1}{(t - t')^{1/3}}h(t')dt'.$$  Then $w(x, -)$ and $w(-, t)$ are continuous for each fixed $x$ and $t$.  Moreover, $w$ solves
\begin{equation}
\label{4.7(1)}
\partial_tw + \partial_x^3w = \delta_0(x)h(t)
\end{equation}
in ${\cal S}'({\mathbb R}^2)$, $w(x, 0) = 0$ and $$w(0, t) = \frac{1}{C_A\Gamma(2/3)}\FracInt_{2/3}(h)(t),$$ where $C_A$ is as defined in Proposition \ref{Pn4.1} and $\FracInt_{2/3}$ is as defined in Definition \ref{Df3.1}.
\end{lemma}

\begin{proof}{}
By Proposition \ref{Pn4.1}, $$|w(x, t)| \leq C\ic\int_0^t \frac{1}{(t - t')^{1/3}}|h(t')|dt',$$ which shows that $w$ is well-defined and that $|w(x, 0)| = 0$.  To check the continuity statements, recall from the proof of Proposition \ref{Pn4.1} that $A_\alpha$ is bounded for $0 \leq \alpha \leq 1/2$, and hence that $A$ is H\"older continuous of order $1/2$.  (See, for instance Proposition 7 of Chapter 5 in \cite{15}.) From this and the formula defining $w$, the continuity statements follow easily.  From them and the fact that $C_A \neq 0$, the formula for $w(0, t)$ follows.  Finally, \pref{4.7(1)} follows from the fact that $$\partial_t\Biggl(\frac{1}{t^{1/3}}A\biggl(\frac{x}{t^{1/3}}\biggr)\Biggr) + \partial_x^3\Biggl(\frac{1}{t^{1/3}}A\biggl(\frac{x}{t^{1/3}}\biggr)\Biggr) = \delta_{(0, 0)}(x, t)$$ in ${\cal S}'({\mathbb R}^2)$.
\end{proof}

\begin{lemma}
\label{Lm4.3}
Let $w$ be as in Lemma \ref{Lm4.2}, with $h \in C_0^\infty({\mathbb R}^+%%It was just $\mathbb R$ before -- should it have a $^+$ now?
)$.  Then the following estimates hold:  %%The numbering is not entirely clear here.  There seem to be equation numbers attached to statements; I have tried to figure out a reasonable way to set this up.
\begin{list}{(\thesublemma)}{\usecounter{sublemma}}
\item\label{4.7}%%This and 4.7(1) above were both labelled 4.7 in the text.  I have guessed regarding which reference refers to which.
We have
{\begin{align*}
\sup_x \adnm{w(x, -)}{(s + 1)/3%%Or s + 1/3?
}{{\mathbb R}_t} \leq C\adnm{h}{(s - 1)/3%%Or s - 1/3?
}{\mathbb R}, && -1 \leq s \leq 1, \\
\esup_x\, \sup_{\gamma \in {\mathbb R}}\, e^{-C|\gamma|}||D_x^{i\gamma}D_x^{s + 1}w(x, -)||_{L^2({\mathbb R}_t)} \leq C\adnm{h}{(s - 1)/3%%Or s - 1/3?
}{\mathbb R}, && -1 \leq s \leq 1, \\
\intertext{and}
\esup_x ||D_x^s\partial_xw(x, -)||_{L^2({\mathbb R}_t)} \leq C\adnm{h}{(s - 1)/3%%Or s - 1/3?
}{\mathbb R}, && -1 \leq s < 1.
\end{align*}}
\item\label{4.8}%%Why isn't this equation number printing?
Suppose that $\supp h \subseteq (0, 1)$ and $\Psi \in C_0^\infty\bigl((-2, 2)\bigr)$ is a cut-off function in $t$.  Then
\begin{align*}
\Psi w(x, -) \in C\bigl((-\infty, +\infty); H^{(s + 1)/3%%Or s + 1/3?
}_0({\mathbb R}_t^+)\bigr), && -1 \leq s \leq 1, \\
\Psi w(x, -) \in C\bigl((-\infty, +\infty); H^{(s + 1)/3%%Or s + 1/3?
}({\mathbb R}_t)\bigr), && -1 \leq s \leq 1, \\
\sup_x \asnm{\Psi w(x, -)}{(s + 1)/3%%Or s + 1/3?
}{{\mathbb R}_t^+}{0} \leq C\asnm{h}{(s - 1)/3%%Or s - 1/3?
}{{\mathbb R}^+}{0} && \\
\intertext{and}
\sup_x \anm{\Psi w(x, -)}{(s + 1)/3%%Or s + 1/3?
}{\mathbb R} \leq C\anm{h}{(s - 1)/3%%Or s - 1/3?
}{{\mathbb R}^+}. &&
\end{align*}
\item\label{4.9}%%Why does this carry the equation numbers for itself and the previous item?
Suppose that $h$, $\Psi$ are as in \pref{4.8}, and that $-1/2 < s \leq 1$.  Then $$\sup_t \adnm{\Psi(t)w(-, t)}{s}{{\mathbb R}_x} \leq C\dnm{h}{(s - 1)/3%%Or s - 1/3?
}$$ and $\Psi(t)w(-, t) \in C\bigl((-\infty, +\infty); H^s({\mathbb R}_x)\bigr)$, with $$\sup_t \anm{\Psi(t)w(-, t)}{s}{{\mathbb R}_x} \leq C\asnm{h}{(s - 1)/3%%Or s - 1/3?
}{{\mathbb R}^+}{0}.$$
\item\label{4.10} The operation $h \mapsto \Psi(t)w(x, t)$ has a natural extension (defined by density) to $H^{(s - 1)/3%%Or s - 1/3?
}_0({\mathbb R}^+) \cap \{h : \supp h \subseteq [0, 1]\}$.  If $s > 1/2$, then $w(x, 0) = 0$ and $$\Psi(t)w(0, t) = \frac{\Psi(t)}{C_A\Gamma(2/3)}\FracInt_{2/3}(h)(t)$$ pointwise.  For other values of $s$ this is valid in the sense of \pref{4.8} and \pref{4.9}.
\item\label{4.11}  We have
\begin{gather*}
\sup_\gamma e^{-C|\gamma|}||D_x^{i\gamma}D_x^{3/4}w||_{L_x^4L_t^\infty} \leq C||h||_{L_t^2} \\
\sup_\gamma\, e^{-C|\gamma|}||D_x^{i\gamma}w||_{L_x^4L_t^\infty} \leq C\dsnm{h}{-1/4}{t} \\
||D_t^{1/4}w||_{L_x^4L_t^\infty} \leq C||h||_{L_t^2}.
\end{gather*}
For $1/4 \leq s \leq 1$, $$\sup_\gamma\, e^{-C|\gamma|}||D_x^{s - 1/4}D_x^{i\gamma}w||_{L_x^4L_t^\infty} \leq C\dsnm{h}{(s - 1)/3%%Or s - 1/3?
}{t}.$$
\item\label{4.12}  For $0 \leq s \leq 1$, $$||D_x^sw||_{L_x^5L_t^{10}} \leq C\dnm{h}{(s - 1)/3%%Or s - 1/3?
}.$$
%\item\label{4.13}  We have $$||\partial_x^2w||_{L_x^\infty \dot{H}^{-1/3}({\mathbb R}_t)} \leq C\adnm{h}{-1/3}{{\mathbb R}_t}.$$
\end{list}
\end{lemma}

\begin{proof}{}
We start out with \pref{4.7}.  We prove the first statement by interpolation between the cases $s = 1$ and $s = -1$.  Let $I$ be an interval.  We claim that
%\addtocounter{equation}{7} %%To bring the equation numbering into line with the enumerate environment above.
\begin{equation}
\label{4.14}
\biggl|\biggl|\int_0^t S(t - t')\biggl(\frac{\chi_I}{|I|}\biggr)h(t')dt'\biggr|\biggr|_{L_x^\infty\dot{H}^{2/3}_t} \leq C||h||_{L^2}
\end{equation}
with $C$ independent of $I$.  In fact, \pref{4.4} shows that $||D_t^{1/3}S(t)\phi||_{L_x^\infty L_t^2} \leq C||\phi||_{L^2}$.  From this estimate, duality gives $$\Bigl|\Bigl|\int_{-\infty}^{+\infty} S(-t)D_t^{1/3}g(-, t)dt\Bigr|\Bigr|_{L_x^2} \leq C||g||_{L_x^1L_t^2}.$$  In turn, arguing as on page 554 of \cite{10}, this implies that $$\Bigl|\Bigl|D_t^{1/3}\int_{-\infty}^{+\infty} S(t - t')D_{t'}^{1/3}g(-, t')dt'\Bigr|\Bigr|_{L_x^\infty L_t^2} \leq C||g||_{L_x^1L_t^2},$$ which is the same as
\begin{equation}
\label{4.15}
\Bigl|\Bigl|D_t^{2/3}\int_{-\infty}^{+\infty} S(t - t')g(-, t')dt'\Bigr|\Bigr|_{L_x^\infty L_t^2} \leq C||g||_{L_x^1L_t^2}.
\end{equation}
Finally, using Lemma 3.4 in \cite{10}, we have
\begin{multline}
\label{4.16}
\lim_{\varepsilon \rightarrow 0} \int\bic\int\blim_{\varepsilon < |\xi^3 - \tau| < 1/\varepsilon} e^{i(x\xi + t\tau)}\frac{1}{\xi^3 - \tau}\hat{g}(\xi, \tau)d\xi\,d\tau = 2\ic\int_0^t S(t - t')g(-, t')dt' \\
- \int_{-\infty}^{+\infty} S(t - t')g(-, t')dt' + 2\ic\int_{-\infty}^0 S(t - t')g(-, t')dt'.
\end{multline}
Let
\begin{equation}
\label{4.17}
A(g) = \lim_{\varepsilon \rightarrow 0} \int\bic\int_{\varepsilon < |\xi^3 - \tau| < 1/\varepsilon} e^{i(x\xi + t\tau)}\frac{1}{\xi^3 - \tau}\hat{g}(\xi, \tau)d\xi\,d\tau.
\end{equation}
Then (3.26) in \cite{10} gives $$||D_t^{2/3}A(g)||_{L_x^\infty L_t^2} \leq C||g||_{L_x^1L_t^2},$$ whcih, combined with \pref{4.15} and \pref{4.16}, yields \pref{4.14}.  Note that \pref{4.15} and \pref{4.14} also give
\begin{equation}
\label{4.18}
\Bigl|\Bigl|D_t^{2/3}\int_t^\infty S(t - t')g(-, t')dt'\Bigr|\Bigr|_{L_x^\infty L_t^2} \leq C||g||_{L_x^1L_t^2}.
\end{equation}
Now, to establish the first estimate in \pref{4.7} for $s = 1$, fix $x_0 \in {\mathbb R}$ and let $\phi \in C_0^\infty({\mathbb R}) \cap \dot{H}^{-2/3}({\mathbb R}_t)$.  We use duality to obtain
\begin{multline}
\label{4.19}
\int_{-\infty}^{+\infty} \biggl(\int_0^t A\biggl(\frac{x_0}{(t - t')^{1/3}}\biggr)\frac{h(t')}{(t - t')^{1/3}}dt'\biggr)\phi(t)dt \\
= \int_{-\infty}^{+\infty} h(t')\biggl(\int_{t'}^{+\infty} \phi(t)A\biggl(\frac{x_0}{(t - t')^{1/3}}\biggr)\frac{1}{(t - t')^{1/3}}dt\biggr)dt'.
\end{multline}
The boundedness and continuity of $A$, together with the dominated convergence theorem, show that this equals
\begin{multline}
\label{4.20}
\int_{-\infty}^{+\infty} h(t')\lim_{|I| \rightarrow 0} \int_{t'}^{+\infty} \phi(t)S(t - t')\biggl(\frac{1}{|I|}\chi_{I - x_0}\biggr)(x_0)dt\,dt' \\
= \lim_{|I| \rightarrow 0} \int_{-\infty}^{+\infty} h(t')\!\!\int_{t'}^{+\infty} \phi(t)S(t - t')\biggl(\frac{1}{|I|}\chi_{I - x_0}\biggr)(x_0)dt\,dt'.
\end{multline}
Changing the order of integration and applying \pref{4.14} establishes the first estimate in \pref{4.7} for $s = 1$.

For the case $s = -1$, similar arguments reduce matters to proving that
\begin{equation}
\label{4.21}
\biggl|\biggl|\int_0^t S(t - t')\biggl(\frac{1}{|I|}\chi_I\biggr)(x)h(t')dt'\biggr|\biggr|_{L_x^\infty L_t^2} \leq C\adnm{h}{-2/3}{{\mathbb R}_t}.
\end{equation}
To do this, we again proceed by duality to see that we need to show that
\begin{equation}
\label{4.22}
\biggl|\biggl|D_{t'}^{2/3}\dic\int_{t'}^\infty S(t - t')\biggl(\frac{1}{|I|}\chi_I\biggr)(x)f(t)dt\biggr|\biggr|_{L_x^\infty L_{t'}^2} \leq C||f||_{L_t^2}
\end{equation}
for $f \in C_0^\infty({\mathbb R})$.  This follows from \pref{4.18}.  Thus, interpolation now establishes the first estimate in \pref{4.7}.

For the second estimate, note that, even though $D_x^{s + 1}A$ need not be continuous when $s > -1/2$, since we are interested in the $\esup$ over $x$, a duality and limit argument similar to the one given above reduces matters to the estimate
\begin{equation}
\label{4.23}
\sup_{\gamma \in {\mathbb R}} e^{-C|\gamma|}\biggl|\biggl|D_x^{i\gamma}D_x^{s + 1}\dic\int_0^t S(t - t')\biggl(\frac{\chi_I}{|I|}\biggr)(x)h(t')dt'\biggr|\biggr|_{L_x^\infty L_t^2} \leq C\adnm{h}{(s - 1)/3%%Or s - 1/3?
}{{\mathbb R}_t}.
\end{equation}
In order to show \pref{4.23}, we proceed as in the proof of \pref{4.14}.  We first note that $$||D_x^{i\gamma/2}D_x^{(s + 1)/2}D_t^{(1 - s)/6}S(t)\phi||_{L_x^\infty L_t^2} \leq C||\phi||_{L^2}.$$  In fact,
\begin{align*}
D_x^{i\gamma/2}D_x^{(s + 1)/2}D_t^{(1 - s)/6}S(t)\phi & = D_t^{(1 - s)/6}\dic\int e^{ix\xi}e^{it\xi^3}|\xi|^{i\gamma/2}|\xi|^{(s + 1)/2}\hat{\phi}(\xi)d\xi \\
                                                      & = \frac{1}{3}D_t^{(1 - s)/6}%%Or should it be $\frac{1}{3}D_t^{(1 - s)/6}?
\dic\int e^{ix\eta^{1/3}}e^{it\eta}|\eta|^{i\gamma/6}|\eta|^{(s + 1)/6}\hat{\phi}(\eta^{1/3})\frac{d\eta}{\eta^{2/3}} \\
                                                      & = \frac{1}{3}\ic\int e^{ix\eta^{1/3}}e^{it\eta}|\eta|^{i\gamma/6}|\eta|^{(s + 1)/6}|\eta|^{(1 - s)/6}\hat{\phi}(\eta^{1/3})\frac{d\eta}{\eta^{2/3}} \\
                                                      & = \frac{1}{3}\ic\int e^{ix\eta^{1/3}}e^{it\eta}|\eta|^{i\gamma/6}|\eta|^{1/3}\hat{\phi}(\eta^{1/3})\frac{d\eta}{\eta^{2/3}},
\end{align*}
so that, by Plancherel, the $L^2$ norm in $t$, squared, equals $$\frac{1}{9}\ic\int |\hat{\phi}(\eta^{1/3})|^2|\eta|^{2%%Or is it 1?
/3}\frac{d\eta}{\eta^{4/3}} = \frac{1}{3}\ic\int |\phi|^2.$$  Once more, duality and the argument on page 554 of \cite{10} give
\begin{equation}
\label{4.24}
\Bigl|\Bigl|D_x^{i\gamma}D_x^{s + 1}D_t^{(1 - s)/3%%Or 1 - s/3?  (Surely not!)
}\dic\int_{-\infty}^{+\infty} S(t - t')g(-, t')dt'\Bigr|\Bigr|_{L_x^\infty L_t^2} \leq C||g||_{L_x^1L_t^2}.
\end{equation}
From \pref{4.24} and the fact that $\supp h \subseteq (0, \infty)$, we see that, in view of \pref{4.16}, \pref{4.23} will follow from
\begin{equation}
\label{4.25}
||D_x^{i\gamma}D_x^{s + 1}D_t^{(1 - s)/3%%Or 1 - s/3?
}A(g)||_{L_x^\infty L_t^2} \leq C%%Same $C$ as is to come?
e^{C|\gamma|}||g||_{L_x^1L_t^2}.
\end{equation}
If one follows then the proofs of Theorem 3.5, (3.8) and (3.26) in \cite{10}, we see that this reduces to showing that
\begin{equation}
\label{4.26}
\biggl|\int e^{ix\xi}\frac{|\xi|^{i\gamma}|\xi|^{s + 1}}{\lambda - \xi^3}d\xi\biggr|\!\cdot\!|\lambda|^{(1 - s)/3} \leq C%%Same $C$ as is to come?
e^{C|\gamma|},
\end{equation}
where $C$ is independent of $x$, $\lambda$.  Scaling reduces matters to showing that
\begin{equation}
\label{4.27}
\biggl|\int e^{ix\xi}\frac{|\xi|^{i\gamma}|\xi|^{s + 1}}{\xi^3 - 1}d\xi\biggr| \leq C%%Same $C$ as is to come?
e^{C|\gamma|}.
\end{equation}
To establish \pref{4.27}, we write $1 = \phi_1(\xi) + \phi_2(\xi) + \phi_3(\xi)$, where each $\phi_i$ is even and smooth and $\supp \phi_1 \subseteq \{\xi : |\xi| < \frac{1}{2}\}$, $\supp \phi_2 \subseteq \{\xi : \frac{1}{4} < |\xi| < \frac{5}{4}\}$ and $\supp \phi_3 \subseteq \{\xi : |\xi| \geq 1 + \frac{1}{8}\}$.  The piece corresponding to $\phi_1$ is the Fourier transform of a bounded function with compact support, and hence is bounded.  For the second piece, we write $\xi^3 - 1 = (\xi - 1)(\xi^2 + \xi + 1)$ and note that $\xi^2 + \xi + 1$ has no real zeroes.  Thus, the piece corresponding to $\phi_2$ corresponds to $$\int e^{ix\xi}\frac{1}{\xi - 1}\theta(\xi)d\xi,$$ where $\theta \in C_0^\infty\bigl(\{\xi : \frac{1}{4} < |\xi| < \frac{5}{4}\}\bigr)$.  The Fourier transform of $\frac{1}{\xi - 1}$ is $e^{ix}\sgn x$, so this last integral equals $$\int e^{iy}(\sgn y)\hat{\theta}(x - y)dy,$$ which is bounded.  Finally, for the third piece we have
\begin{multline*}
\int e^{ix\xi}\frac{|\xi|^{i\gamma}|\xi|^{s + 1}}{\xi^3 - 1}\phi_3(\xi)d\xi \\
= \int e^{ix\xi}|\xi|^{i\gamma}|\xi|^{s + 1}\biggl\{\frac{1}{\xi^3 - 1} - \frac{1}{\xi^3}\biggr\}\phi_3(\xi)d\xi + \int e^{ix\xi}\frac{|\xi|^{i\gamma}|\xi|^{s + 1}}{\xi^3}\phi_3(\xi)d\xi
\end{multline*}
Since $\bigl|\frac{1}{\xi^3 - 1} - \frac{1}{\xi^3}\bigr| = \bigl|\frac{1}{\xi^3(\xi^3 - 1)}\bigr| \leq \frac{C}{|\xi|^6}$ on the support of $\phi_3$, the first integral is clearly bounded.  So is the second one, as long as $s < 1$.  When $s = 1$, we are left with $\int e^{ix\xi}\frac{|\xi|^{i\gamma}}{\xi}\phi_3(\xi)d\xi$, which is bounded by $C%%Same $C$ as is to come?
e^{C|\gamma|}$.

Finally, note that the third inequality in \pref{4.7} is proved in a similar manner.  The restriction to $s < 1$ comes from the fact that, when $\gamma = 0$, the integral $\int e^{ix\xi}\phi_3(\xi)\frac{d\xi}{|\xi|}$ is not bounded.

To establish \pref{4.8}, note that, since $\supp h \subseteq (0, \infty)$, \pref{4.16} implies that $w = w_1 + w_2$, where $$w_1(x, t) = \frac{1}{2}A(\delta_0 \otimes h)(x, t)$$ and $$w_2(x, t) = \frac{1}{2}\int_{-\infty}^{+\infty} S(t - t')(\delta_0)(x)h(t')dt'.$$  Note also that, in obtaining the bounds in \pref{4.8}, since $(s - 1)/3%%Or s - 1/3?
\leq 0 \leq (s + 1)/3$, for each fixed $x$, the norms on both sides of the inequalities are complex interpolation scales (by Corollary \ref{Co2.1} and Proposition \ref{Pn2.2}).  Thus, it is enough to establish the bounds for $s \neq 1/2$.  Next, note that, since, by Lemma \ref{Lm4.2}, $w(x, 0) = 0$, Proposition \ref{Pn2.6} shows that, for $s \neq 1/2$, we have $$\asnm{\Psi w(x, -)}{(s + 1)/3%%Or s + 1/3?  I'm done asking this.  Future occurrences of this can be found by searching.  If it doesn't end in )/3 I'll put a comment, once.
}{{\mathbb R}_t^+}{0} \approx \anm{\Psi w(x, -)}{(s + 1)/3}{{\mathbb R}_t^+} \leq \anm{\Psi w(x, -)}{(s + 1)/3}{\mathbb R}.$$  Thus, to establish the estimates in \pref{4.8}, it suffices to show that
\begin{equation}
\label{4.28}
\anm{\Psi w_i(x, -)}{(s + 1)/3}{\mathbb R} \leq C\anm{h}{(s - 1)/3}{\mathbb R}, \quad i = 1, 2.
\end{equation}
In order to establish \pref{4.28}, write $h = h_1 + h_2$, where $h_i \in {\cal S}({\mathbb R})$ and $\supp \hat{h}_1 \subseteq \{\tau : |\tau| < 2\}$ and $\supp \hat{h}_2 \subseteq \{\tau : |\tau| > 1\}$.  We first estimate
\begin{multline*}
\Psi(t)A(h_1)(x, t) = \Psi(t)\ic\int\bic\int e^{ix\xi}e^{it\tau}\frac{\hat{h}_1(\tau)}{\tau - \xi^3}d\tau\,d\xi \\
= \Psi(t)\ic\int e^{it\tau}\hat{h}_1(\tau)\biggl(\int e^{ix\xi}\frac{d\xi}{\tau - \xi^3}\biggr)d\tau = \Psi(t)\ic\int e^{it\tau}\frac{\hat{h}_1(\tau)}{\tau^{2/3}}K(x, \tau)d\tau,
\end{multline*}
where $|K(x, \tau)| \leq C$.  (See the proof of (3.26) in \cite{10}.)  Next, we claim that $|\hat{h}(\tau)| \leq C\nm{h}{(s - 1)/3}$ for $|\tau| \leq 2$.  In fact, since $\supp h \subseteq (0, 1)$, we may choose $\mu \in C_0^\infty$ s.t.\ $\mu \equiv 1$ on $(0, 1)$ and $h = \mu h$.  From this we compute that $$\hat{h}(\tau) = \int \hat{\mu}(\tau - \lambda)\hat{h}(\lambda)d\lambda = \int \hat{\mu}(\tau - \lambda)(1 + |\lambda|)^{(1 - s)/3}\frac{\hat{h}(\lambda)}{(1 + |\lambda|)^{(1 - s)/3}}d\lambda,$$ so that $$|\hat{h}(\tau)| \leq \Bigl(\int |\hat{\mu}(\tau - \lambda)|^2(1 + |\lambda|)^{2(1 - s)/3}d\lambda\Bigr)^{1/2}\nm{h}{(s - 1)/3}.$$  Since $\hat{\mu} \in {\cal S}$ and $|\tau| \leq 2$, the claim follows; but then $|\Psi(t)A(h_1)(x, t)| \leq C\Psi(t)\nm{h}{(1 - s)/3}$, and hence $\Psi A(h_1)(x, -) \in L^2({\mathbb R}_t)$, uniformly in $x$.  A similar argument gives us that $\pd{}{t}\bigl(\Psi A(h_1)(x, -)\bigr) \in L^2({\mathbb R}_t)$ uniformly in $x$, and hence we have the estimate $$\nm{\Psi A(h_1)(x, -)}{(1 + s)/3} \leq C\nm{h}{(1 - s)/3}.$$  For the corresponding term in $w_2$, we have
\begin{multline*}
\Psi(t)\ic\int\bic\int e^{ix\xi}e^{i(t - t')\xi^3}h_1(t')dt'\,d\xi = \Psi(t)\ic\int\bic\int e^{ix\xi}e^{it\tau}\hat{h}_1(\tau)\delta_0(\tau + \xi^3)d\tau\,d\xi \\
= \Psi(t)\ic\int e^{it\tau}\hat{h}_1(\tau)\Bigl(\int e^{ix\xi}\delta_0(\tau + \xi^3)d\xi\Bigr)d\tau \\
= \Psi(t)\ic\int e^{it\tau}\hat{h}_1(\tau)\biggl(\int e^{ix\eta^{1/3}}\delta_0(\tau + \eta)\frac{d\eta}{\eta^{2/3}}\biggr)d\tau = \Psi(t)\ic\int e^{it\tau}\frac{\hat{h}_1(\tau)}{\tau^{2/3}}e^{-ix\tau^{1/3}}d\tau,
\end{multline*}
which is handled in the same way.

Finally, we turn to the terms corresponding to $h_2$.  Note that $\dnm{h_2}{(s - 1)/3} \leq \nm{h}{(s - 1)/3}$, and that $\dnm{h_2}{-1} \leq \nm{h}{(s - 1)/3}$.  We then use the proof of the first estimate in \pref{4.7} to obtain the bounds
\begin{gather*}
||D_t^{(s + 1)/3}A(h_2)(x, -)||_{L_t^2} \leq C\nm{h}{(s - 1)/3}, \\
\Bigl|\Bigl|D_t^{(s + 1)/3}\dic\int_{-\infty}^{+\infty} S(t - t')(x)h_2(t')dt'\Bigr|\Bigr|_{L_t^2} \leq C\nm{h}{(s - 1)/3}, \\
||A(h_2)(x, -)||_{L_t^2} \leq C\nm{h}{(s - 1)/3}\\
\intertext{and}
\Bigl|\Bigl|\int_{-\infty}^{+\infty} S(t - t')(x)h_2(t')dt'\Bigr|\Bigr|_{L_t^2} \leq C\nm{h}{(s - 1)/3}.
\end{gather*}
These bounds, combined with the Leibniz rule and the bound $|D^\alpha\Psi| \leq C$, give us \pref{4.28}.  Finally, the continuity statement follows from Lemma \ref{Lm4.2} and the estimates already proved.  

We now turn to \pref{4.9}.  In view of \pref{4.16}, we are again reduced to studying two terms.  For the term corresponding to $$A(\delta_0 \otimes h)(x, t) = \int\bic\int e^{ix\xi}e^{it\lambda}\frac{\hat{h}(\lambda)}{\lambda - \xi^3}d\xi\,d\lambda,$$ we have $$D_x^sA(\delta_0 \otimes h)(x, t) = \int e^{ix\xi}\biggl(|\xi|^s\ic\int e^{it\lambda}\frac{\hat{h}(\lambda)}{\lambda - \xi^3}d\lambda\biggr)d\xi.$$  By Plancherel, we are reduced to bounding (with $\cal H$ the Hilbert transform)
\begin{multline*}
\biggl|\biggl||\xi|^s\int e^{it\lambda}\frac{\hat{h}(\lambda)}{\lambda - \xi^3}d\lambda\biggr|\biggr|_{L_\xi^2} = \biggl|\biggl||\xi|^s{\cal H}(e^{it-}\hat{h})(\xi^3)\biggr|\biggr|_{L_\xi^2} \\
= \biggl(\int |\eta|^{2s/3}|{\cal H}(e^{it-}\hat{h})|^2(\eta)\frac{d\eta}{\eta^{2/3}}\biggr)^{\!1/2}.
\end{multline*}
Now, since $-1 < 2s/3 - 2/3$ (because $s > -1/2$) and $2s/3 - 2/3 \leq 0$ (because $s \leq 1$), the weight $|\eta|^{2s/3 - 2/3}$ is an $A_2$ weight (see, for instance, Chapter V of \cite{16}), and hence the Hilbert transform is bounded in $L^2$ with this weight.  Thus we get that the last expression is bounded by $$C\biggl(\int |\eta|^{2s/3}|e^{it\eta}\hat{h}(\eta)|^2\frac{d\eta}{\eta^{2/3}}\biggr)^{1/2} = \dnm{h}{(s - 1)/3}.$$  To estimate the term $D_x^s\int_{-\infty}^{+\infty} S(t - t')(\delta_0)(x)h(t')dt'$, we are reduced by the group property to estimating
\begin{equation}
\label{4.29}
\Bigl|\Bigl|D_x^s\ic\int_{-\infty}^{+\infty} S(-t')(\delta_0)(x)h(t')dt'\Bigr|\Bigr|_{L^2} \leq C\dnm{h}{(s - 1)/3}.
\end{equation}
We will establish \pref{4.29} by duality.  Let $f \in C_0^\infty({\mathbb R}) \cap L^2({\mathbb R})$, and consider $$\int D_x^s\ic\int_{-\infty}^{+\infty} S(-t')(\delta_0)(x)h(t')dt'\,\overline{f(x)}dx = \int_{-\infty}^{+\infty} \int \delta_0(x)h(t')\overline{D_x^sS(t')f(x)%%Should the $f(x)$ be conjugated also?
}dx\,dt'.$$  Note that $D_x^sS(t')f$ is a continuous function of $x$, since $D_x^sf \in H^\ell$ for all $\ell \geq 0$ (here we use $-1/2 < s$).  Hence, this expression equals $\int_{-\infty}^{+\infty} h(t')\overline{D_x^sS(t')f(x)}\bigr|_{x = 0}dt'$.  We thus need to show that
\begin{equation}
\label{4.30}
||D_x^sS(t)f||_{L_x^\infty\dot{H}^{(1 - s)/3}({\mathbb R}_t)} \leq C||f||_{L^2},
\end{equation}
but this is contained in the inequality after \pref{4.23}.

To establish the second inequality in \pref{4.9}, note first that, since $s > -1/2$, $-1/2 < (s - 1)/3 \leq 0$, and hence, since $\supp h \subseteq (0, 1)$, Remark \ref{Rk2.4} and Proposition \ref{Pn2.7} show that $\adsnm{h}{\alpha}{{\mathbb R}^+}{0} \simeq \anm{h}{\alpha}{{\mathbb R}^+} \simeq \asnm{h}{\alpha}{{\mathbb R}^+}{0}$, where $\alpha = (s - 1)/3$.  Next, we distinguish two cases:  $-1/2 < s \leq 0$, and $0 \leq s \leq 1$.  In the second case, by the first inequality, $$||D_x^sw(-, t)||_{L^2} \leq C\dnm{h}{(s - 1)/3} \simeq C\asnm{h}{(s - 1)/3}{{\mathbb R}^+}{0}.$$  Also $$||w(-, t)||_{L^2} \leq C\dnm{h}{-1/3} \simeq C\anm{h}{-1/3}{{\mathbb R}^+} \leq C\anm{h}{(s - 1)/3}{{\mathbb R}^+} \leq C\asnm{h}{(s - 1)/3}{{\mathbb R}^+}{0},$$ which gives the estimate.  In the first case, we simply note that $$\nm{w(-, t)}{s} \leq \dnm{w(-, t)}{s} \leq C\dnm{h}{(s - 1)/3} \simeq C\snm{h}{(s - 1)/3}{0}.$$  Finally, the continutity statement follows from the estimate, together with Lemma \ref{Lm4.2}.

Next, note that \pref{4.10} is a direct consequence of \pref{4.8}, \pref{4.9} and Lemma \ref{Lm4.2}.

We turn to \pref{4.11}.  The first estimate there is constained in Proposition 3.5(3) of \cite{7} and also in Proposition 2.3 of \cite{13}.  We will present an alternate proof, using the method in \cite{5'}.  The first part of the proof consists in establishing the inequality
\begin{equation}
\label{4.31}
\Bigl|\Bigl|D_x^{i\gamma}D_x^{3/4}\dic\int_{-\infty}^{+\infty} S(t - t')(\delta_0)(x)h(t')dt'\Bigr|\Bigr|_{L_x^4L_t^\infty} \leq C||h||_{L^2}.
\end{equation}
To prove \pref{4.31}, we use the `$T^*T$ method'.  Thus, testing against an $f$ in $L_x^{4/3}L_t^1$, we see that \pref{4.31} follows from
\begin{equation}
\label{4.32}
\Bigl|\Bigl|D_x^{i\gamma}D_x^{3/4}\dic\int_{-\infty}^{+\infty} S(t - t')f(x, t)dt\Bigr|\Bigr|_{L_x^\infty L_{t'}^2} \leq C||f||_{L_x^{4/3}L_t^1}.
\end{equation}
Fix $x_0$, and compute the $L_{t'}^2$ norm squared; we obtain
\begin{equation}
\label{4.33}
\int\bic\int f(x, t)\Bigl(\int\bic\int K_{x_0}(x, y, t, s)\overline{f(y, s)}dy\,ds\Bigr)dx\,dt,
\end{equation}
where $$K_{x_0}(x, y, t, s) = \int D_x^{i\gamma}D_x^{3/4}A(x_0 - x, t - t')D_y^{-i\gamma}D_y^{3/4}A(y - x_0, t' - s)dt'$$ and $A(x, t) = \frac{1}{t^{1/3}}A(x/t^{1/3})$, where $A$ is as defined in Definition \ref{Df4.1}.  We now proceed to compute $K_{x_0}$, using the fact that $A(x, t) = \int e^{ix\xi}e^{it\xi^3}d\xi$.  We obtain
\begin{multline*}
K_{x_0}(x, y, t, s) \\
= \int\bic\int e^{i(x_0 - x)\xi}e^{i(t - t')\xi^3}|\xi|^{i\gamma}|\xi|^{3/4}d\xi\!\cdot\!\!\int e^{i(y - x_0)\eta}e^{i(t' - s)\eta^3}|\eta|^{-i\gamma}|\eta|^{3/4}d\eta\,dt' \\
= \int\bic\int e^{i(x_0 - x)\xi}e^{it\xi^3}e^{-is\eta^3}|\xi|^{i\gamma}|\xi|^{3/4}e^{i(y - x_0)\eta}|\eta|^{-i\gamma}|\eta|^{3/4}\!\cdot\!\!\int e^{it'(\eta^3 - \xi^3)}dt'\,d\xi\,d\eta \\
= \int\bic\int e^{i(x_0 - x)\xi}e^{it\xi^3}e^{-is\eta^3}|\xi|^{i\gamma}|\xi|^{3/4}e^{i(y - x_0)\eta}|\eta|^{-i\gamma}|\eta|^{3/4}\delta_0(\eta^3 - \xi^3)d\eta\,d\xi.
\end{multline*}
Making the change of variables $\eta^3 = \alpha$ %%Why?  The symbol $\alpha$ doesn't appear again.
and computing the integral over $\eta$%%Or $\alpha$?
, we obtain $$K_{x_0}(x, y, t, s) = \frac{1}{3}\ic\int e^{-i(x - y)\xi}e^{i(t - s)\xi^3}\frac{|\xi|^{3/2}}{|\xi|^2}d\xi = \frac{1}{3}\ic\int e^{-i(x - y)\xi}e^{i(t - s)\xi^3}\frac{d\xi}{|\xi|^{1/2}}.$$  A well-known bound (see, for instance, Lemma 3.6 in \cite{10}) now gives
\begin{equation}
\label{4.33(2)}%%This is the second 4.33.  I think there are no further references to it, so I have let all references in the paper point to the first one.
|K_{x_0}(x, y, t, s)| \leq \frac{C}{|x - y|^{1/2}},
\end{equation}
with $C$ independent of $x_0$, $t$ and $s$, which, combined with \pref{4.33} and the theorem on fractional integration, yields \pref{4.32}.  We will now use the proof in \cite{5'} to deduce our estimate from \pref{4.31}.  This is a general procedure.  First, it is easy to see that it suffices to prove $$||\sup_n w(x, t_n)||_{L_x^4} \leq C||h||_{L^2}$$ for any sequence $0 < t_1 < t_2 < \cdots < t_n < t_{n + 1} < \cdots$, with a bound which does not depend on the particular sequence.  Next, note that it suffices to show that, for any measurable function $N : {\mathbb R} \rightarrow {\mathbb N}$ with finite range, the linear operator $T^{(N)}$ defined by $$T^{(N)}(h)(x) = w(x, t_{N(x)})$$ satisfies the estimate
\begin{equation}
\label{4.35}
||T^{(N)}(h)||_{L_x^4} \leq C||h||_{L^2},
\end{equation}
with $C$ independent of $N$ and $\{t_n\}$.  Assume that $||h||_{L^2} = 1$, and define a probability measure $\lambda$ on $(0, \infty)$ by $\lambda(S) = \int_S |h|^2dt$.  The main tool is Lemma 2.1 in \cite{5'}:  Define $Y_n = [0, t_n]$, so that $Y_n \subseteq Y_{n + 1}$.  Then there exists a collection $\{B_j^m\}$ of measurable subsets of $(0, \infty)$, indexed by $m \in \{0, 1, 2, \cdots\}$ and $1 \leq j \leq 2^m$, satisfying:
\begin{enumerate}
\item For each $m$, $\{B_j^m : 1 \leq j \leq 2^m%%Should it really be $m$?
\}$ is a partition of $(0, \infty)$ into disjoint measurable subsets.
\item Each $B_j^m$ is a union of precisely two sets $B_{j_1}^{m + 1}$, $B_{j_2}^{m + 1}$.
\item $\lambda(B_j^m) = 2^{-m}$ for all $m, j$.
\item Each set $Y_n$ may be decomposed, modulo $\lambda$-null sets, as an empty, finite or countably infinite union $$Y_n = \bigcup_{i \geq 1} B_{j_i}^{m_i},$$ with $m_1 < m_2 < \cdots$.  This decomposition may not be unique, but we fix one such decomposition for each $n$.
\end{enumerate}
Next, define $$T^*(h)(x) = \sup_t \Bigl|\int_0^\infty D_x^{i\gamma}D_x^{3/4}S(t - t')(\delta_0)(x)h(t')dt'\Bigr|,$$ and note that, for $h$ satisfying $\supp h \subseteq (0, \infty)$, \pref{4.31} implies that $||T^*(h)||_{L_x^4} \leq C||h||_{L^2}$.  Also note that we can write $$T^{(N)}(h)(x) = \int_0^{t_{N(x)}} D_x^{i\gamma}D_x^{3/4}S(t_{N(x)} - t')(\delta_0)(x)h(t')dt'.$$  Define now $A_n = \{x : N(x) = n\}$.  Let $$R = \{(m, j, n) : B_j^m\text{ is one of the sets in the decomposition of }Y_n\}$$ and $D_j^m = \bigcup\limits_{(m, j, n) \in R} A_n$.  Note that, for a fixed $m$, the sets $D_j^m$ are disjoint.  In fact, if $D_j^m \cap D_i^m \neq \varnothing$, then, since the $A_n$s are disjoint, there exists an $n$ such that $A_n$ intersects, and hence is contained in, both $D_j^m$ and $D_i^m$.  Hence, since $(m, j, n), (m, i, n) \in R$, $B_j^m$ and $B_i^m$ both occur in the decomposition of $Y_n$ -- but then, by construction, $i = j$.  Define now $h_j^m = h\chi_{B_j^m}$, so that $h\chi_{Y_n} = \sum\limits_{(m, j, n) \in R} h_j^m$.  We then have
\begin{multline*}
\int_0^{t_{N(x)}} D_x^{i\gamma}D_x^{3/4}S(t_{N(x)} - t')(\delta_0)(x)h(t')dt' \\
= \sum_n \chi_{A_n}(x)\ic\int_0^{t_n} D_x^{i\gamma}D_x^{3/4}S(t_n - t')(\delta_0)(x)h(t')dt' \\
= \sum_n \chi_{A_n}(x)\ic\int_0^\infty D_x^{i\gamma}D_x^{3/4}S(t_n - t')(\delta_0)(x)\chi_{(0, t_n)}(t')h(t')dt' \\
= \sum_n \sum_{(m, j, n) \in R} \chi_{A_n}(x)\ic\int_0^\infty D_x^{i\gamma}D_x^{3/4}S(t_n - t')(\delta_0)(x)h_j^m(t')dt'
\end{multline*}
-- but then $$|T^{(N)}(h)(x)| \leq \sum_n \sum_{(m, j, n) \in R} \chi_{A_n}(x)T^*(h_j^m)(x) = \sum_m \sum_j \chi_{D_j^m}(x)T^*(h_j^m)(x).$$  Now fix $m$, and compute (recalling that the $D_j^m$ are disjoint) that
\begin{multline*}
\Bigl|\Bigl|\sum_j \chi_{D_j^m}(x)T^*(h_j^m)\Bigr|\Bigr|_{L_x^4}^4 = \sum_j \int_{D_j^m} |T^*(h_j^m)|^4 \\
\leq C\sum_j \Bigl(\int |h_j^m|^2\Bigr)^2 = C\sum_j 2^{-m}\ic\int_{B_j^m} |h|^2 \leq C2^{-m},
\end{multline*}
since the $B_j^m$ are disjoint.  Thus $$||T^{(N)}(h)||_{L^4} \leq \sum_m 2^{-m/4} \leq C,$$ as desired.  Unfortunately, this elegant general method does not seem to apply to the proof of the other two inequalities in \pref{4.11}, due to the presence of $t$ (fractional) derivatives.  We will therefore give a different type of proof.  We use formula \pref{4.16} to reduce matters to the corresponding two estimates.  The first one, for the second inequality in \pref{4.11}, is
\begin{equation}
\label{4.36}
\Bigl|\Bigl|D_x^{i\gamma}\dic\int_{-\infty}^{+\infty} S(t - t')(\delta_0)(x)h(t')dt'\Bigr|\Bigr|_{L_x^4L_t^\infty} \leq C\dnm{h}{-1/4}.
\end{equation}
Since this is a convolution in the $t$ variable and $D_x^{i\gamma}S(t - t')(\delta_0)(x)\linebreak = D_t^{i\gamma/3}S(t - t')(\delta_0)(x)$, the $D_x^{i\gamma}$ can be moved to $h$, and, since the $\dot{H}^{-1/4}$ norm is not changed by the action of $D_t^{i\gamma/3}$, we can take $\gamma = 0$.  Writing $h = D_t^{1/4}g$ for $g \in L^2$, we see that \pref{4.36} becomes
\begin{equation}
\label{4.37}
\Bigl|\Bigl|D_t^{1/4}\dic\int_{-\infty}^{+\infty} S(t - t')(\delta_0)(x)g(t')dt'\Bigr|\Bigr|_{L_x^4L_t^\infty} \leq C||g||_{L^2}.
\end{equation}
Following the argument used in the proof of \pref{4.31}, matters are reduced to estimating the kernel $$H_{x_0}(x, y, t, s) = \int D_t^{1/4}A(x_0 - x, t - t')D_s^{1/4}A(y - x_0, t' - s)dt'.$$  Recall that $$A(x, t) = \int e^{ix\xi}e^{it\xi^3}d\xi = \frac{1}{3}\int e^{ix\eta^{1/3}}e^{it\eta}\frac{d\eta}{\eta^{2/3}},$$ so that $$D_t^{1/4}A(x, t) = \frac{1}{3}\int e^{ix\eta^{1/3}}e^{it\eta}|\eta|^{1/4}\frac{d\eta}{|\eta|^{2/3}} = \frac{1}{3}\int e^{ix\xi}e^{it\xi^3}|\xi|^{3/4}d\xi,$$ and $H_{x_0}$ is just the kernel $K_{x_0}$ studied before.  We next need an estimate for $A(\delta_0 \otimes h)$.  We will first treat the case $\gamma = 0$.  Writing $h = D_t^{1/4}g$ for $g \in L^2$, we see that this boils down to the estimate
\begin{equation}
\label{4.38}
||D_t^{1/4}A(\delta_0 \otimes g)||_{L_x^4L_t^\infty} \leq C||g||_{L^2}.
\end{equation}
%%  What was supposed to be here? 
Now, $$D_t^{1/4}A(\delta_0 \otimes g)(x, t) = \int\bic\int e^{ix\xi}e^{it\tau}\frac{|\tau|^{1/4}}{\tau - \xi^3}\hat{g}(\tau)d\xi\,d\tau.$$  We first compute $\int e^{ix\xi}\frac{d\xi}{\tau - \xi^3}$.  This has been carried out on page 562 of \cite{10}, where it is shown that $$\int e^{ix\xi}\frac{d\xi}{\xi^3 - 1} = \alpha_1e^{-ix}\sgn x + \alpha_2e^{ix/2}e^{-a_1|x|} + \alpha_3e^{ix/2}(\sgn x)e^{-a_1|x|}$$ with $\alpha_i \in {\mathbb C}$ and $a_1 > 0$.  Moreover, since $$\int e^{ix\xi}\frac{d\xi}{\tau - \xi^3} = \int e^{i(x\tau^{1/3})\eta}\frac{d\eta}{1 - \eta^3}\mult\frac{1}{\tau^{2/3}},$$ we see that we have a sum of three terms, which we will label $T_1$, $T_2$ and $T_3$.  We start out with
\begin{multline*}
T_1(g)(x, t) = C_1\ic\int e^{it\tau}|\tau|^{1/4}e^{-ix\tau^{1/3}}\sgn(x\tau^{1/3})\hat{g}(\tau)\frac{d\tau}{\tau^{2/3}} \\
= C_1(\sgn x)\ic\int e^{it\tau}e^{-ix\tau^{1/3}}(\sgn \tau^{1/3})\hat{g}(\tau)\frac{d\tau}{|\tau|^{5/12}}.
\end{multline*}
Since $g \in L^2$, we are left with estimating $$\tilde{T}_1(h)(x, t) = \int e^{it\tau}e^{-ix\tau^{1/3}}\hat{h}(\tau)\frac{d\tau}{|\tau|^{5/12}}$$ in the $L_x^4L_t^\infty$ norm, in terms of the $L^2$ norm of $h$.  This follows from Theorem 2.5 in \cite{12}.

\old{Next, note that, since $\sgn(x\tau^{2/3}) = (\sgn x)(\sgn \tau^{2/3})$, t}The remaining two terms are similar.  We turn to $T_2$:  $$T_2(g)(x, t) = C_2\ic\int e^{it\tau}e^{ix\tau^{1/3}/2%%Is this right?
}e^{-a_1|x|\,|\tau^{1/3}|}\hat{g}(\tau)\frac{d\tau}{|\tau|^{5/12}}.$$  Following the proof of Theorem 2.5 in \cite{12}, we see that, to obtain the desired estimate for $T_2$, we need merely to use the estimate mentioned on page 564 of \cite{10} (see also Corollary 2.9 of \cite{12}).  Let $f \in L_x^{4/3}L_t^1$.  Then $$\int\bic\int \overline{f(x, t)}T_2(g)(x, t)dx\,dt = C_2\ic\int \hat{g}(\tau)\,\overline{\int\bic\int \frac{e^{-it\tau}e^{-ix\tau^{1/3}/2}e^{-a_1|x|\,|\tau|^{1/3}}}{|\tau|^{5/12}}f(x, t)dxdt}\,d\tau.$$  Thus, we need to estimate $$\biggl|\biggl|\int\bic\int e^{-it\tau}e^{-x\tau^{1/3}/2}e^{-a_1|x|\,|\tau|^{1/3}}\frac{f(x, t)}{|\tau|^{5/12}}dx\,dt\biggr|\biggr|_{L_\tau^2} \leq C||f||_{L_x^{4/3}L_t^1}.$$  If we write out the $L^2$ norm, we obtain
\begin{multline*}
\int\int\bic\int e^{-it\tau}e^{-ix\tau^{1/3}/2}e^{-a_1|x|\,|\tau|^{1/3}}\frac{f(x, t)dx\,dt}{|\tau|^{5/12}} \\
\times \int\bic\int e^{-is\tau}e^{iy\tau^{1/3}/2}e^{-a_1|y|\,|\tau|^{1/3}}\frac{\overline{f(y, s)}dy\,ds}{|\tau|^{5/12}}\,d\tau \\
= \int\bic\int f(x, t)\ic\int\bic\int K(x, y, s, t)\overline{f(y, s)}dyds\,dxdt
\end{multline*}
where $$K(x, y, s, t) = \int e^{i\tau(s - t)}e^{i(y - x)\tau^{1/3}/2}e^{-a_1|x|\,|\tau|^{1/3}}e^{-a_1|y|\,|\tau|^{1/3}}\frac{d\tau}{|\tau|^{5/6}}.$$  The estimates mentioned above give the bound $$|K(x, y, s, t)| \leq \frac{C}{|x - y|^{1/2}},$$ from which our estimate follows easily.

Finally, note that, given \pref{4.16}, the third estimate in \pref{4.11} follows from \pref{4.37} and \pref{4.38}.  To conclude the proof of \pref{4.11}, we are left with showing that
\begin{equation}
\label{4.39}
||D_x^{i\gamma}D_t^{1/4}A(\delta_0 \otimes g)||_{L_x^4L_t^\infty} \leq Ce^{C|\gamma|}%%Is this right?
||g||_{L^2}, \quad \gamma \neq 0
\end{equation}
The proof of this is a variant of the proof already given.  We start out by introducing an operator $S_\eta$, given by
\begin{equation}
\label{4.40}
S_\eta g(x, t) = \int_{-\infty}^{+\infty} e^{it\tau}e^{ix\tau^{1/3}\eta}\frac{|\tau|^{1/4}}{\tau^{2/3}}\hat{g}(\tau)d\tau.
\end{equation}
We claim that
\begin{equation}
\label{4.41}
||S_\eta g||_{L_x^4L_t^\infty} \leq \frac{C}{|\eta|^{1/4}}||g||_{L^2}.
\end{equation}
This follows from scaling, and the bound for $\tilde{T}_1$ sketched above.  Now, in order to obtain \pref{4.39}, write $$D_x^{i\gamma}D_t^{1/4}A(\delta_0 \otimes g) = \int e^{it\tau}\hat{g}(\tau)|\tau|^{i\gamma/3}\frac{|\tau|^{1/4}}{\tau^{2/3}}\biggl(\int e^{ix\tau^{1/3}\eta}|\eta|^{i\gamma}\frac{d\eta}{1 - \eta^3}%%Is this right?
\biggr)d\tau.$$  Writing $1 - \eta^3 = (1 - \eta)(1 + \eta + \eta^2)$, we see that $$K(x\tau^{1/3}, \eta) = e^{ix\tau^{1/3}\eta}|\eta|^{i\gamma}\frac{1}{1 - \eta^3}%%Is this right?
= e^{ix\tau^{1/3}\eta}\frac{|\eta|^{i\gamma}}{(1 - \eta)(1 + \eta + \eta^2)}.$$  Let now $\varphi_j \in C^\infty%%$_\cdot$?
({\mathbb R})$ satisfy $1 = \varphi_1(\eta) + \varphi_2(\eta) + \varphi_3(\eta)$ and $\supp \varphi_1 \subseteq \{\eta : |\eta| < \frac{3}{4}\}$, $\supp \varphi_2 \subseteq \{\eta : \frac{1}{2} < \eta < \frac{3}{2}\}$ and $\supp \varphi_3 \subseteq \{\eta : \eta < -\frac{1}{2}\} \cup \{\eta : \eta > \frac{5}{4}\}$, and consider the corresponding kernels $K_j$, so that our operator is $L_1(g) + L_2(g) + L_3(g)$, where
\begin{equation}
\label{4.42}
L_j%%This is $i$ in the text, but that seems to clash with $i = \sqrt{-1}$.
(g)(x, t) = \int e^{it\tau}\hat{g}(\tau)\frac{|\tau|^{i\gamma/3}}{|\tau|^{5/12}}K_j(x\tau^{1/3}, \eta)d\tau\,d\eta.
\end{equation}
Let us consider $$L_1(g) = \int |\eta|^{i\gamma}\frac{\varphi_1(\eta)}{(1 - \eta)(1 + \eta + \eta^2)}S_\eta(D_t^{i\gamma/3}g)(x, t)d\eta.$$  Using the bound \pref{4.41}, the fact that $$||D_t^{i\gamma/3}g||_{L^2} = ||g||_{L^2}$$ and the fact that, on $\supp \varphi_1$, we have that $\bigl|\frac{1}{(1 - \eta)(1 + \eta + \eta^2)}\bigr| \leq C$, we see that $$||L_1(g)||_{L_x^4L_t^\infty} \leq C\ic\int\blim_{|\eta| < 3/4} \frac{1}{|\eta|^{1/4}}||g||_{L^2}d\eta \leq C||g||_{L^2},$$ as desired.  For $L_3$, we have $$L_3(g) = \int |\eta|^{i\gamma}\frac{\varphi_3(\eta)}{(1 - \eta)(1 + \eta + \eta^2)}S_\eta(D_t^{i\gamma/3}g)(x, t)d\eta,$$ and, on $\supp \varphi_3$, $\bigl|\frac{1}{(1 - \eta)(1 + \eta + \eta^2)}\bigr| \leq \frac{C}{(1 + |\eta|)^3}$, thus yielding a similar bound.  We next turn to $L_2$.  We have
\begin{multline*}
L_2(g) = \int |\eta|^{i\gamma}\frac{\varphi_2(\eta)}{(1 - \eta)(1 + \eta + \eta^2)}S_\eta(D_t^{i\gamma/3}g)(x, t)d\eta \\
= \int (|\eta|^{i\gamma} - 1)\frac{\varphi_2(\eta)}{(1 - \eta)(1 + \eta + \eta^2)}S_\eta(D_t^{i\gamma/3}g)(x, t)d\eta \\
+ \int \frac{\varphi_2(\eta)}{(1 - \eta)(1 + \eta + \eta^2)}S_\eta(D_t^{i\gamma/3}g)(x, t)d\eta \\
= \int \frac{|\eta|^{i\gamma} - 1}{(1 - \eta)(1 + \eta + \eta^2)}\varphi_2(\eta)S_\eta(D_t^{i\gamma/3}g)(x, t)d\eta \\
+ \int \frac{\varphi_2(\eta) - 1%%Whence comes the $-1$?
}{(1 - \eta)(1 + \eta + \eta^2)}S_\eta(D_t^{i\gamma/3}g)(x, t)d\eta \\
+ \int \frac{1}{1 - \eta^3}S_\eta(D_t^{i\gamma/3}g)(x, t)d\eta
\end{multline*}
For the first term, note that, on $\supp \varphi_2$, we have that $\bigl|\frac{|\eta|^{i\gamma} - 1}{(1 - \eta)(1 + \eta + \eta^2)}\bigr| \leq C|\gamma|$, and hence the argument used for $L_1$ applies.  For the second one, note that, on $\supp(\varphi_2 - 1)$, $\bigl|\frac{1}{(1 - \eta)(1 + \eta + \eta^2)}\bigr| \leq \frac{C}{(1 + |\eta|)^3}$, and hence the argument we used for $L_3$ applies.  Finally, the third term is identical to the one handled for the case $\gamma = 0$, and thus \pref{4.39} follows.  Finally, to obtain the last estimate in \pref{4.11}, \old{consider the analytic function $D_x^{z - 1/4}D_x^{i\gamma}w$ for $0 \leq \Re z \leq 1$, and test its $L_x^4L_t^\infty$ norm by duality.  Thus, we consider $$F(z) = \int\bic\int D_x^{z - 1/4}D_x^{i\gamma}w(x, t)\overline{f(x, t)}dx\,dt,$$ and }we estimate it by a well-known variant of the three lines theorem (see Lemma 4.2 in Chapter V of \cite{17}, for instance), using the first two bounds.  To establish \pref{4.12}, we use a similar argument, this time using the last estimate in \pref{4.11} and the second estimate in \pref{4.7}.  Thus, we consider the analytic family $$\{D_x^{z(s - 1/4) + (1 - z)(s + 1)}w(x, t)\}.$$  When $\Re z = 0$, we use the estimate in \pref{4.7}, while, when $\Re z = 1$, we use the one in \pref{4.11}.  Since $s = \frac{4}{5}\bigl(s - \frac{1}{4}\bigr) + \frac{1}{5}(s + 1)$, and $\frac{1}{5} = \frac{4}{5}\mult\frac{1}{4} + \frac{1}{5}\mult 0$, $\frac{1}{10} = \frac{4}{5}\mult 0 + \frac{1}{5}\mult\frac{1}{2}$, the estimate follows.
\old{
Finally, we turn to \pref{4.13}.  Recall from \pref{4.16} that we need to handle only two terms, since $\supp h \subseteq (0, \infty)$.  The term $$D_t^{-1/3}\partial_x^2\ic\int_{-\infty}^{+\infty} S(t - t')(\delta_0)(x)h(t')dt' = \int_{-\infty}^{+\infty} D_t^{-1/3}\partial_x^2S(t - t')(\delta_0)(x)h(t')dt'.$$  Moreover,
\begin{multline*}
D_t^{-1/3}\partial_x^2S(t - t')(\delta_0) = D_t^{-1/3}\dic\int e^{ix\xi}e^{it\xi^3}\xi^2d\xi \\
= \frac{1}{3}\ic\int e^{ix\eta^{1/3}}e^{it\eta}\eta^{2/3}\frac{1}{|\eta|^{1/3}}\frac{d\eta}{\eta^{2/3}} = \frac{1}{3}D_xS(t - t')(\delta_0)(x),
\end{multline*}
and so the desired estimate for this term follows from the one for the corresponding term in the case $s = 0$ of the second estimate in \pref{4.7}.  We are left then with $D_t^{-1/3}\partial_x^2A(\delta_0 \otimes h)(x, t)$, $h \in \dot{H}^{-1/3}$, and hence, since $h = D_t^{1/3}g$ with $g \in L^2$, we have $\partial_x^2A(\delta_0 \otimes g)$, which is estimated by the corresponding term in the case $s = 0$ of the second estimate in \pref{4.7}.}
\end{proof}

\subsection{Inhomogeneous Duhamel term estimates}

\begin{lemma}
\label{Lm4.4}
Let $$w(x, t) = \int_0^t S(t - t')h(x, t')dt'.$$  Then the following estimates hold:
\begin{gather}
\label{4.43}||\partial_xw||_{L_t^\infty L_x^2} \leq C||h||_{L_x^1L_t^2}, \\
\label{4.44}\sup_\gamma ||D_x^{i\gamma}\partial_x^2w||_{L_x^\infty L_t^2} \leq Ce^{C|\gamma|}||h||_{L_x^1L_t^2}, \\
\label{4.45}||w||_{L_x^\infty\dot{H}^{2/3}_t} \leq C||h||_{L_x^1L_t^2}, \\
\label{4.46}||D_x^{-1/2}D_x^{i\gamma}w||_{L_x^4L_t^\infty} \leq C||h||_{L_x^{4/3}L_t^1} \\
\intertext{and}
\label{4.47}||w||_{L_x^5L_t^{10}} \leq C||h||_{L_x^{5/4}L_t^{10/9}}.
\end{gather}
\end{lemma}

\begin{proof}{}
\pref{4.43} is (3.7) in \cite{10}.  \pref{4.44} is (3.8) in \cite{10} when $\gamma = 0$.  The general case is contained in the proof of the second inequality in \pref{4.7} (for $s = 1$).  \pref{4.45} follows from the proof of the case $s = 1$ of the first inequality in \pref{4.7}.  \pref{4.46} is (3.10) in \cite{10}, while \pref{4.47} is (3.12) in \cite{10}.
\end{proof}

\begin{lemma}
\label{Lm4.5}
Let $w$ be as in Lemma \ref{Lm4.4}.  Then the following estimates hold (for $h \in C_0^\infty\bigl({\mathbb R} \times [-T, T]\bigr)$:

For $0 \leq s \leq 1$ and $\Psi$ satisfying $\supp \Psi \subseteq [-T, T]$ and $|\Psi'(t)| \leq C/T$, we have (for $T < 1$, with $\beta(s) = \frac{1}{3} - \frac{s}{6}$)
\begin{gather}
\label{4.48}\begin{split}
\asnm{\Psi w(0, -)}{(s + 1)/3}{{\mathbb R}^+}{0} \leq CT^{\beta(s)}||h||_{L_T^2H_x^s} \\
\anm{\Psi w(0, -)}{(s + 1)/3}{\mathbb R} \leq CT^{\beta(s)}||h||_{L_T^2H_x^s}, \\
\end{split} \\
\label{4.49}\adnm{w(0, -)}{2/3}{{\mathbb R}_t} \leq C||h||_{L_x^1L_t^2}, \\
\label{4.50}\adnm{w(0, -)}{1/4}{{\mathbb R}_t} \leq C||h||_{L_x^{4/3}L_t^1%%Or is it $s$?
} \\
\intertext{and}
\label{4.51}\adnm{w(0, -)}{1/3}{{\mathbb R}_t} \leq C||h||_{L_x^{5/4}L_t^{10/9}}.
\end{gather}
\end{lemma}

\begin{proof}{}
For the estimates in \pref{4.48}, first note that, for $h \in C_0^\infty\bigl({\mathbb R} \times [-T, T]\bigr)$, it is easy to see that $w$ satisfies $$\sup_{t \in [-T, T]} \anm{w(-, t)}{s}{{\mathbb R}_x} < +\infty$$ for all $s$, and hence $w \in C^\infty\bigl({\mathbb R} \times (-T, T)\bigr)$ and $w(-, 0) \equiv 0$.  Next, note that, because of this, if the second estimate in \pref{4.48} is established for $s \neq 1/2$, then the first will follow by Propositions \ref{Pn2.4} and \ref{Pn2.6}, and hence, by Corollary \ref{Co2.1}, this will be the case for all $0 \leq s \leq 1$.  Note also that, by Remark \ref{Rk2.3} (applied to the interval $[-T, T]$), it suffices to estimate $\adnm{\Psi w(0, -)}{(s + 1)/3}{\mathbb R}$.  Thus we will show that
\begin{equation}
\label{4.52}
\biggl|\biggl|\Psi(t)\ic\int_0^t S(t - t')h(x, t')dt'\Bigr|_{x = 0}\biggr|\biggr|_{\dot{H}^{(s + 1)/3}_t({\mathbb R})} \leq CT^{\beta(s)}||h||_{L_T^2H_x^s}.
\end{equation}
We first consider the left-hand side of \pref{4.52} for $s = -1$.  By Minkowski's integral inequality, we can bound it by
\begin{equation}
\label{4.53}
\int_{-T}^T ||S(t)S(-t')h(x, t')||_{L_x^\infty L_T^2}dt' \leq CT^{1/2}||h||_{L_T^2H^{-1}({\mathbb R})},
\end{equation}
where the second inequality follows from \pref{4.3}.

Next, consider the left-hand side in \pref{4.52} when $s = 2$.  From the Leibniz rule, we have that
\begin{multline*}
\partial_t\Bigl(\Psi(t)\ic\int_0^t S(t - t')h(x, t')dt'\Bigr) = \partial_t\Psi(t)\ic\int_0^t S(t - t')h(x, t')dt' + \Psi(t)h(x, t) \\
- \Psi(t)\ic\int_0^t \partial_x^3S(t - t')h(x, t')dt',
\end{multline*}
and each of these terms is to be controlled in $L_x^\infty L_T^2$.  For the first term, we have the estimate
\begin{multline*}
\frac{C}{T}\Bigl|\Bigl|\int_0^t S(t - t')h(x, t')dt'\Bigr|\Bigr|_{L_x^\infty L_T^2} \\
\leq \frac{C}{T}\Bigl|\Bigl|\int_0^t S(t - t')h(x, t')dt'\Bigr|\Bigr|_{L_x^2L_T^2} + \frac{C}{T}\Bigl|\Bigl|\partial_x\ic\int_0^t S(t - t')h(x, t')dt'\Bigr|\Bigr|_{L_x^2L_T^2} \\
\leq \frac{C}{T}\mult CT||h||_{L_T^2L_x^2} + \frac{C}{T}\mult CT||\partial_x h||_{L_T^2L_x^2} \leq C||h||_{L_T^2H_x^2}.
\end{multline*}
For the second term we use the bound $$||h||_{L_x^\infty L_T^2} \leq C||h||_{L_T^2H_x^2}$$ to obtain a similar bound.

Finally, for the third term, using Minkowski's integral inequality and \pref{4.2}, we obtain $$\int_{-T}^T ||S(-t')\partial_x^2h(x, t')||_{L_x^2}dt' \leq CT^{1/2}||h||_{L_T^2H_x^2},$$ so that
\begin{equation}
\label{4.54}
\biggl|\biggl|\Psi(t)\ic\int_0^t S(t - t')h(x, t')dt'\Bigr|_{x = 0}\biggr|\biggr|_{\dot{H}^1_t({\mathbb R})} \leq C||h||_{L_T^2H_x^2}.
\end{equation}
Hence, \pref{4.52} follows by interpolation of \pref{4.53} and \pref{4.54}.

In the proofs of \pref{4.49}--\pref{4.51}, we use \pref{4.16} and the observation that $\chi_{(-\infty, 0)}h$ is in the same space as $h$ to reduce matters to the corresponding estimates for $A(h)$ and $\int_{-\infty}^{+\infty} S(t - t')h(x, t')dt'$.  In the case of \pref{4.49}, this is contained in the proof of the first inequality in \pref{4.7}.  In the case of \pref{4.51}, duality reduces matters to the estimates proved in establishing the third inequality in \pref{4.11}.  Finally, \pref{4.50} can be obtained either by interpolation of \pref{4.49} and \pref{4.51} or by duality from the proof of \pref{4.12}.
\end{proof}

\begin{remark}
\label{Rk4.2}
The restriction to $x = 0$ in the right-hand sides of \pref{4.48}--\pref{4.51} is not significant.  In fact, if $X$ is the space of functions in $t$ in which the estimate has been made, the left-hand side can be replaced by $C({\mathbb R}; X)$.  This is because of translation invariance and the first remark in the proof of Lemma \ref{Lm4.5}.
\end{remark}

\nw{\newsect{Some estimates for the group and its associated Duhamel terms in Bourgain's spaces}{5}

This section introduces a modification of the spaces introduced by Bourgain in \cite{5}.  We also establish useful estimates for the linear solution group, the Duhamel forcing term and the inhomogeneous Duhamel term in these spaces.  The results of this section are combined with a bilinear estimate in \ref{7} to prove local well-posedness of \pref{1.2} for data $(\phi, f) \in L_x^2 \times H_t^{1/3}$ in the standard KdV setting, $\ksub = 1$.

\subsection{Bourgain's spaces with a low frequency modification}

\begin{definition}
\label{Df5.1}
Let $f \in {\cal S}'({\mathbb R}^2)$.  We say that $f \in X_b$ if $||f||_{X_b} < \infty$, where $||\!\cdot\!||_{X_b}$ is defined by
\begin{multline}
\label{5.1}
||f||_{X_b} = \Bigl(\int\bic\int (1 + |\lambda - \xi^3|)^{2b}|\hat{f}(\xi, \lambda)|^2d\xi\,d\lambda\Bigr)^{1/2} \\
+ \Bigl(\int\bic\int\blim_{|\xi| < 1} (1 + |\lambda|)^{2\alpha} \hat{f}(\xi, \lambda)|^2d\xi\,d\lambda\Bigr)^{1/2},
\end{multline}
where $\frac{1}{2} < \alpha < \frac{2}{3}$ is fixed, and $0 < b < \frac{1}{2}$.
\end{definition}

\begin{remark}
\label{Rk5.1}
The space $X_b$ also depends upon the parameter $\alpha$, but we have chosen to suppress this dependence in the notation.  This space and the space $Y_b$ below are introduced for studying the case $\ksub = 1$ of \pref{1.2} for data $(\phi, f) \in L_x^2 \times H_t^{1/3}$.  A natural extension of $X_b$ for problems with data in $H_x^s \times H_t^{(s + 1)/3}$ may be defined by including the spatial Sobolev weight $(1 + |\xi|)^{2s}$ in the two integrals in \pref{5.1}.  We shall refer to this extension as $X_{s, b}$.  Similar comments apply to the space we define next.
\end{remark}

\begin{definition}
\label{Df5.2}
Let $f \in {\cal S}'({\mathbb R}^2)$.  We say that $f \in Y_b$ if $||f||_{Y_b} < +\infty$, where $||\!\cdot\!||_{Y_b}$ is defined by
\begin{multline}
\label{5.2}
||f||_{Y_b} = \Bigl(\int\bic\int \frac{|\hat{f}(\xi, \lambda)|^2}{(1 + |\lambda - \xi^3|)^{2b}}d\xi\,d\lambda\Bigr)^{1/2} \\
+ \Bigl(\int\bic\int\blim_{|\xi| < 1} \frac{|\hat{f}(\xi, \lambda)|^2}{(1 + |\lambda|)^{2(1 - \alpha)}}d\xi\,d\lambda\Bigr)^{1/2} \\
+ \biggl(\int \biggl(\int \frac{|\hat{f}(\xi, \lambda)|}{1 + |\lambda - \xi^3|}d\lambda\biggr)^2d\xi\biggr)^{1/2},
\end{multline}
where $0 < b < \frac{1}{2}$.
\end{definition}

We start off with a useful property of the space $X_b$:  stability under multiplication by smooth time cutoffs.

\begin{lemma}
\label{Lm5.1}
Let $\theta \in C_0^\infty({\mathbb R})$ and $f \in X_b$.  Then, for $b \in (0, 1)$\comment{Shouldn't it be $(0, 1/2)$?}, we have $\theta(t)f(x, t) \in X_b$ and
\begin{equation}
\label{5.3}
||\theta f||_{X_b} \leq C||f||_{X_b}.
\end{equation}
\end{lemma}

\begin{proof}
There are two terms in \pref{5.1} that we must control.  We first consider the low frequency term.  Note that $\int\tbic\int\tblim_{|\xi| \leq 1} |\hat{f}(\xi, \lambda)|^2d\xi\,d\lambda \leq ||f||_{X_b}^2$.  Next, $$\int\bic\int\blim_{|\xi| \leq 1} |\lambda|^{2\alpha}|\widehat{\theta f}|^2 = \int\bic\int\blim_{|\xi| \leq 1} \bigl|D_t^\alpha\bigl(\theta\cdot\tilde{f}(\xi, t)\bigr)\bigr|^2d\xi\,dt,$$ where $\tilde{f}$ denotes the spatial Fourier transform of $f$.  We now use $$||D^s\comment{Is this right?}(\theta v\comment{Or $\gamma$?})||_{L^2} \leq C(||v||_{L^2} + ||D^s v||_{L^2})$$ to get
\begin{multline*}
\int\bic\int\blim_{|\xi| \leq 1} |\lambda|^{2\alpha}|\widehat{\theta f}|^2 \leq C\int\bic\int\blim_{|\xi| \leq 1} (|\tilde{f}(\xi, t)|^2 + |D^\alpha\tilde{f}(\xi, t)|)dt\,d\xi \\
\leq C\int\bic\int\blim_{|\xi| \leq 1} (1 + |\lambda|)^{2\alpha}|\hat{f}(\xi, \lambda)|^2d\xi\,d\lambda.
\end{multline*}
We have shown that $$\int\bic\int_{|\xi| \leq 1} (1 + |\lambda|)^{2\alpha}|\widehat{\theta f}(\xi, \lambda)|^2d\xi\,d\lambda \leq C\int\bic\int\blim_{|\xi| \leq 1} (1 + |\lambda|)^{2\alpha}|\hat{f}(\xi, \lambda)|^2d\xi\,d\lambda.$$

It remains to estimate the other term in \pref{5.1}.  Note that $\widehat{\theta f} = \hat{\theta} *_\lambda \hat{f}$ and that we have to show that, $\forall a \in {\mathbb R}$, $$\int |\hat{f} *_\lambda \hat{\theta}|^2(1 + |\lambda - a|)^{2b}d\lambda \leq C\int |\hat{f}|^2(1 + |\lambda - a|)^{2b}d\lambda,$$ uniformly in $a$.  Since $\int |\hat{\theta}|d\lambda \leq C$, we have $$\int |\hat{f} *_\lambda \hat{\theta}|^2 \leq C\int |\hat{f}|^2(1 + |\lambda - a|)^{2b}d\lambda,$$ so we need to handle $$\int |\lambda - a|^{2b}|\hat{f} *_\lambda \hat{\theta}|^2d\lambda = \int |D_t^b(e^{ia-}\theta\tilde{f})|^2dt.$$  The Leibniz rule for fractional derivatives \cite{10} gives $$||D_t^b(e^{ia-}\tilde{f}\theta) - \theta D_t^b(e^{ia-}\tilde{f}) - D_t^b(\theta)e^{ia-}\tilde{f}||_{L^2} \leq C||\theta||_{L^\infty}||D_t^b(e^{ia-}\tilde{f})||_{L^2}.$$  Since $||D_t^b(\theta)||_{L^\infty} \leq C$ and $||\tilde{f}||_{L^2}$ is in the right-hand side of \pref{5.3}, we are done.
\end{proof}

\subsection{Group estimates}

The following result supplements Lemma \ref{Lm4.1}.

\begin{lemma}
\label{Lm5.2}
Let $\theta \in C_0^\infty({\mathbb R})$ be a cutoff function adapted to $[-T, T]$, with $T < 1$.  For $\phi \in L^2({\mathbb R})$, we have
\begin{equation}
\label{5.4}
||\theta(t)S(t)\phi||_{X_b} \leq C\Bigl(\int (1 + |\lambda|)^{2\alpha}|\hat{\theta}(\lambda)|^2d\lambda\Bigr)^{1/2}||\phi||_{L^2}.
\end{equation}
\end{lemma}

\begin{proof}
We write $(\theta(t)S(t)\phi)\widehat{\ \ }(\xi, \lambda) = \hat{\theta}(\lambda - \xi^3)\hat{\phi}(\xi)$, and use this expression in the two terms of \comment{The whole comprises the parts; the parts compose the whole, or rather the whole is composed of the parts.} the $X_b$ norm.  For the low-frequency term, we have $$\int\bic\int\blim_{|\xi| \leq 1} |\lambda|^{2\alpha}|\hat{\theta}(\lambda - \xi^3)|^2|\hat{\phi}(\xi)|^2d\xi\,d\lambda = \int\blim_{|\xi| \leq 1} |\hat{\phi}(\xi)|^2\Bigl(\int |\lambda|^{2\alpha}|\hat{\theta}(\lambda - \xi^3)|^2d\lambda\Bigr)d\xi.$$  We change variables to see that the above expression is equal to
\begin{multline*}
\int\blim_{|\xi| \leq 1} |\hat{\phi}(\xi)|^2\Bigl(\int |\lambda + \xi^3|^{2\alpha}|\hat{\theta}(\lambda)|^2d\lambda\Bigr)d\xi \\
\leq C\int_{|\xi| \leq 1} |\hat{\phi}(\xi)|^2\Bigl(\int (1 + |\lambda|)^{2\alpha}|\hat{\theta}(\lambda)|^2d\lambda\Bigr)d\xi,
\end{multline*}
so the low-frequency term is fine.  The remaining term in \pref{5.1} is
\begin{multline*}
\int\bic\int (1 + |\lambda - \xi^3|)^{2b}|\hat{\theta}(\lambda - \xi^3)|^2|\hat{\phi}(\xi)|^2d\lambda\,d\xi \\
= \int |\hat{\phi}(\xi)|^2\Bigl(\int (1 + |\lambda - \xi^3|)^{2b}|\hat{\theta}(\lambda - \xi^3)|^2d\lambda\Bigr)d\xi,
\end{multline*}
and we are done since $b < \frac{1}{2}$.
\end{proof}

\subsection{Duhamel forcing term estimate}

We supplement Lemmas \ref{Lm4.2} and \ref{Lm4.3} concerning the Duhamel forcing term with the following estimate.

\begin{lemma}
\label{Lm5.3}
Let $g \in H^{-1/3}({\mathbb R})$ and let $\theta \in C_0^\infty({\mathbb R})$ be a cutoff function.  We have
\begin{equation}
\label{5.5}
\Bigl|\Bigl|\theta(t)\ic\int_0^t S(t - t')\delta_0(x)g(t')dt'\Bigr|\Bigr|_{X_b} \leq C\nm{g}{-1/3}.
\end{equation}
\end{lemma}

\begin{proof}
Let $\hat{w}(\xi, \lambda)$ denote the space-time Fourier transform of $\delta_0(x)g(t)$.  The multiplier representation of the operator $S(t)$ (see \pref{1.6}) allows one to show that the object to be estimated may be written as
\begin{equation}
\label{5.6}
u(x, t) = \theta(t)\int\bic\int e^{ix\xi}\frac{e^{it\lambda} - e^{it\xi^3}}{\lambda - \xi^3}\hat{w}(\xi, \lambda)d\xi\,d\lambda.
\end{equation}
We make a useful decomposition of this function.  Let $\psi \in C_0^\infty({\mathbb R})$ satisfy $\psi = 1$ near $x = 0$ and $\supp \psi \subseteq \{x : |x| < 1\}$.  We first break up \pref{5.6} into two pieces, one near $\lambda - \xi^3 = 0$ and the other far from $\lambda - \xi^3 = 0$, by writing $u = u_1 + u_2$, where
\begin{equation}
\label{5.7}
u_1(x, t) = \theta(t)\int\bic\int e^{ix\xi}\frac{e^{it\lambda} - e^{it\xi^3}}{\lambda - \xi^3}\psi(\lambda - \xi^3)\hat{w}(\xi, \lambda)d\xi\,d\lambda
\end{equation}
and
\begin{equation}
\label{5.8}
u_2(x, t) = \theta(t)\int\bic\int e^{ix\xi}\frac{e^{it\lambda} - e^{it\xi^3}}{\lambda - \xi^3}\bigl[1 - \psi(\lambda - \xi^3)\bigr]\hat{w}(\xi, \lambda)d\xi\,d\lambda.
\end{equation}
Next, we take the Taylor expansion of the exponential to observe that
\begin{equation}
\label{5.9}
u_1(x, t) = \theta(t)\sum_{k = 1}^\infty \frac{i^kt^k}{k!}\int\bic\int e^{ix\xi}e^{it\xi^3}\psi(\lambda - \xi^3)(\lambda - \xi^3)^{k - 1}\hat{w}(\xi, \lambda)d\xi\,d\lambda.
\end{equation}
Let $\theta_k(t) = t^k\theta(t)$, so that $||\theta_k||_{L^2} \leq ||\theta||_{L^2}$ for all $k$; then $\theta_k'(t) = kt^{k - 1}\theta(t) + t^k\theta'(t)$, so that $||\theta_k||_{L^2} + ||\theta_k'||_{L^2} \leq C(1 + |k|)$.  Let $\hat{\phi}_k(\xi) = \int \psi(\lambda - \xi^3)(\lambda - \xi^3)^{k - 1}\hat{w}(\xi, \lambda)d\lambda.$  This notation allows us to reexpress $u_1$ as
\begin{equation}
\label{5.10}
u_1(x, t) = \sum_{k = 1}^\infty \frac{i^k}{k!}\theta_k(t)\int e^{ix\xi}e^{it\xi^3}\hat{\phi}_k(\xi)d\xi.
\end{equation}
In light of Lemma \ref{Lm5.2}, to control $u_1$ as claimed, it suffices to verify that
\begin{equation}
\label{5.11}
||\hat{\phi}_k(\xi)||_{L_\xi^2} = \Bigl|\Bigl|\int \psi(\lambda - \xi^3)(\lambda - \xi^3)^{k - 1}\hat{g}(\lambda)d\lambda\Bigr|\Bigr|_{L_\xi^2} \leq C\nm{g}{-1/3},
\end{equation}
where we have used $\hat{w}(\xi, \lambda) = \hat{g}(\lambda)$\comment{ Is this possible?}.  The support property of $\psi$ shows that $|\hat{\phi}_k(\xi)|$ is controlled by $\int\tblim_{|\lambda - \xi^3| \leq 1} \hat{g}(\lambda)d\lambda$.  We break the $L_\xi^2$ norm into a low- and a high-frequency piece.  Consider first
\begin{multline*}
\int\blim_{|\xi| \leq 2} \Bigl(\int\blim_{|\lambda - \xi^3| \leq 1} |\hat{g}(\lambda)|d\lambda\Bigr)^2d\xi = \int\blim_{|\xi| \leq 2} \Bigl(\int\blim_{\substack{|\lambda - \xi^3| \leq 1 \\ |\lambda| \leq 9\comment{9?}}} |\hat{g}(\lambda)|d\lambda\Bigr)^2d\xi \\
\leq \int\blim_{|\xi| \leq 2} \Bigl(\int\blim_{|\lambda| \leq 9} \frac{|\hat{g}(\lambda)|}{(1 + |\lambda|)^{1/3}}(1 + |\lambda|)^{1/3}d\lambda\Bigr)^2d\xi \leq C\nm{g}{-1/3}^2.
\end{multline*}
Next, we consider $$\int\blim_{|\xi| \geq 2} \Bigl(\int\blim_{|\lambda - \xi^3| \leq 1} |\hat{g}(\lambda)|d\lambda\Bigr)^2d\xi = \int\blim_{|\xi| \geq 2} \Bigl(\int\blim_{\substack{|\lambda - \xi^3| \leq 1 \\ |\lambda| \geq 1}} |\hat{g}(\lambda)|d\lambda\Bigr)^2d\xi.$$  Write $\eta = \xi^3$ and change variables to find that the above expression is equal to
\begin{equation}
\label{5.12}
\int\blim_{|\eta| \geq 8} \Bigl(\int\blim_{\substack{|\lambda - \eta| \leq 1 \\ |\lambda| \leq 1\comment{Should this condition be omitted?}}} |\hat{g}(\lambda)|d\lambda\Bigr)^2\frac{d\eta}{\eta^{2/3}}.
\end{equation}
Now, by H\"older, $$\Bigl(\int\blim_{\substack{|\lambda - \eta| \leq 1 \\ |\lambda| \geq 1}} |\hat{g}(\lambda)|d\lambda\Bigr)^2 \leq \biggl\{\int\blim_{|\lambda - \eta| \leq 1} \biggl[\frac{|\hat{g}(\lambda)|}{(1 + |\lambda|)^{1/3}}\biggr]^{3/2}d\lambda\biggr\}^{4/3}\Bigl\{\int\blim_{\substack{|\lambda - \eta| \leq 1 \\ |\lambda| \geq 1}} (1 + |\lambda|)d\lambda\Bigr\}^{2/3}.$$  The last expression on the preceding line is comparable to $(1 + |\eta|)^{2/3}$.  Going back to \pref{5.10} leaves $$\int\blim_{|\eta| \geq 8} \biggl\{\int\blim_{|\lambda - \eta| \leq 1} \biggl[\frac{|\hat{g}(\lambda)|}{(1 + |\lambda|)^{1/3}}\biggr]^{3/2}d\lambda\biggr\}^{4/3}d\eta.$$  Another application of H\"older, using the constraint on the $\lambda$-integration, proves \pref{5.11} and therefore $$||u_1||_{X_b} \leq C\nm{g}{-1/3}.$$

We turn our attention to the term $u_2$.  There are two pieces,
\begin{equation}
\label{5.13}
u_{2, 1}(x, t) = \theta(t)\int\bic\int e^{ix\xi}e^{it\lambda}\frac{1 - \psi(\lambda - \xi^3)}{\lambda - \xi^3}\hat{w}(\xi, \lambda)d\xi\,d\lambda
\end{equation}
and
\begin{equation}
\label{5.14}
u_{2, 2}(x, t) = \theta(t)\int\bic\int e^{ix\xi}e^{it\xi^3}\frac{1 - \psi(\lambda - \xi^3)}{\lambda - \xi^3}\hat{w}(\xi, \lambda)d\xi\,d\lambda,
\end{equation}
satisfying $u_2 = u_{2, 1} - u_{2, 2}$.  The term $u_{2, 2}$ may be handled using Lemma \ref{Lm5.2} once we establish that
\begin{equation}
\label{5.15}
\biggl|\biggl|\int \frac{1 - \psi(\lambda - \xi^3)}{\lambda - \xi^3}\hat{g}(\lambda)d\lambda\biggr|\biggr|_{L_\xi^2} \leq C\nm{g}{-1/3}.
\end{equation}
The (square of the) $|\xi| \leq 1$ contribution to the $L_\xi^2$ norm appearing in \pref{5.15} is controlled by
\begin{multline*}
  \int\blim_{|\xi| \leq 1} \biggl(\int\blim_{|\lambda - \xi^3| \geq 1/2} \frac{1}{1 + |\lambda - \xi^3|}|\hat{g}(\lambda)|d\lambda\biggr)^2d\xi \\
\leq \int\blim_{|\xi| \leq 1} \Bigl(\int\blim_{|\lambda| \leq 2} |\hat{g}(\lambda)|d\lambda\Bigr)^2d\xi + \int\blim_{|\xi| \leq 1} \Bigl(\int\blim_{|\lambda| \geq 2} \frac{1}{1 + |\lambda|}|\hat{g}(\lambda)|d\lambda\Bigr)^2d\xi \\
\leq C\nm{g}{-1/3}^2.
\end{multline*}
When $|\xi| \geq 1$ and $|\lambda| \leq \frac{1}{2}$, $|\lambda - \xi^3| \sim |\xi^3|$ and we get $$\int\blim_{|\xi| \geq 1} \frac{1}{|\xi|^6}\Bigl(\int\blim_{|\lambda| \leq 1/2} |\hat{g}(\lambda)|d\lambda\Bigr)^2d\xi \leq \nm{g}{-1/3}^2.$$  When $|\xi| \geq 1$ and $|\lambda| \geq \frac{1}{2}$, we are left with
\begin{multline*}
\int\blim_{|\xi| \geq 1} \biggl(\int\blim_{|\lambda| \geq 1/2} \frac{1 - \psi(\lambda - \xi^3)}{\lambda - \xi^3}\hat{g}(\lambda)d\lambda\Bigr)d\xi \\
= \int\blim_{|\eta| \geq 1} \Bigl(\int\blim_{|\lambda| \geq 1/2} \frac{1 - \psi(\lambda - \eta)}{\lambda - \eta}\hat{g}(\lambda)d\lambda\Bigr)^2\frac{d\eta}{\eta^{2/3}}.
\end{multline*}
We now use the fact that $\frac{1}{\eta^{2/3}}$ is an $A_2$ weight \cite{16} to get that the above expression is no greater than $$C\int\blim_{|\lambda| \geq 1/2} \frac{|\hat{g}(\lambda)|^2}{|\lambda|^{2/3}}d\lambda \leq C\nm{g}{-1/3}^2.$$  This completes the treatment of $u_{2, 2}$; all that remains to be considered is $u_{2, 1}$.  Lemma \ref{Lm5.1} shows that we may ignore the time cutoff $\theta(t)$ in \pref{5.13}.  We shall consider the low-frequency term in \pref{5.1} first.  We therefore look at the expression
\begin{equation}
\label{5.16}
\int\bic\int\blim_{|\xi| \leq 1} \frac{|\lambda|^{2\alpha}}{(1 + |\lambda - \xi^3|)^2}|\hat{g}(\lambda)^2d\lambda\,d\xi.
\end{equation}
Note that $\int\tblim_{|\xi| \leq 1} \frac{d\xi}{(1 + |\lambda - \xi^3|)^2} \leq \frac{C}{1 + |\lambda|^2}$ and \pref{5.16} is bounded by $$\int \frac{1}{(1 + |\lambda|)^{2(1 - \alpha)}}|\hat{g}(\lambda)|^2d\lambda,$$ which is fine as long as $1 - \alpha > \frac{1}{3}$, i.e., $\alpha < \frac{2}{3}$, which we have assumed in the definition of $X_b$.  The other part of the $X_b$ norm is controlled as follows.  We look at
\begin{multline}
\label{5.17}
\int\bic\int (1 + |\lambda - \xi^3|)^{2b}\frac{1}{(1 + |\lambda - \xi^3|)^2}|\hat{g}(\lambda)|^2d\lambda\,d\xi \\
= \int |\hat{g}(\lambda)|^2\biggl(\int \frac{d\xi}{(1 + |\lambda - \xi^3|)^{2(1 - b)}}\biggr)d\lambda.
\end{multline}

\begin{claim}
\label{5.18}
Assuming $b < \frac{1}{2}$\comment{ but this was already assumed in the definition!}, we have $$\int \frac{d\xi}{(1 + |\lambda - \xi^3|)^{2(1 - b)}} \leq \frac{C}{(1 + |\lambda|)^{2/3}}.$$
\end{claim}

\begin{proof}[Proof of \pref{5.18}]
Note that our assumption $b < \frac{1}{2}$ guarantees that $2(1 - b) = 1 + \varepsilon$ for some $\varepsilon > 0$.  Therefore, the claim is fine in case $|\lambda| \leq 2$.  For $|\lambda| \geq 2$, we write $\eta = \xi^3$ and change variables to get $$\int \frac{d\xi}{(1 + |\lambda - \xi^3|)^{2(1 - b)}} = \int \frac{d\eta}{|\eta|^{2/3}(1 + |\lambda - \eta|)^{1 + \varepsilon}}.$$  The $\eta$-integral is split into three regions.
\begin{list}{\Roman{enumi}.}{\usecounter{enumi}}
\item $|\eta| < \frac{1}{2}|\lambda|$.

Here we have $(1 + |\lambda - \eta|) \sim (1 + |\lambda|)$, so we write $$\int \frac{d\eta}{|\eta|^{2/3}(1 + |\lambda- \eta|)^{1 + \varepsilon}} \leq \frac{1}{(1 + |\lambda|)^{(2/3) + (\varepsilon/2)}}\int \frac{d\eta}{|\eta|^{2/3}(1 + |\eta|)^{(1/3) + (\varepsilon/2)}}$$ to observe the claim.
\item $|\eta| \sim |\lambda|$.

Here we get $$\frac{1}{|\lambda|^{2/3}}\int\blim_{|\lambda|/2 \leq |\eta| \leq 2|\lambda|} \frac{d\eta}{(1 + |\lambda - \eta|)^{1 + \varepsilon}} \leq \frac{C}{|\lambda|^{2/3}}.$$
\item $2|\lambda| < |\eta|$.

In this region, $(1 + |\lambda - \eta|) \sim (1 + |\eta|)$, so we get $$\int\blim_{|\eta| \geq 2|\lambda|} \frac{d\eta}{|\eta|^{2/3}(1 + |\eta|)^{1 + \varepsilon}} \leq \int\blim_{|\eta| \geq 2|\lambda|} \frac{d\eta}{(1 + |\eta|)^{1 + \varepsilon + (2/3)}} \leq \frac{C}{(1 + |\lambda|)^{(2/3) + \varepsilon}}.$$
\end{list}
This completes the proof of the claim.
\end{proof}

Returning to \pref{5.17} and using the claim finishes off our estimate of $u_{2, 1}$, and with it the proof of Lemma \ref{Lm5.3}.
\end{proof}

\begin{remark}
\label{Rk5.1'}\comment{Notice that this is the second Remark 5.1 in the text.}
The assumption $b < \frac{1}{2}$ is crucial in the proof of the second case of the claim.  It was this that forced us to introduce these modified Bourgain spaces.
\end{remark}

\subsection{Inhomogeneous Duhamel term estimates}

This subsection contains two lemmas concerning the inhomogeneous Duhamel term $\int_0^t S(t - t')w(x, t')dt'$.  We begin by showing that, for $\theta \in C_0^\infty({\mathbb R})$, the formula $w(x, t) \mapsto \theta(t)\int_0^t S(t - t')w(x, t')dt'$ defines a bounded map $Y_b \rightarrow X_b$.  Then we show that a (time-localized) inhomogeneous Duhamel term may be restricted to $\{x = 0\}$ as an $H^{1/3}({\mathbb R}_t)$ function when the inhomogeneity $w \in Y_b$.

\begin{lemma}
\label{Lm5.4}
For $w \in Y_b$ and $\theta \in C_0^\infty$, we have that
\begin{equation}
\label{5.19}
\Bigl|\Bigl|\theta(t)\int_0^t S(t - t')w(x, t')dt'\Bigr|\Bigr|_{X_b} \leq C||w||_{Y_b}.
\end{equation}
\end{lemma}

\begin{proof}
The expression to be controlled in $X_b$ was considered in \pref{5.6} during the proof of Lemma \ref{Lm5.3}.  As in the discussion there, we decompose $u$ into $u_1 + u_{2, 1} -\comment{Or should it really be $+$?} u_{2, 2}$ and estimate the terms separately.  The $u_1$ term was reexpressed in \pref{5.10}, from which we see that, to control it appropriately, it suffices to show that $$||\phi_k||_{L^2} \leq ||w||_{Y_b}.$$  The definition of $\hat{\phi}_k(\xi)$, which appeared between \pref{5.9} and \pref{5.10}, shows that we may consider $$||\phi_k||_{L^2} \leq \Bigl(\int \Bigl(\int\blim_{|\lambda - \xi^3| \leq 1} |\hat{w}(\xi, \lambda)|d\lambda\Bigr)^2d\xi\Bigr)^{1/2},$$ and this is bounded by $$C\biggl(\int \biggl(\int \frac{|\hat{w}(\xi, \lambda)|}{1 + |\lambda - \xi^3|}d\lambda\biggr)^2d\xi\biggr)^{1/2},$$ which is a part of the $Y_b$ norm.

We next focus on the first term in \pref{5.1} for $u_{2, 1}$ and $u_{2, 2}$.

For $u_{2, 2}$, we use Lemma \ref{Lm5.2} and the fact that $$\hat{\phi}(\xi) = \int \frac{1 - \psi(\lambda - \xi^3)}{\lambda - \xi^3}\hat{w}(\xi, \lambda)d\lambda$$ satisfies $$|\hat{\phi}(\xi)| \leq \int \frac{1}{1 + |\lambda - \xi^3|}|\hat{w}(\xi, \lambda)|d\lambda,$$ so that $||\phi||_{L^2} \leq ||w||_{Y_b}$.  For $u_{2, 1}$, we use Lemma \ref{Lm5.1} to ignore the time cutoff $\theta(t)$ in \pref{5.13} and reduce matters to controlling $$\int\bic\int \frac{[1 - \psi(\lambda - \xi^3)]^2}{|\lambda - \xi^3|^2}(1 + |\lambda - \xi^3|)^{2b}|\hat{w}(\xi, \lambda)|^2d\xi\,d\lambda.$$  The support properties of $\psi$ tell us that the above quantity is no greater than
\begin{multline*}
C\int\bic\int \frac{1}{(1 + |\lambda - \xi^3|)^{2(1 - b)}}|\hat{w}(\xi, \lambda)|^2d\xi\,d\lambda \\
= C\int\bic\int \frac{1}{(1 + |\lambda - \xi^3|)^{2(1 - 2b)}}\,\frac{|\hat{w}(\xi, \lambda)|^2}{(1 + |\lambda - \xi^3|)^{2b}}d\xi\,d\lambda.
\end{multline*}
Since $1 - 2b > 0$, $(1 + |\lambda - \xi^3|)^{2(1 - 2b)} \geq 1$, so we bound by $||w||_{Y_b}$\comment{ This sentence does not seem to make much sense}.

It remains to estimate the low-frequency part of \pref{5.1},
\begin{multline*}
\int\blim_{|\xi| \leq 1} \int |\lambda|^{2\alpha}\frac{[1 - \psi(\lambda - \xi^3)]^2}{|\lambda - \xi^3|^2}|\hat{w}(\xi, \lambda)|^2d\xi\,d\lambda \\
\leq C\int\blim_{|\xi| \leq 1} \int \frac{|\lambda|^{2\alpha}}{(1 + |\lambda|)^2}|\hat{w}(\xi, \lambda)|^2d\xi\,d\lambda \\
\leq C\int\blim_{|\xi| \leq 1} \frac{1}{(1 + |\lambda|)^{2(1 - \alpha)}}|\hat{w}(\xi, \lambda)|^2d\xi\,d\lambda \leq C||w||_{Y_b}.
\end{multline*}
\end{proof}

The inhomogeneous Duhamel term defines ``boundary values''.

\begin{lemma}
\label{Lm5.5}
Let $h \in Y_b$ and $\theta \in C_0^\infty({\mathbb R})$.  We have that $$\Bigl|\Bigl|\theta(t)\int_0^t S(t - t')h(x, t')dt'\Bigr|_{\{x = 0\}}\Bigr|\Bigr|_{H_t^{1/3}} \leq C||h||_{Y_b}.$$
\end{lemma}

\begin{proof}
We decompose the object under consideration into $u_1$, $u_{2, 1}$ and $u_{2, 2}$.  As before, $u_1$ and $u_{2, 2}$ rely on the estimate for a term of the form $\theta(t)S(t)\phi$ which we established in Lemma \ref{Lm4.1} (see \pref{4.(4.5)} and Remark \ref{Rk4.(1.5)}).  It remains to estimate $$u_{2, 1}(0, t) = \theta(t)\int\bic\int e^{it\lambda}\hat{h}(\xi, \lambda)\frac{1 - \psi(\lambda - \xi^3)}{\lambda - \xi^3}d\lambda\,d\xi = \theta(t)\beta(t).$$  By Proposition \ref{Pn2.8}, $$\nm{\theta\beta}{1/3} \sim \dnm{\theta\beta}{1/3} = ||D_t^{1/3}(\theta\beta)||_{L^2}.$$  By the Leibniz rule from \cite{10}, $$||D_t^{1/3}(\theta\beta) - \theta D_t^{1/3}(\beta) - D_t^{1/3}(\theta)\beta||_{L^2} \leq C||\theta||_{L^\infty}||D_t^{1/3}(\beta)||_{L^2}.$$  Also, $||\beta||_{L^6} \leq C||D_t^{1/3}(\beta)||_{L^2}$ and $||D_t^{1/3}(\theta)||_{L^3} \leq C$, so it suffices to show that $$||D_t^{1/3}(\beta)||_{L^2} \leq C||w||_{Y_b}.$$  We have $$D_t^{1/3}(\beta) = \int e^{it\lambda}|\lambda|^{1/3}\ic\int \hat{h}(\xi, \lambda)\frac{1 - \psi(\lambda - \xi^3)}{\lambda - \xi^3}d\xi\,d\lambda,$$ so $$||D_t^{1/3}(\beta)||_{L^2} = \int |\lambda|^{2/3}\biggl|\int \hat{h}(\xi, \lambda)\frac{1 - \psi(\lambda - \xi^3)}{\lambda - \xi^3}d\xi\biggr|^2d\lambda.$$  The $\xi$-integral is controlled by
\begin{multline*}
\biggl|\int \hat{h}(\xi, \lambda)\frac{1 - \psi(\lambda - \xi^3)}{\lambda - \xi^3}d\xi\biggr| \leq \int |\hat{h}(\xi, \lambda)|\frac{1}{1 + |\lambda - \xi^3|}d\xi \\
\leq \int \hat{h}(\xi, \lambda)\frac{1}{(1 + |\lambda - \xi^3|)^b}\,\frac{1}{1 + |\lambda - \xi^3|)^{1 - b}}d\xi \\
\leq \biggl(\int \frac{|\hat{h}(\xi, \lambda)|^2}{(1 + |\lambda - \xi^3|)^{2b}}d\xi\biggr)^{1/2}\biggl(\int \frac{1}{(1 + |\lambda - \xi^3|)^{2(1 - b)}}d\xi\biggr)^{1/2}.
\end{multline*}
We apply \pref{5.18}, cancel $|\lambda|^{2/3}$ and bound by $||h||_{Y_b}$ as desired.
\end{proof}
}

\newsect{Existence and uniqueness results for the homogeneous and inhomogeneous linear quarter-plane problems, with data in Sobolev spaces}{6}

We begin by discussing some spaces of functions of space-time which will be used in our study of the nonlinear problems.  Next, we construct and prove estimates on the solutions of the linear homogeneous analogue of \pref{1.2}.  The corresponding construction and estimation of solutions of the linear inhomogeneous problem concludes the section.

Our construction of solutions of \pref{1.1} relies on a contraction mapping argument for solving the forced initial value problem \pref{1.2}.  This procedure requires us to interpret the traces along $\{t = 0\}$ and $\{x = 0\}$ of the solution of \pref{1.2} in order to validate the boundary conditions in \pref{7.1}\comment{really?}.  The time-localized solutions $w$ of \pref{1.2} that we construct will have {\em good $H_x^s \times H_t^{(s + 1)/3}$ traces} in the sense that, for some $T > 0$,
%\vspace{1mm}
\begin{equation*}
\tag{Good Traces}
%\labeqn{(Good Traces)}{
w \in C\bigl((-T, T); H^s({\mathbb R}_x)\bigr) \cap C\bigl((-\infty, \infty); H^{(s + 1)/3}({\mathbb R}_t)\bigr).%}
\end{equation*}
%\vspace{1mm}
The contraction estimate is established in a space closely related to the linearization of \pref{1.2}.  For the general case $k \geq 2$ (with optimizations presented in the modified KdV $k = 2$ and $L^2$ critical $k = 4$ settings), we employ the mixed-norm spaces used in \cite{10}.  These spaces are based on the {\em local smoothing effect} $$||\partial_xS(t)\phi||_{L_x^\infty L_t^2} \leq C||\phi||_{L_x^2}$$ and the {\em maximal function estimate} $$||S(t)\phi||_{L_x^4L_t^\infty} \leq C||D^{1/4}\phi||_{L_x^2}.$$  Recall from the introduction that the local smoothing effect is closely related to the good traces property (Good~Traces), above.  In the standard KdV setting $k = 1$, we establish the contraction estimate in the space $X_b$ defined in \ref{5}.

In this section, we use the results in the previous sections to construct solutions of the homogeneous and inhomogeneous linear analogues of \pref{1.2}.  We also prove a uniqueness result.

\subsection{Homogeneous solution operator}

We consider the linear homogeneous initial-boundary value problem
\begin{equation}
\label{6.1}
\begin{cases}
\partial_tw + \partial_x^3w = 0, & x > 0, t \in (0, T_0) \\
w(x, 0) = \phi(x),               & x > 0                 \\
w(0, t) = f(t),                  & t \in [0, T_0],
\end{cases}
\end{equation}
where $\phi \in H^s({\mathbb R}^+)$, $f \in H^{(s + 1)/3}({\mathbb R}^+)$, $0 \leq s \leq 1$, $w \in C\bigl([0, +\infty); H^{(s + 1)/3}((0, T_0))\bigr) \cap C\bigl([0, T_0]; H^s({\mathbb R}_x^+)\bigr)$.  The equation holds in the sense of ${\cal D}'\bigl({\mathbb R}^+ \times (0, T_0)\bigr)$ and the initial value $\phi$ is taken in the sense of $C\bigl([0, +\infty); H^s({\mathbb R}_x^+)\bigr)$, while the `lateral value' $f$ is taken in the sense of $C\bigl([0, +\infty); H^{(s + 1)/3}((0, T_0))\bigr)$.  We will split our considerations into the cases $1/2 < s \leq 1$ and $0 \leq s < 1/2$.

\begin{theorem}
\label{Th6.1}
Let $1/2 < s \leq 1$.  Given $T_0 > 0$, there exists a linear operator $\HS = \HS_{T_0}$ (the homogeneous solution operator) on the subspace of $H^{(s + 1)/3}({\mathbb R}^+) \times H^s({\mathbb R}^+)$ of functions with the property that $f(0) = \phi(0)$ (in the sense of Proposition \ref{Pn2.3}) such that $$w = \HS(f, \phi) \in C\bigl((-\infty, +\infty); H^{(s + 1)/3}({\mathbb R}_t)\bigr) \cap C\bigl((-\infty, +\infty); H^s({\mathbb R}_x)\bigr)$$ and $w$ solves \pref{6.1} in the sense described above.  Moreover, $w(x, -)$ and $w(-, t)$ are continuous for each $x$ and $t$ and $w(x, 0) = \phi(x)$ and $w(0, t) = f(t)$ in the sense of Proposition \ref{Pn2.3} for $x \geq 0$ and $0 \leq t \leq T_0$.  In addition, $w$ satisfies the following estimates:
\begin{gather}
\label{6.2}\begin{split}
\sup_{\gamma \in {\mathbb R}} e^{-C|\gamma|}||D_x^{i\gamma}D_x^{s + 1}w||_{L_x^\infty L_t^2} \leq C||(f, \phi)||_{H^{(s + 1)/3} \times H^s} \\
||D_x^s\partial_xw||_{L_x^\infty L_t^2} \leq C||(f, \phi)||_{H^{(s + 1)/3} \times H^s}, \quad s < 1
\end{split} \\
\label{6.3}\sup_{\gamma \in {\mathbb R}} e^{-C|\gamma|}||D_x^{s - 1/4}D_x^{i\gamma}w||_{L_x^4L_t^\infty} \leq C||(f, \phi)||_{H^{(s + 1)/3} \times H^s} \\
\label{6.4}||D_x^sw||_{L_x^5L_t^{10}} \leq C||(f, \phi)||_{H^{(s + 1)/3} \times H^s} \\
\label{6.5}|||w|||_{X_{s, b}} \leq C||(f, \phi)||_{H^{(s + 1)/3} \times H^s}.
\end{gather}
(Recall that the space $X_{s, b}$ was defined in Remark \ref{Rk5.1}.)  The constants in all the estimates above depend only on $s$ and $T_0$.
\end{theorem}

\begin{proof}{}
We start out with the construction of $w$.  Let $\Psi_1$, $\Psi_2$ and $\Psi_3$ be cut-off functions, all supported on $[-2T_0, 2T_0]$ and identically 1 on $[-T_0, T_0]$, satisfying $\Psi_1.\Psi_2 = \Psi_1$ and $\Psi_2.\Psi_3 = \Psi_2$.  Let $\tilde{\phi}$ denote an extension of $\phi$ to all of $\mathbb R$ satisfying $\anm{\tilde{\phi}}{s}{\mathbb R} \leq C\anm{\phi}{s}{{\mathbb R}^+}$ (see Remark \ref{Rk2.1}, for instance, for the existence of a linear extension operator).  Let $\alpha(t) = S(t)\tilde{\phi}\bigl|_{x = 0}$ and $\tilde{\alpha} = \Psi_1.\alpha$.  By \pref{4.4}, $\tilde{\alpha} \in H^{(s + 1)/3}({\mathbb R}_t)$ and $\tilde{\alpha}(0) = \tilde{\phi}(0) = \phi(0)$.  Let $\tilde{f}$ be an extension of $f$ to all of $\mathbb R$ satisfying $\anm{\tilde{f}}{(s + 1)/3}{\mathbb R} \leq C\anm{f}{(s + 1)/3}{{\mathbb R}^+}$ and let $f_1 = \Psi_1.\tilde{f} - \tilde{\alpha}$.  Then, by Propositions \ref{Pn2.4} and \ref{Pn2.6}, we have that $$\anm{f_1}{(s + 1)/3}{\mathbb R} \leq C||(f, \phi)||_{H^{(s + 1)/3} \times H^s},$$ and, given our compatibility condition, we also have
\begin{equation}
\label{6.6}
\asnm{f_1}{(s + 1)/3}{{\mathbb R}^+%%Should it really be $_+$?
}{0} \leq C||(f, \phi)||_{H^{(s + 1)/3} \times H^s}.
\end{equation}
Let now $C_A$ be the constant in Proposition \ref{Pn4.1}, and $\FracInt_\alpha$ the Riemann-Liouville fractional integral studied in \ref{3}.  We define
\begin{equation}
\label{6.6(2)}%%This is the second (6.6) in the paper.  Henceforth I shall guess which reference refers to which.
\tilde{h}(t) = \frac{1}{C_A\Gamma(2/3)}\FracInt_{-2/3}(\tilde{f}_1)(t), \quad t \in {\mathbb R}^+
\end{equation}
By Propositions \ref{Pn3.1} and \ref{Pn3.2}, we have that
\begin{equation}
\label{6.7}
\asnm{\tilde{h}}{(s - 1)/3}{{\mathbb R}^+}{0} \leq C\asnm{\tilde{f}_1}{(s + 1)/3}{{\mathbb R}^+}{0} \leq C||(f, \phi)||_{H^{(s + 1)/3} \times H^s}.
\end{equation}
Finally, let $h = \Psi_3.\tilde{h}$.  By Propositions \ref{Pn2.5}, \ref{Pn2.6} and \ref{Pn2.8}, we have (setting $h(t) = 0$ for $t < 0$)
\begin{equation}
\label{6.8}
\adsnm{h}{(s - 1)/3}{{\mathbb R}^+}{0} \simeq \anm{h}{(s - 1)/3}{{\mathbb R}^+} \simeq \asnm{h}{(s - 1)/3}{\mathbb R}{0} \leq C||(f, \phi)||_{H^{(s + 1)/3} \times H^s}.
\end{equation}
We claim that
\begin{equation}
\label{6.9}
\Psi_2.\FracInt_{2/3}(h) = \frac{1}{C_A\Gamma(2/3)}\tilde{f}_1
\end{equation}
In fact, we know that $\FracInt_{2/3}(\tilde{h}) = \frac{1}{C_A\Gamma(2/3)}\tilde{f}_1$ from the remarks after Definition \ref{Df3.1} and density considerations.  Hence $\Psi_2.\FracInt_{2/3}(\tilde{h}) = \frac{1}{C_A\Gamma(2/3)}\tilde{f}_1$, and so, to verify the claim, all that we need to show is that $\Psi_2\FracInt_{2/3}(\tilde{h}) = \Psi_2\FracInt_{2/3}(h)$.  This is true, by the definition of $\FracInt_{2/3}$ and the properties of $\Psi_2$ and $\Psi_3$, for $h \in C_0^\infty({\mathbb R}^+)$, and hence in our case by density considerations.

We now define
\begin{equation}
\label{6.10}
w(x, t) = \Psi_2(t)\Bigl\{\int_0^t S(t - t')(\delta_0)(x)h(t')dt' + S(t)\tilde{\phi}(x)\Bigr\}.
\end{equation}
In view of \pref{6.8} and Lemmas \ref{Lm4.1} and \ref{Lm4.3}, we see that \pref{6.2}--\pref{6.4} hold.  \pref{6.5} holds because of Lemma \ref{Lm5.2} and Remark \ref{Rk5.1}.  The fact that $w \in C\bigl((-\infty, +\infty); H^{(s + 1)/3}({\mathbb R}_t)\bigr) \cap C\bigl((-\infty, +\infty); H^s({\mathbb R}_x%%Or should it really be ${\mathbb R}^\times$?
)\bigr)$ follows also from \pref{6.8} and Lemmas \ref{Lm4.1} and \ref{Lm4.3}.  The same is true about the continuity statements.  The fact that $w(x, 0) = \phi(x)$ for $x \geq 0$ follows from \pref{4.8} and Proposition \ref{Pn2.4}, while the fact that $w(0, t) = f(t)$ for $0 \leq t \leq T_0$ follows from Lemma \ref{Lm4.2}, density, \pref{6.9}, our definition of $f_1$ and the fact that $\Psi_i \equiv 1$ on $[0, T_0]$.  Finally, the fact that the equations in \pref{6.1} hold follows from \pref{4.7} and density.
\end{proof}

\begin{theorem}
\label{Th6.2}
Let $0 \leq s < 1/2$.  Given $T_0 > 0$, there exists a linear operator $\HS = \HS_{T_0}$ on $H^{(s + 1)/3}({\mathbb R}^+) \times H^s({\mathbb R}^+)$ such that $$w = \HS(f, \phi) \in C\bigl((-\infty, +\infty); H^{(s + 1)/3}({\mathbb R}_t)\bigr) \cap C\bigl((-\infty, +\infty); H^s({\mathbb R}_x)\bigr),$$ and $w$ solves \pref{6.1} in the sense described above.  Moreover, $w$ satisfies the estimates \pref{6.2}, \pref{6.4} and \pref{6.5} (and \pref{6.3} if $s \geq 1/4$).
\end{theorem}

\begin{proof}{}
The construction and proof are identical to that of Theorem \ref{Th6.1}, once one invokes Lemma \ref{Lm2.3} and Proposition \ref{Pn2.6}, so that $H^\alpha_0({\mathbb R}^+) = H^\alpha({\mathbb R}^+)$ for $0 \leq \alpha < 1/2$.  This explains the fact that there is no compatibility condition in this range.
\end{proof}

\begin{remark}
\label{Rk6.1}
Global (in time) versions of the results in Theorem \ref{Th6.2} can be obtained by using the homogeneous Sobolev spaces $\dot{H}^{(s + 1)/3}({\mathbb R}^+)$ and $\dot{H}^{(s + 1)/3}_0({\mathbb R}^+)$, the fact that, for $0 \leq s < 1/2$, $\dot{H}^{(s + 1)/3}({\mathbb R}^+) = \dot{H}^{(s + 1)/3}_0({\mathbb R}^+)$ and the fact that the operators $\FracInt_\alpha$ act well on these homogeneous global spaces.  For instance, when $s = 0$, given $(f, \phi) \in \dot{H}^{1/3}({\mathbb R}_t^+) \times L^2({\mathbb R}_x^+)$, we can find $w$ in $C\bigl((-\infty, +\infty); \dot{H}^{1/3}({\mathbb R}_t)\bigr) \cap C\bigl((-\infty, +\infty); L^2({\mathbb R}_x)\bigr)$ such that $w(x, 0) = \phi(x)$ and $w(0, t) = f(t)$ for $x, t > 0$, $\partial_tw%%Or is it $v$?
+ \partial_x^3w = 0$ in ${\cal D}'({\mathbb R}^+ \times {\mathbb R}^+)$ and $w$ satisfies the estimate:
\begin{equation}
\label{6.11}
\sup_\gamma e^{-C|\gamma|}||D_x^{i\gamma}D_xw||_{L_x^\infty L_t^2} + ||\partial_xw||_{L_x^\infty L_t^2} + ||w||_{L_x^5L_t^{10}} \leq C||(f, \phi)||_{\dot{H}^{1/3} \times L^2}.
\end{equation}
This follows as in the proof of Theorem \ref{Th6.2}, leaving out the cut-off functions.  Similar results hold on $\dot{H}^{(s + 1)/3}({\mathbb R}_t^+) \times \dot{H}^s({\mathbb R}_x^+)$ for $0 \leq s < 1/2$.
\end{remark}

We will now turn to our uniqueness theorem, which will imply that the solutions constructed in Theorems \ref{Th6.1} and \ref{Th6.2} and Remark \ref{Rk6.1} are unique (when restricted to ${\mathbb R}_x^+ \times [0, T_0]$).

\begin{theorem}
\label{Th6.3}
Let $w$ be a solution to \pref{6.1}, with $s = 0$ and $\phi$ and $f$ identically 0.  Then $w \equiv 0$ on ${\mathbb R}_x^+ \times [0, T_0]$.
\end{theorem}

\begin{proof}{}
We begin by regularizing our solutions, so that certain identities can be established.  First, let $$\tilde{w}(x, t) = \begin{cases}
w(x, t) & (x, t) \in {\mathbb R}^+ \times [0, T_0]     \\
0       & (x, t) \not\in {\mathbb R}^+ \times [0, T_0]
\end{cases}$$ and fix $\theta \in C_0^\infty({\mathbb R})$ satisfying $\theta \geq 0$, $\int \theta = 1$ and $\supp \theta \subseteq [-1, -1/2]$.  For $\varepsilon, \delta > 0$, define
\begin{equation}
\label{6.12}
w_{\varepsilon, \delta}(x, t) = \int\bic\int \tilde{w}(y, t')\frac{\theta}{\varepsilon}\biggl(\frac{x - y}{\varepsilon}\biggr)\mult\frac{\theta}{\delta}\biggl(\frac{t - t'}{\delta}\biggr)dy\,dt'.
\end{equation}
We start out by making a few observations about $w_{\varepsilon, \delta}$.  Fix $0 < T < T_0$, and consider only $\delta < T_0 - T$.  Suppose that $(x, t) \in {\mathbb R}^+ \times (0, T)$.  Then, the integration in the definition of $w_{\varepsilon, \delta}$ takes place on $\{y : \frac{\varepsilon}{2} + x < y < x + \varepsilon\} \times \{t' : \frac{\delta}{2} + t < t' < t + \delta\}$, and hence, for such $(x, t)$, $\tilde{w}(y, t') = w(y, t')$%%Where'd the $\varepsilon$ and $\delta$ go?
.  Thus, since $w_{\varepsilon, \delta} \in C^\infty({\mathbb R} \times {\mathbb R})$, $\partial_tw_{\varepsilon, \delta} + \partial_x^3w_{\varepsilon, \delta} \equiv 0$ on ${\mathbb R}^+ \times (0, T)$.  Notice also that $w_{\varepsilon, \delta} \in C\bigl([0, \infty); H^{1/3}((0, T))\bigr) \cap C\bigl([0, T]; L^2({\mathbb R}^+)\bigr)$, uniformly in $(\varepsilon, \delta)$.  We will now establish the identity
\begin{equation}
\label{6.13}
\partial_t\ic\int\blim_{x > x_0} w_{\varepsilon, \delta}^2(x, t)dx = 2\partial_x^2w_{\varepsilon, \delta}(x_0, t)\!\cdot\!w_{\varepsilon, \delta}(x_0, t) - \bigl(\partial_xw_{\varepsilon, \delta}(x_0, t)\bigr)^2
\end{equation}
for $x_0 \geq 0$ and $0 \leq t \leq T$.  In fact, for $\varepsilon$, $\delta$ fixed, $\partial_tw_{\varepsilon, \delta}, \partial_x^3w_{\varepsilon, \delta} \in L^2({\mathbb R}^+)$ for each such fixed $t$, and hence the left-hand side of \pref{6.13} equals $$2\ic\int\blim_{x > x_0} \partial_tw_{\varepsilon, \delta}(x, t)\!\cdot\!w_{\varepsilon, \delta}(x, t)dx = -2\ic\int_{x > x_0} \partial_x^3w_{\varepsilon, \delta}(x, t)\!\cdot\!w_{\varepsilon, \delta}(x, t)dx.$$  Note that, since $\partial_x^3w_{\varepsilon, \delta}, \partial_x^2w_{\varepsilon, \delta}, \partial_xw_{\varepsilon, \delta}, w_{\varepsilon, \delta} \in L^2({\mathbb R}^+)$, we can certainly find $x_n \rightarrow +\infty$ so that $\partial_x^jw_{\varepsilon, \delta}(x_n, t) \rightarrow 0$ for all $j = 0, 1, 2, 3$.  Then we see that the right-hand side above equals
\begin{multline*}
2\partial_x^2w_{\varepsilon, \delta}(x_0, t)\!\cdot\!w_{\varepsilon, \delta}(x_0, t) + 2\ic\int\blim_{x > x_0} \partial_x^2w_{\varepsilon, \delta}(x, t)\!\cdot\!\partial_xw_{\varepsilon, \delta}(x, t)dx \\
= 2\partial_x^2w_{\varepsilon, \delta}(x_0, t)\!\cdot\!w_{\varepsilon, \delta}(x_0, t) - \bigl(\partial_xw_{\varepsilon, \delta}(x_0, t)\bigr)^2,
\end{multline*}
as claimed.  Fix $0 < t_0 \leq T$, and integrate \pref{6.13} between $0$ and $t_0$ to obtain
\begin{multline}
\label{6.14}
\int\blim_{x > x_0} w_{\varepsilon, \delta}^2(x, t_0)dx - \int\blim_{x > x_0} w_{\varepsilon, \delta}^2(x, 0)dx \\
= 2\ic\int_0^{t_0} \partial_x^2w_{\varepsilon, \delta}(x_0, t)\!\cdot\!w_{\varepsilon, \delta}(x_0, t)dt - \int_0^{t_0} \bigl(\partial_xw_{\varepsilon, \delta}(x_0, t)\bigr)^2dt \\
\leq 2\int_0^{t_0} \partial_x^2w_{\varepsilon, \delta}(x_0, t)\!\cdot\!w_{\varepsilon, \delta}(x_0, t)dt.
\end{multline}

Fix now $0 < \delta < T_0 - T$, and define $w_\delta(x, t) = \int \tilde{w}(x, t')\frac{\theta}{\delta}\bigl(\frac{t - t'}{\delta}\bigr)dt'$.

\begin{claim}
\label{6.15}
$w_{\varepsilon, \delta}, \partial_xw_{\varepsilon, \delta}, \partial_x^2w_{\varepsilon, \delta} \in C\bigl([0, 1]; L^2([0, T])\bigr)$, uniformly in $\varepsilon > 0$, for $\delta > 0$ fixed.  In fact, since $w \in C\bigl([0, 2]; H^{1/3}((0, T_0))\bigr)$, this easily follows for $w_{\varepsilon, \delta}$, and for $\partial_tw_{\varepsilon, \delta}$.  By the equation, this also follows for $\partial_x^3w_{\varepsilon, \delta}$.  Next, it is easy to check that $\partial_x^2w_{\varepsilon, \delta} \in L^\infty\bigl([0, 1]; L^2([0, T])\bigr)$, uniformly in $\varepsilon$.  In fact, let $f \in L^2\bigl([0, T]\bigr)$ satisfy $||f||_{L^2} = 1$, and consider $F(x) = \int_0^T w_{\varepsilon, \delta}(x, t)f(t)dt$.  Then, as is well-known, $$F(x_0 + h) + F(x_0 - h) - 2F(x_0) = F''(\overline{x})h^2$$ for some $\overline{x} \in (x_0 - h, x_0 + h)$, for each $h > 0$.  Choose $x_0 = 1/2$ and $h = 1/4$ to conclude that $$\Bigl|\int_0^T \partial_x^2w_{\varepsilon, \delta}(\overline{x}, t)f(t)dt\Bigr| \leq C_\delta$$ since $w_{\varepsilon, \delta} \in C\bigl([0, 1]; L^2([0, T])\bigr)$ -- but, if we consider any fixed $x \in [0, 1]$,
\begin{multline*}
\int_0^T \partial_x^2w_{\varepsilon, \delta}(x, t)f(t)dt = \int_0^T \bigl[\partial_x^2w_{\varepsilon, \delta}(x, t) - \partial_x^2w_{\varepsilon, \delta}(\overline{x}, t)\bigr]f(t)dt \\
+ \int_0^T \partial_x^2w_{\varepsilon, \delta}(\overline{x}, t)f(t)dt,
\end{multline*}
and, using the fact that $\partial_x^3w_{\varepsilon, \delta} \in C\bigl([0, 1]; L^2([0, T%%Or should it really be 1?
])\bigr)$, our claim follows.  Once we know this, our estimate on $\partial_x^3w_{\varepsilon, \delta}$ shows that $\partial_x^2w_{\varepsilon, \delta} \in C\bigl([0, 1]; L^2([0, T])\bigr)$, uniformly in $\varepsilon$, and in fact they are equicontinuous in $\varepsilon$.  The statement for $\partial_xw_{\varepsilon, \delta}$ follows similarly.
\end{claim}

\begin{claim}
\label{6.16}
For any fixed $0 \leq x_0 \leq 1$, $w_{\varepsilon, \delta}(x_0, -) \xrightarrow{\varepsilon \rightarrow 0} w_\delta(x_0, -)$ in $L^2\bigl([0, T]\bigr)$.  Note that this is immediate from the definitions and the fact that $w_\delta \in C\bigl([0, 1]; L^2([0, T])\bigr)$ since $w$ is in $C\bigl([0, 1]; L^2([0, T_0])\bigr)$.
\end{claim}

\begin{claim}
\label{6.17}
For any fixed $0 \leq x_0 \leq 1$ and $0 \leq t_0 \leq T$, $$\int\blim_{x > x_0} w_{\varepsilon, \delta}^2(x, t_0)dx \xrightarrow{\varepsilon \rightarrow 0} \int\blim_{x > x_0} w_\delta^2(x, t_0)dx.$$  This follows easily from the fact that $w_\delta \in C\bigl([0, T]; L^2({\mathbb R}^+)\bigr)$.  Now, using \pref{6.14}--\pref{6.16%%Do we really want (6.17) here?  That's *this* claim.
}, we see that, for $x_0 \in [0, 1]$,
\begin{equation}
\label{6.18}
\int\blim_{x > x_0} w_\delta^2(x, t_0)dx - \int\blim_{x > x_0} w_\delta^2(x, 0)dx \leq C_\delta||w_\delta(x_0, -)||_{L^2([0, T])}.
\end{equation}
Next, note that, as $x_0 \rightarrow 0$, $||w_\delta(x_0, -)||_{L^2([0, T])} \rightarrow ||w_\delta(0, -)||_{L^2([0, T])}$ since $w_\delta \in C\bigl([0, 1]; L^2([0, T])\bigr)$, which follows from $w \in C\bigl([0, 1]; L^2([0, T])\bigr)$.  By our assumption and the definition of $w_\delta$, however, $w_\delta(0, -) \equiv 0$ on $[0, T]$.  Thus, it is easy to see that \pref{6.18} yields
\begin{equation}
\label{6.19}
\int\blim_{x > 0} w_\delta^2(x, t_0)dx \leq \int\blim_{x > 0} w_\delta^2(x, 0)dx.
\end{equation}
Since, however, $w \in C\bigl([0, T]; L^2({\mathbb R}^+)\bigr)$, as $\delta \rightarrow 0$ we obtain $$\int\blim_{x > 0} w^2(x, t_0)dx \leq \int\blim_{x > 0} w^2(x, 0)dx = 0,$$ and hence $w \equiv 0$ on ${\mathbb R}^+ \times [0, T]$ as desired.
\end{claim}
\end{proof}

\subsection{Inhomogeneous linear solution operator}

Next, we turn our attention to constructing solutions to the inhomogeneous linear quarter-plane problem
\begin{equation}
\label{6.20}
\begin{cases}
\partial_tw + \partial_x^3w = h, & x > 0, t \in (0, T_0) \\
w(x, 0) = 0,                     & x > 0                 \\
w(0, t) = 0,                     & t \in [0, T_0].
\end{cases}
\end{equation}

\begin{theorem}
\label{Th6.4}
Assume that $0 \leq s \leq 1$, $s \neq 1/2$ and $h \in L^2_{T_0%%Isn't $T_0$ a constant?
}H^s({\mathbb R}_x^+)$.  Given $T_0$, there exists a linear operator $\IHS = \IHS_{T_0}$ (the inhomogeneous solution operator) such that $$w = \IHS(h) \in C\bigl((-\infty, +\infty); H^{(s + 1)/3}({\mathbb R}_t)\bigr) \cap C\bigl((-\infty, +\infty); H^s({\mathbb R}_x)\bigr),$$ and $w$ solves \pref{6.20} in the sense that the equation holds in ${\cal D}'\bigl({\mathbb R}^+ \times (0, T_0)\bigr)$, and the `lateral values' are taken in the sense of $C\bigl((-\infty, +\infty); H^{(s + 1)/3}((0, T_0))\bigr)$ and the initial ones in the sense of $C\bigl([0, T_0]; H^s({\mathbb R}_x^+)\bigr)$.  (If $s > 1/2$, this also holds in the pointwise sense.)  In addition, for $T_0 \leq 1$, $w$ satisfies the following estimates:
\setcounter{subtheorem}{0}\comment{How to fix this numbering?}
\begin{list}{(\thesubtheorem)}{\usecounter{subtheorem}}%%We want to tie this to the equation numbering.
\item\label{6.21}  We have $$\sup_t \anm{w(-, t)}{s}{{\mathbb R}_x} \leq CT_0^{\beta(s)}||h||_{L_{T_0}^2H^s({\mathbb R}_x^+)}.$$
\item\label{6.22}  We have $$\sup_x \anm{w(x, -)}{(s + 1)/3}{{\mathbb R}_t} \leq CT_0^{\beta(s)}||h||_{L_{T_0}^2H^s({\mathbb R}_x^+)}.$$
\item\label{6.23}  We have $$\sup_x \asnm{w(x, -)}{(s + 1)/3}{{\mathbb R}_t^+}{0} \leq CT_0^{\beta(s)}||h||_{L_{T_0}^2H^s({\mathbb R}_x^+)}.$$
\item\label{6.24}  The quantities on the left-hand sides of \pref{6.2}, \pref{6.3} and \pref{6.4}, in case $s \geq 1/4$, are controlled by $CT_0^{\beta(s)}||h||_{L_{T_0}^2H^s({\mathbb R}_x^+)}$,
\end{list}
where $\beta(s) = \frac{1}{3} - \frac{s}{6}$.
\end{theorem}

\begin{proof}{}
Let $$w_1%%Is this a 1 or an $s$?
(x, t) = \int_0^t S(t - t')\tilde{h}(x, t')dt',$$ where $\tilde{h} \in L^2\bigl((-\infty, +\infty); H^s({\mathbb R}_x)\bigr)$ and $\tilde{h}$ satisfies $$\tilde{h}(x, t) = \begin{cases}
h(x, t), & (x, t) \in {\mathbb R}^+ \times [0, T_0] \\
0,       & t > T_0\text{ or }t < 0
\end{cases}$$ and $||\tilde{h}||_{L_t^2H^s({\mathbb R})} \leq C||h||_{L_{T_0}^2H^s({\mathbb R}_x^+)}$.  (Moreover, by Remark \ref{Rk2.1}, we can take $\tilde{h} = E(h)$, where $E$ is a linear extension operator in the $x$ variable.)  We know that $\partial_tw_1 + \partial_x^3w_1 = \tilde{h}$ in ${\cal D}'({\mathbb R} \times {\mathbb R})$.  We next consider the estimates \pref{6.21}--\pref{6.24}\comment{Fix these references.} for $w_1$.  Minkowski's integral inequality and the corresponding estimates for the group (\pref{4.1}, \pref{4.2}, \pref{4.5}) show that (in the case of \pref{6.21}\comment{Fix this reference.}, for instance), with $\Psi \equiv 1$ on $[-T_0, T_0]$ and $\supp \Psi \subseteq [-2T_0, 2T_0]$,
\begin{multline}
\label{6.25}
\Bigl|\Bigl|\Psi(t)\ic\int_0^t S(t - t')\tilde{h}(x, t')dt'\Bigr|\Bigr|_{H^s({\mathbb R}_x)} \leq \int_{-2T_0}^{2T_0} \anm{S(t - t')\tilde{h}(x, t')}{s}{{\mathbb R}_x}dt' \\
\leq C\ic\int_0^{T_0} \anm{h(-, t')}{s}{{\mathbb R}_x^+}dt' \leq CT_0^{1/2}||h||_{L_{T_0}^2H^s({\mathbb R}_x^+)}
\end{multline}
and, since $\beta(s) \leq 1/2$, $\Psi(t)w_1(x, t)$ satisfies the estimates \pref{6.21} and \pref{6.24}.  Finally, translation invariance and \pref{4.48} give \pref{6.22} and \pref{6.23}\comment{Fix the references in this sentence and the last one.} for $\Psi(t)w_1(x, t)$, as well as the correct continuity and trace statements (see Remark \ref{Rk4.2}).  Next, let $f(t) = \Psi(t)w_1(0, t)$, $T = \max\{2T_0, 1\}$ and $w_2 = \HS_T(f, 0)$, where $\HS_T$ is the operator constructed in Theorems \ref{Th6.1} and \ref{Th6.2}.  (Note that \pref{4.48} shows that $f \in H^{(s + 1)/3}_0({\mathbb R}^+)$, and hence, for $s > 1/2$, the compatibility condition in Theorem \ref{6.1} is satisfied.)  Finally, let $w(x, t) = \Psi(t)w_1(x, t) - w_2(x, t)$.  Estimate \pref{4.48}, together with Theorems \ref{Th6.1} and \ref{Th6.2} and the above estimates for $w_1$, now gives the desired result.
\end{proof}

\old{\begin{theorem}
\label{Th6.5}
Suppose that $h = \tilde{h}\big|_{{\mathbb R}_x^+ \times [0, T_0]}$, with $\tilde{h} \in Y^{s, \alpha}$ for some $0 \leq s \leq 1$, $s \neq 1/2$, and $\alpha > 1/2$.  Then $w = \IHS_{T_0}(\tilde{h})$ belongs to $C\bigl((-\infty, +\infty); H^{(s + 1)/3}({\mathbb R}_t)\bigr) \cap C\bigl((-\infty, +\infty); H^s({\mathbb R}_x)\bigr)$, and $w$ solves \pref{6.20} in the sense of Theorem \ref{Th6.4}.  Moreover, $w$ satisfies the following estimates:
\begin{gather}
\label{6.26}\sup_t \anm{w(-, t)}{s}{{\mathbb R}_x} \leq C_{T_0}|||\tilde{h}|||'_{s, \alpha}, \\
\label{6.27}\sup_x \anm{w(x, -)}{(s + 1)/3}{{\mathbb R}_t} \leq C_{T_0}|||\tilde{h}|||'_{s, \alpha}, \\
\label{6.28}\sup_x \asnm{w(x, -)}{(s + 1)/3}{{\mathbb R}_t}{0} \leq C_{T_0}|||\tilde{h}|||'_{s, \alpha} \\
\intertext{and}
\label{6.29}|||w|||_{s, \alpha} \leq C_{T_0}|||\tilde{h}|||'_{s, \alpha}.
\end{gather}
\end{theorem}

\begin{proof}{}
By Lemma \ref{Lm5.2}, $\Psi.w_1$ satisfies \pref{6.26}--\pref{6.29}.  Moreover, $f \in H^{(s + 1)/3}_0({\mathbb R}_t^+)$ by \pref{5.6}, and hence $w_2$ satisfies \pref{6.26}--\pref{6.29} by Theorem \ref{Th6.1}.  This concludes the proof.
\end{proof}}

\begin{theorem}
\label{Th6.6}
Let $w = \IHS_{T_0}(h)$.  Then, if $h \in L_x^{5/4}L_t^{10/9}$, then $$w \in C\bigl((-\infty, +\infty); L^2({\mathbb R}_x)\bigr) \cap C\bigl((-\infty, +\infty); H^{1/3}({\mathbb R}_t)\bigr),$$ $\pd{w}{x} \in L_x^\infty L_t^2$, $w \in L_x^5L_t^{10}$ (all with norm control) and $w$ solves \pref{6.20} in the sense of Theorem \ref{Th6.4}, with $s = 0$.
\end{theorem}

\begin{proof}{}
Clearly, it suffices to establish the estimates.  Let $w_1(x, t) = \int_0^t S(t - t')h(x, t')dt'$.  Then, by \pref{4.47}, $$||w_1||_{L_x^5L_t^{10}} \leq C||h||_{L_x^{5/4}L_t^{10/9}}.$$  To show that $\pd{w}{x} \in L_x^\infty L_t^2$, we first note that $$\pd{w}{x} = {\cal H}D_xw = D_x\ic\int_0^t S(t - t'){\cal H}(h)(x, t')dt',$$ with $\cal H$ the Hilbert transform on functions in the $x$-variable.  Since $\cal H$ is bounded on $L_x^{5/4}L_t^{10/9}$ (see \cite{15}), it suffices to show that $D_xw \in L_x^\infty L_t^2$.  Using \pref{4.16} and duality reduces matters to checking that each one of the terms in \pref{4.16} maps $L_x^1L_t^2$ into $L_x^5L_t^{10}$.  This, in turn, follows by complex interpolation between the second inequality for $s = 1$ in \pref{4.7} and the first inequality in \pref{4.11}, together with translation invariance (to pass from $\delta_0 \otimes h$ to $\delta_{x_0} \otimes h$) and Minkowski's integral inequality.  (The proof also follows from Proposition 2.3 in \cite{13}.)  The fact that $w_1 \in C\bigl((-\infty, +\infty); \dot{H}^{1/3}({\mathbb R}_t)\bigr)$ follows from a familiar density argument, together with the inequality
\begin{equation}
\label{6.30}\comment{Apparently this and the next equation are numbered badly.}
||w_1||_{L_x^\infty\dot{H}^{1/3}_t} \leq C||h||_{L_x^{5/4}L_t^{10/9}}.
\end{equation}
Once again, \pref{6.30} follows from \pref{4.16}, duality and the first inequality in \pref{4.7}, together with translation invariance and Minkowski's integral equality.  The fact that $w_1 \in C\bigl((-\infty, +\infty); L^2({\mathbb R}_x)\bigr)$ follows by duality from \pref{4.5} with $s = 0$, in the same manner as \pref{4.43} %%Is this a local reference, or a reference to [10]?
is established in (3.7) of \cite{10}.  The fact that $H^{1/3}({\mathbb R}_t^+) = H^{1/3}_0({\mathbb R}_t^+)$ (see Proposition \ref{Pn2.6}) and the Hardy-Littlewood-Sobolev inequality
\begin{equation}
\label{6.31}
||f||_{L^6({\mathbb R})} \leq C||D_t^{1/3}f||_{L^2({\mathbb R})},
\end{equation}
together with the Leibniz rule (Theorem A.12 in \cite{10}), now allow us to establish the desired bounds for $\Psi.w_1$.  Finally, since $f = \Psi w_1(0, -) \in H^{1/3}_0({\mathbb R}_t^+)$, the corresponding bounds for $w_2$ follow from Theorem \ref{Th6.2}.
\end{proof}

\begin{remark}
\label{Rk6.2}
A global (in time) version of Theorem \ref{Th6.6} can be established by using the homogeneous Sobolev space $\dot{H}^{1/3}$, together with Remark \ref{Rk6.1} and the proof above, and omitting the cut-off function.  We obtain a solution $w$ defined on %%Or do we really want `in'?
${\mathbb R}^+ \times {\mathbb R}^+%%Do we really want those $+$s?
$ which lies in $$C\bigl((-\infty, +\infty); L^2({\mathbb R}_x)\bigr) \cap C\bigl((-\infty, +\infty); \dot{H}^{1/3}({\mathbb R}_t)\bigr) \cap L_x^5L_t^{10}({\mathbb R} \times {\mathbb R})$$ and satisfies $\pd{w}{x} \in L_x^\infty L_t^2$.
\end{remark}

The final result of this section is a Bourgain space analogue of Theorem \ref{Th6.4}.  We show that the inhomogeneous solution operator $\IHS$ sends the inhomgeneity $w \in Y_b$ into $X_b$ functions with well-defined traces along $\{x = 0\}$ and $\{t = 0\}$ as $H_t^{1/3}$ and $L_x^2$ functions.

\begin{theorem}
\label{Th6.6'}
Suppose that $h = \tilde{h}\bigr|_{{\mathbb R}_x^+ \times [0, T_0]}$, with $\tilde{h} \in Y_{s, b}$ for some $0 \leq s \leq 1$, $s \neq \frac{1}{2}$.  Then $w = \IHS_{T_0}(\tilde{h}\comment{Or just $h$?})$ belongs to $C\bigl((-\infty, \infty); H^{(s + 1)/3}({\mathbb R}_t)\bigr) \cap C\bigl((-\infty, \infty); H^s({\mathbb R}_x)\bigr)$, and $w$ solves \pref{6.20} in the sense of Theorem \ref{Th6.4}.  Moreover, $w$ satisfies the following estimates:
\begin{gather}
\label{6.25'} \sup_t \anm{w(-, t)}{s}{{\mathbb R}_x} \leq C_{T_0}||h||_{Y_{s, b}} \\
\label{6.26'} \sup_x \anm{w(x, -)}{(s + 1)/3}{{\mathbb R}_t} \leq C_{T_0}||h||_{Y_{s, b}} \\
\label{6.27'} \sup_x \asnm{w(x, -)}{(s + 1)/3}{{\mathbb R}_t}{0} \leq C_{T_0}||h||_{Y_{s, b}} \\
\intertext{and}
\label{6.28'} ||w||_{X_{s, b}} \leq C_{T_0}||h||_{Y_{s, b}}.
\end{gather}
\end{theorem}

\begin{proof}
Lemma \ref{Lm5.4} implies \pref{6.28'}.  Lemma \ref{Lm5.5} and the results of \ref{2} imply \pref{6.25'}--\pref{6.27'}.
\end{proof}

\begin{remark}
\label{Rk6.3}
Theorem \ref{Th6.3} implies the uniqueness of the restriction to ${\mathbb R}^+ \times [0, T_0]$ of the solutions constructed in Theorems \ref{Th6.4}--\ref{Th6.6'}.
\end{remark}

\nw{\begin{remark}
\label{Rk6.4}
The methods in the previous section also apply, with minor modifications, to the homogeneous and inhomogeneous problems with a transport term $c\partial_xu$, $c \neq 0$.  For example, for the homogeneous problem, $$\begin{cases}
\partial_tw + \partial_x^3w + c\partial_xw = 0, & x > 0, t \in (0, T_0) \\
w(x, 0) = \phi(x),                              & x > 0                 \\
w(0, t) = f(t),                                 &
\end{cases}$$ one still introduces the forced initial value problem $$\begin{cases}
\partial_t\tilde{w} + \partial_x^3\tilde{w} + c\partial_x\tilde{w} = \delta_0(x)g(t), & x \in {\mathbb R}, t \in (0, T_0) \\
\tilde{w}(x, 0) = \tilde{\phi}(x), &
\end{cases}$$ where the forcing function $g$ is selected to ensure that $\tilde{w}(0, t) = \tilde{f}(t)$, $t \in (0, T_0)$.  To see that this can be done (for $T_0$ small), let \comment{Should there be a $~$ over the $\phi$s below?} $$\tilde{S}(t)\phi(x) = \int e^{i(x\xi + t(\xi^3 + c\xi))}\hat{\phi}(\xi)d\xi$$ and consider the integral equation $$\int_0^t \tilde{S}(t - t')\delta_0(x)g(t')dt'\Bigr|_{\{x = 0\}} = \tilde{f}(t) - \tilde{S}(t)\phi(x)\bigr|_{x = 0} = \Tilde{\Tilde{f}}(t).$$  Using the notation of \ref{2}, \ref{3} and \ref{4}, the left-hand side becomes $$\int_0^t \frac{1}{(t - t')^{1/3}}A\bigl((t - t')^{2/3}c\bigr)g(t')dt' = A_c(g)(t).$$  Note that $A(0) \neq 0$, and $A$ is differentiable at $0$, since $\int e^{i\xi^3}\xi\,d\xi = \int e^{i\eta}\frac{\eta^{1/3}}{\eta^{2/3}}d\eta$.  Thus, it is easy to see that $A_c(g) - A(0)\int_0^t \frac{g(t')dt'}{(t - t')^{1/3}}$ maps $H_0^{-1}({\mathbb R}^+)$ into $H_0^{1/3}({\mathbb R}^+)$ and $L^2({\mathbb R}^+)$ into $H_0^1\bigl((0, T_0]\bigr)$, with norm small in $T_0$.  This, combined with the results in \ref{2}, \ref{3} and \ref{4}, gives the invertibility of $A_c(g)$ from $H_0^{(s - 1)/3\comment{Or $s - 1/3$?}}\bigl((0, T_0]\bigr)$ onto $H_0^{(s + 1)/3\comment{Or $s + 1/3$?}}\bigl((0, T_0]\bigr)$ for small $T_0$.  The oscillatory integral estimates in \ref{4} (for small time) for the multipliers $e^{it\xi^3}$ and $e^{it(\xi^3 + c\xi)}$ are identical, while the estimates in \ref{5}, in the Bourgain spaces
\begin{multline*}
\tilde{X}_b = \Bigl\{f \in {\cal S}'({\mathbb R}^2) : \Bigl(\int\bic\int (1 + |\lambda - (\xi^3 + c\xi)|)^{2b}|\hat{f}(\xi, \lambda)|^2d\xi\,d\lambda\Bigr)^{1/2} \\
+ \Bigl(\int\bic\int\blim_{|\xi| < 1} (1 + |\lambda|)^{2\alpha}|\hat{f}(\xi, \lambda)|^2d\xi\,d\lambda\Bigr)^{1/2}\Bigr\},
\end{multline*}
are also identical.  Thus, the results in Theorems \ref{Th6.1}, \ref{Th6.2}, \ref{Th6.4}, \ref{Th6.6} and \ref{Th6.6'} extend to this setting.  The uniqueness result Theorem \ref{Th6.3} does hold too (with similar proof).%, but only in the case when $c > 0$\comment{It's written $\geq$ in the text, but didn't we prohibit $c = 0$?}.
\end{remark}}

\newsect{Well-posedness results for the nonlinear problems}{7}

In this section we will deal with well-posedness results for the nonlinear quarter-plane problem
\begin{equation}
\label{7.1}
\begin{cases}
\partial_tu + \partial_x^3u + u^k\partial_xu, & x > 0, t \in [0, T], k \in {\mathbb N} \\
u(x, 0) = \phi(x),                            & x > 0         \\
u(0, t) = f(t),                               & t \in [0, T].
\end{cases}
\end{equation}
Using the solutions to the corresponding linear problems and their estimates, as presented in \ref{6}, our results will follow from the methods in \cite{5} and \cite{10}, yielding identical results.  To avoid a lengthy discussion, we have decided to detail only certain ``highlight'' results.  The techniques presented in \ref{6} and the ones in this section do, however, yield in full the results in \cite{5} and \cite{10}.

\subsection{Well-posedness of generalized KdV in mixed-norm spaces}

We start out by discussing the case $k \geq 2$, first presenting a relatively straightforward result which gives the local well-posedness for $\phi \in H^s({\mathbb R}^+)$, $f \in H^{(s + 1)/3}({\mathbb R}^+)$, $1/2 < s \leq 1$ and $k \geq 2$.

\begin{theorem}
\label{Th7.1}
If $k \geq 2$ and $1/2 < s \leq 1$, then, given $(f, \phi) \in H^{(s + 1)/3}({\mathbb R}^+) \times H^s({\mathbb R}^+)$ satisfying the compatibility condition $f(0) = \phi(0)$, there exists $$T = T(||(f, \phi)||_{H^{(s + 1)/3}({\mathbb R}^+) \times H^s({\mathbb R}^+)})$$ such that the `integral equation'
\begin{equation}
\label{7.2}
w = \HS_1\comment{Or is it $T$?}(f, \phi) + \IHS_{2T}(\Psi.w^k.\partial_xw),
\end{equation}
where $\Psi \equiv 1$ on $[-T, T]$ and $\supp \Psi \subseteq [-2T, 2T]$, has a unique fixed point in the space $$B = \bigl\{w \in C\bigl((-\infty, +\infty); H^s({\mathbb R})\bigr) \cap C\bigl((-\infty, +\infty); H^{(s + 1)/3}({\mathbb R})\bigr) : \Lambda(w) < \infty\},$$ where $\Lambda(w) = \max\inlim_{1 \leq i \leq 6} \lambda_i(w)$ for
\begin{align*}
\lambda_1(w) & = \sup_x \anm{w(x, -)}{(s + 1)/3}{\mathbb R}, \\
\lambda_2(w) & = \sup_t \anm{w(-, t)}{s}{\mathbb R}, \\
\lambda_3(w) & = \sup_\gamma e^{-C|\gamma|}||D_x^{i\gamma}D_x^{s + 1}w||_{L_x^\infty L_t^2} + ||D_x^s\partial_xw||_{L_x^\infty L_t^2} \\
\intertext{(where the second term is added only when $s < 1$),}
\lambda_4(w) & = \sup_\gamma e^{-C|\gamma|}||D_x^{i\gamma}D_xw||_{L_x^\infty L_t^2} + ||\partial_xw||_{L_x^\infty L_t^2,} \\
\lambda_5(w) & = \sup_\gamma e^{-C|\gamma|}||D_x^{s - 1/4}D_x^{i\gamma}w||_{L_x^4L_t^\infty} \\
\intertext{and}
\lambda_6(w) & = \sup_\gamma e^{-C|\gamma|}||D_x^{i\gamma}w||_{L_x^4L_t^\infty}.
\end{align*}
The resulting $w$ solves \pref{7.1} in ${\mathbb R}^+ \times [0, T]$.
\end{theorem}

\begin{proof}{}
We will show that, for $T$ small, the mapping $$M_{(f, \phi)}(w) = \HS_1%%Or should it be $T$?
(f, \phi) + \IHS_{2T}(\Psi.w^k.\partial_xw)$$ is a contraction on $B_a = \{w \in B : \Lambda(w) \leq a\}$ for suitable $a$.  First note that Theorems \ref{Th6.1} and \ref{Th6.2} show that $\HS_1(f, \phi) \in B$ and $$\Lambda\bigl(\HS_1(f, \phi)\bigr) \leq C||(f, \phi)||_{H^{(s + 1)/3} \times H^s}.$$  Now Theorem \ref{Th6.4} reduces matters to the `non-linear estimate'
\begin{equation}
\label{7.2(2)}%%This is the second (7.2) in this paper.  Henceforth I will guess which reference refers to which.
||\Psi.w^k.\partial_xw||_{L_{2T}^2H^s({\mathbb R})} \leq C\Lambda(w)^{k + 1},
\end{equation}
which we proceed to establish.  The left-hand side of \pref{7.2(2)} is controlled by $$||\Psi.w^k.\partial_xw||_{L_{2T}^2L_x^2} + ||\Psi.D_x^s(w^k.\partial_xw)||_{L_{2T}^2L_x^2} = I + II.$$  For $I$, we use the bound $||w^{k - 2}||_{L_x^\infty} \leq C\snm{w}{s}{x}^{k - 2}$, which follows from Sobolev embedding and the fact that $s > 1/2$, to estimate $$I \leq C\lambda_2(w)^{k - 2}||\Psi.w^2.\partial_xw||_{L_{2T}^2L_x^2} \leq C\lambda_2(w)^{k - 2}\lambda_6(w)^2\lambda_4(w) \leq C\Lambda(w)^{k + 1}.$$  In order to bound $II$, we use the Leibniz rule, Theorem A.8 of \cite{10}%%Or do we really want [2]?
, and the chain rule, Theorem A.6 of \cite{10}, to find that
\begin{multline*}
II \leq ||D_x^s(w^k.\partial_xw)||_{L_t^2L_x^2} \leq ||D_x^s(w^k.\partial_xw) - w^k.D_x^s\partial_xw - D_x^s(w^k)\partial_xw||_{L_x^2L_t^2} \\
+ ||w^kD_x^s\partial_xw||_{L_x^2L_t^2} + ||D_x^s(w^k)\partial_xw||_{L_x^2L_T^2} =: II_1 + II_2 + II_3.
\end{multline*}
We have $$II_2 \leq C\lambda_2(w)^{k - 2}||w^2.D_x^s\partial_xw||_{L_x^2L_t^2} \leq C\lambda_2(w)^{k - 2}\lambda_6(w)^2\lambda_3(w).$$  For $II_3$, we use H\"older's inequality to see that
\begin{multline*}
II_3 \leq ||\partial_xw||_{L_x^{4(s + 1)/s}L_t^{2(s + 1)}}||D_x^s(w^k)||_{L_x^{4(s + 1)/(s + 2)}L_t^{2(s + 1)/s}} \\
\overset{*}{\leq} C||\partial_xw||_{L_x^{4(s + 1)/s}L_t^{2(s + 1)}}||D_x^s(w)||_{L_x^{4(s + 1)}L_t^{2(s + 1)/s}}||w^{k - 1}||_{L_x^4L_t^\infty} \\
C\lambda_2(w)^{k - 2}\lambda_6(w)||\partial_xw||_{L_x^{4(s + 1)/s}L_t^{2(s + 1)}}||D_x^sw||_{L_x^{4(s + 1)}L_t^{2(s + 1)/s}}.
\end{multline*}
We have used Theorem A.6 of \cite{10} to justify the inequality marked with a $*$.  Note, however, that $$||\partial_xw||_{L_x^{4(s + 1)/s}L_t^{2(s + 1)}} \leq C\lambda_3(w)^\theta\lambda_6(w)^{1 - \theta},$$ where $\theta = 1/(s + 1)$, and $$||D_x^s(w)||_{L_x^{4(s + 1)}L_t^{2(s + 1)/s}} \leq C\lambda_3(w)^{\theta'}\lambda_6(w)^{1 - \theta'},$$ where $\theta' = s/(s + 1)$, by a well-known variant of the three-lines theorem (see, for instance, Lemma 4.2 in Chapter V of \cite{17}).  (In the case of the first inequality, we also need to pass from $D_xw$ to $\partial_xw$, which involves the boundedness of the Hilbert transform in mixed-norm spaces.)  The end result of all this is that $II_3 \leq C\Lambda(w)^{k + 1}$.  Finally, Theorem A.8 of \cite{10} shows that $$II_1 \leq C||D_x^s(w^k)||_{L_x^{4(s + 1)/(s + 2)}L_t^{2(s + 1)/s}}||\partial_xw||_{L_x^{4(s + 1)/s}L_t^{2(s + 1)}},$$ and, proceeding as in the proof of the bound for $II_3$, we arrive at \pref{7.2(2)}.
\end{proof}

We now give, for each $k \geq 2$, a precise result, following \cite{10}.  We will carry out the details only in the cases $k = 2$ and $k = 4$.

\subsubsection{The case $k = 2$, modified KdV}

\begin{theorem}
\label{Th7.2}
For $k = 2$ and $1/4 \leq s < 1/2$, given $(f, \phi) \in H^{(s + 1)/3}({\mathbb R}^+) \times H^s({\mathbb R}^+)$, there exists $$T = T(||(f, \phi)||_{H^{(s + 1)/3}({\mathbb R}^+) \times H^s({\mathbb R}^+)})$$ such that the `integral equation' \pref{7.2} (for $k = 2$) has a unique fixed point in the space $B$ defined in Theorem \ref{Th7.1}.  The resulting $w$ solves \pref{7.1} in ${\mathbb R}^+ \times [0, T]$.
\end{theorem}

The proof is identical to that of Theorem \ref{Th7.1}.  The only place where $s > 1/2$ played a role there was in the bound $||w^{k - 2}||_{L_x^\infty L_T^\infty%%Or $L_t^2$?
} \leq C\lambda_2(w)^{k - 2}$, which, when $k = 2$, is not needed.  Note that, even in the half-plane case, it is known that one cannot take $s < 1/4$ (\cite{14'}).  In our proof, this is reflected in the exponent $s - 1/4$ in \pref{4.5}.

\begin{remark}
\label{Rk7.1}
One can construct solutions in Theorems \ref{Th7.1} and \ref{Th7.2} for arbitrarily large $T$ by making $||(f, \phi)||_{H^{(s + 1)/3} \times H^s}$ small.  We need only replace, for $T > 1$, $\HS_1$ by $\HS_{2T}$.
\end{remark}

\subsubsection{The case $k = 4$}

We have chosen to do in detail the case $k = 4$ because its proof is slightly simpler, since it does not involve fractional derivatives.  (See \cite{10} and \cite{13} for proofs in the half-plane case.)

\begin{definition}
\label{Df7.1}
We define $\HS_\infty$ and $\IHS_\infty$ to be the operators arising in Remarks \ref{Rk6.1} and \ref{Rk6.2}, respectively.
\end{definition}

\begin{theorem}
\label{Th7.3}
For $k = 4$, given $(f, \phi) \in \dot{H}^{1/3}({\mathbb R}^+) \times L^2({\mathbb R}^+)$ such that $$||(f, \phi)||_{\dot{H}^{1/3} \times L^2} \leq \delta$$ (with $\delta > 0$ an absolute constant), the integral equation
\begin{equation}
\label{7.3}
w = \HS_\infty(f, \phi) + \IHS_\infty(w^4%%It's not exactly clear that this is what the text says, but I think it is.
\partial_xw)
\end{equation}
has a unique fixed point in the space $$B = \bigl\{w \in C\bigl((-\infty, +\infty); \dot{H}^{1/3}({\mathbb R})\bigr) \cap C\bigl((-\infty, +\infty); L^2({\mathbb R})\bigr) : \Lambda(w) < \infty\},$$ where $\Lambda(w) = \max\inlim_{1 \leq i \leq 4} \lambda_i(w)$ for
\begin{align*}
\lambda_1(w) & = \sup_x \adnm{w(x, -)}{1/3}{\mathbb R}, \\
\lambda_2(w) & = \sup_t ||w(-, t)||_{L^2({\mathbb R})}, \\
\lambda_3(w) & = ||\partial_xw||_{L_x^\infty L_t^2}     \\
\intertext{and}
\lambda_4(w) & = ||w||_{L_x^5L_t^{10}}.
\end{align*}
The resulting $w$ solves \pref{7.1} (for $k = 4$) in ${\mathbb R}^+ \times (0, \infty)$.
\end{theorem}

\begin{proof}{}
Note that \pref{6.11} shows that $\Lambda\bigl(\HS_\infty(f, \phi)\bigr) \leq C\delta$.  Moreover, note that, if $w \in B$, then H\"older's inequality shows that $w^4\partial_xw \in L_x^{5/4}L_t^{10/9}$.  Thus, Remark \ref{Rk6.2} shows that $\Lambda\bigl(\IHS_\infty(w^4\partial_xw)\bigr) \leq C\Lambda(w)^5$, and so, if $B_a = \{w \in B : \Lambda(w) \leq a\}$ and $w \in B_a$, we have that $$\Lambda\bigl(\HS_\infty(f, \phi) + \IHS_\infty(w^4\partial_xw)\bigr) \leq C\delta + Ca^5 \leq C\delta + \frac{1}{2}a$$ whenever $Ca^4 < 1/2$.  If we now choose $\delta$ so that $C\delta \leq \frac{1}{2}a$, our mapping sends $B_a$ into $B_a$.  Similar reasoning gives that it is a contraction, and this establishes the theorem.
\end{proof}

\begin{theorem}
\label{Th7.4}
For $k = 4$, given $(f, \phi) \in H^{1/3}({\mathbb R}^+) \times L^2({\mathbb R}^+)$, there exists $T = T(f, \phi) < 1$ such that the `integral equation' \pref{7.4} below has a unique fixed point in the space $B_T$ defined below.  The resulting $w$ solves \pref{7.1} in ${\mathbb R}^+ \times (0, T)$.  Here,
\begin{equation}
\label{7.4}
w(x, t) = \Psi(t)\HS_1(f, \phi)(x, t) + \IHS_1(\Psi.w^4\partial_xw)(x, t),
\end{equation}
where $\Psi(t) = 1$ for $|t| \leq T$ and $\supp \Psi \subseteq \{t : |t| < 2T\}$, and $$B_T = \bigl\{w \in C\bigl([-2T, 2T]; L^2({\mathbb R})\bigr) \cap C\bigl((-\infty, +\infty); H^{1/3}((-2T, 2T))\bigr) : \Lambda_T(w) < \infty\},$$ where $\Lambda_T(w) = \max\inlim_{i \leq 4} \lambda_i(w)$ for
\begin{align*}
\lambda_1(w) & = \sup_x \anm{w - \Psi.\HS_1(f, \phi)}{1/3}{(-2T, 2T)}, \\
\lambda_2(w) & = \sup_t ||w - \Psi.\HS_1(f, \phi)||_{L^2},             \\
\lambda_3(w) & = ||\partial_xw||_{L_x^\infty L_{2T}^2}                \\
\intertext{and}
\lambda_4(w) & = ||w||_{L_x^5L_{2T}^{10}}.
\end{align*}
\end{theorem}

\begin{proof}{}
The extra ingredients in this proof are the following estimates:  Given $(f_0, \phi_0) \in H^{1/3}({\mathbb R}^+) \times L^2({\mathbb R}^+)$, for any $\varepsilon > 0$ there exists $T = T(f_0, \phi_0)$ and $\delta = \delta(f_0, \phi_0)$ such that, if $||(f, \phi) - (f_0, \phi_0)||_{H^{1/3} \times L^2} \leq \delta$, then
\begin{gather}
\label{7.5}||\partial_x\HS_1(f, \phi)||_{L_x^\infty L_T^2} < \varepsilon \\
\intertext{ and }
\label{7.6}||\HS_1(f, \phi)||_{L_x^5L_T^{10}} < \varepsilon.
\end{gather}
To establish \pref{7.5}, note that, by \pref{6.2}, it suffices to show it for $(f_0, \phi_0)$, with $\varepsilon$ replaced by $\varepsilon/2$.  Next, note that the corresponding result for $S(t)\tilde{\phi}_0$ is (5.11) in \cite{10}.  We thus have to show that, if $h \in H^{1/3}_0({\mathbb R}^+)$ and $\supp h \subseteq [0, 1]$, then the corresponding estimate holds for $\int_0^t S(t - t')(\delta_0)(x)h(t')dt'$.  Pick now $h_j \in C_0^\infty\bigl((0, 1)\bigr)$ satisfying $h_j \rightarrow h$ in $H^{1/3}$.  Because of \pref{4.7} and Propositions \ref{Pn2.7} and \ref{Pn2.8}, it suffices to establish the result for $$w_0(x, t) = \int_0^t S(t - t')(\delta_0)(x)h_{j_0}(t')dt',$$ with $j_0$ fixed large.  Since, however, $h_{j_0} \in C_0^\infty\bigl((0, 1)\bigr)$, $$\pd{w_0}{t}(x, t) = \int_0^t S(t - t')(\delta_0)(x)\pd{h_{j_0}}{t'}(t')dt',$$ and so $\bigl|\bigl|\pd{}{x}\pd{}{t}w_0\bigr|\bigr|_{L_x^\infty L_t^2} \leq C_{j_0},$ by \pref{4.7}.  Then, since $w_0(-, 0) \equiv 0$, $$\int_0^T \biggl|\pd{w_0}{x}(x, t)\biggr|^2dt \leq \int_0^T \biggl|\int_0^t \pd{}{t'}\pd{}{x}w_0(x, t')dt'\biggr|^2dt \leq C_{j_0}T^2,$$ and our result follows.  Next note that \pref{7.6} follows from a similar but simpler argument.

Note that \pref{7.5} and \pref{7.6} imply that $\Lambda_T\bigl(\Psi.\HS_1(f, \phi)\bigr) \leq \varepsilon$, where $\varepsilon$ is any fixed positive number, for $T$ sufficiently small.  Also, note that $$\Lambda_T\bigl(\IHS_1(\Psi.w^4\partial_xw)\bigr) \leq C||\Psi.w^4\partial_xw||_{L_x^{5/4}L_t^{10/9}} + C\varepsilon,$$ by virtue of Theorem \ref{Th6.6}, and that H\"older's inequality shows that $||\Psi.w^4\partial_xw||_{L_x^{5/4}L_t^{10/9}}\linebreak \leq \Lambda_T^5(w)$.  Thus $$\Lambda_T\bigl(\Psi.\HS_1(f, \phi) + \IHS_1(\Psi.w^4\partial_xw)\bigr) \leq C\varepsilon + C\Lambda_T^5(w),$$ so that, for $a$ satisfying $Ca^4 \leq 1/2$ and $C\varepsilon \leq a/2$, we have that our mapping sends $B_{T, a} = \{w : \Lambda_T(w) \leq a\}$ into itself, and a similar argument yields the contraction property.  This establishes Theorem \ref{Th7.4}.
\end{proof}

\begin{remark}
\label{Rk7.2}
In a manner similar to the one used to prove Corollary 2.11 in \cite{10}, working now with the integral equation
\begin{equation}
\label{7.7}
w = \HS_1(f, \phi) + \IHS_1(\Psi.w^4\partial_xw),
\end{equation}
one can show that, for $k = 4$ and $0 < s < 1/2$, \pref{7.1} is well-posed in $H^{(s + 1)/3} \times H^s$ in an interval $[0, T]$, with $$T = T(||(f, \phi)||_{H^{(s + 1)/3} \times H^s}).$$  Moreover, if $||(f, \phi)||_{H^{(s + 1)/3} \times H^s}$ is small and $0 < s < 1/2$, one can construct solutions for all $T$.
\end{remark}

\begin{remark}
\label{Rk7.3}
Following the proof of Theorem 2.6 in \cite{10} and the ideas used in the proof of Theorem \ref{Th7.2}, one can show local well-posedness for \pref{7.1} (for $k = 3$) in $H^{(s + 1)/3} \times H^s$, with $1/12 \leq s < 1/2$.  Moreover, following the proofs of Theorems \ref{Th7.3} and \ref{Th7.4}, Remark \ref{Rk7.2} and Theorems 2.15 and 2.17 and Corollary 2.18 in \cite{10}, we can extend Theorems \ref{Th7.3} and \ref{Th7.4} and Remark \ref{Rk7.2} to the case of $k > 4$, where the role of $L^2$ is played by $H^{s_k}$ and the role of $H^{1/3}$ is played by $H^{(s_k + 1)/3}$ (where $s_k = \frac{1}{2} - \frac{2}{k}$).
\end{remark}

\subsection{Local well-posedness of KdV using Bourgain's spaces}

We next turn our attention to the bilinear case, $k = 1$.  Here we will obtain results analogous to Bourgain's in \cite{5}.  We start out with the main bilinear estimate, Appendix 2 in \cite{5}.

\begin{lemma}
\label{Lm7.1}
For $v, w \in X_b$, we have, for suitable $b < \frac{1}{2}$ and $\alpha > \frac{1}{2}$, the bilinear estimate
\begin{equation}
\label{7.8'}
||\partial_x(vw)||_{Y_b} \leq C||v||_{X_b}||w||_{X_b}.
\end{equation}
\end{lemma}

\begin{proof}
We recall three fundamental estimates.  If $\hat{F}_\rho(\xi, \lambda) = \frac{f(\xi, \lambda)}{(1 + |\lambda - \xi^3|)^\rho}$, then, for $\rho > \frac{3}{8}$ and $0 \leq \theta \leq \frac{1}{8}$, we have (\cite{11'})
\begin{equation}
\label{7.9'}
||D_x^\theta F_\rho||_{L_x^4L_t^4} \leq C||f||_{L_\lambda^2}.
\end{equation}

If $\rho > \frac{1}{4}$, we have
\begin{equation}
\label{7.10'}
||D_x^{1/2}F_\rho||_{L_x^4L_t^2} \leq C||f||_{L_\lambda^2}.
\end{equation}
The estimate \pref{7.10'} is an interpolant between the local smoothing estimate and a trivial estimate; see \cite{5}.

If $\hat{H}_\rho(\xi, \lambda) = \frac{f(\xi, \lambda)}{(1 + |\lambda|)^\rho}$ with $\rho > \frac{1}{2}$, then Sobolev's inequality implies
\begin{equation}
\label{7.11'}
||H_\rho||_{L_x^2L_t^\infty} \leq C||f||_{L_{\xi\lambda}^2}.
\end{equation}

To prove \pref{7.8'}, we estimate the three pieces of the $Y_b$ norm appearing in \pref{5.2}.  Before turning to the analysis of the first piece, we introduce notation related to the $X_b$ norm, allowing us to reexpress \pref{7.8'} with $L^2$ norms on the right side.  Let $\chi = \chi_{[-1, 1]}$ and define $\beta(\xi, \lambda) = (1 + |\lambda - \xi^3|)^b + \chi(\xi)|\lambda|^\alpha$.  Introduce $g_1(\xi, \lambda) = \beta(\xi, \lambda)\hat{u}(\xi, \lambda)$ and $g_2(\xi, \lambda) = \beta(\xi, \lambda)\hat{v}(\xi, \lambda)$.  Note that $g_1 \in L_{\xi\lambda}^2 \Leftrightarrow u \in X_b$.

By duality, the %{\bf
first term in \pref{5.2}%}
is appropriately controlled if we show
\begin{equation}
\label{7.12'}
\int\blim_{\substack{\xi = \xi_1 + \xi_2 \\ \lambda = \lambda_1 + \lambda_2}} \frac{d(\xi, \lambda)|\xi|}{(1 + |\lambda - \xi^3|)^b}\frac{\hat{g}_1(\xi_1, \lambda_1)}{\beta(\xi_1, \lambda_1)}\frac{\hat{g}_2(\xi_2, \lambda_2)}{\beta(\xi_2, \lambda_2)} \leq C||d||_{L^2}||g_1||_{L^2}||g_2||_{L^2}.
\end{equation}
Symmetry allows us to assume that $|\xi_2| \geq |\xi_1|$.  Without loss of generality, we may assume $d, g_1, g_2 \geq 0$.

\begin{list}{Case \Alph{enumi}.}{\setlength{\leftmargin}{\leftmargini}\usecounter{enumi}}
\item $|\xi| \leq 1$.

In this region, \pref{7.12'} is bounded by $$\int\blim_{\substack{\xi = \xi_1 + \xi_2 \\ \lambda = \lambda_1 + \lambda_2}} d(\xi, \lambda)\frac{\hat{g}_1(\xi_1, \lambda_1)}{(1 + |\lambda_1 - \xi_1^3|)^b}\comment{Should it really be $|\lambda - \xi^3|$ instead of $|\lambda_1 - \xi_1^3|$?}\frac{\hat{g}_2(\xi_2, \lambda_2)}{(1 + |\lambda_2 - \xi_2^3|)^b}.$$  Let $\hat{D}(\xi, \lambda) = d(\xi, \lambda)$ and $\hat{G}_i(\xi, \lambda) = \frac{\hat{g}_i(\xi, \lambda)}{(1 + |\lambda - \xi^3|)^b}$.  Fourier transform properties permit us to rewrite the above expression as $\int D.G_1.G_2$, which we bound using H\"older by $||D||_{L_{xt}^2}||G_1||_{L_{xt}^4}||G_2||_{L_{xt}^4}$.  We apply \pref{7.9'} with $\theta = 0$ to the $L_{xt}^4$ norms; this requires assuming $\frac{3}{8} < b$.  Note that this case also addresses the second term in \pref{5.2}
\item $|\xi| > 1$.

The convolution constraints $\xi = \xi_1 + \xi_2$ and $\lambda = \lambda_1 + \lambda_2$ imply
\begin{equation}
\label{7.13'}
M = \max\{|\lambda - \xi^3|, |\lambda_1 - \xi_1^3|, |\lambda_2 - \xi_2^3|\} > |\xi|\,|\xi_1|\,|\xi_2|.
\end{equation}
\begin{list}{Case \Alph{enumi}\arabic{enumii}.}{\setlength{\leftmargin}{\leftmarginii}\usecounter{enumii}}
\item $|\xi_1| \geq 1$.  (Recall that we have assumed $|\xi_2| \geq |\xi_1|$.)

In this case, $|\xi|\,|\xi_1|\,|\xi_2| \geq \frac{1}{2}|\xi|^2$, so the maximum appearing in \pref{7.13'} is powerful.  We consider three subcases of this case.
\begin{list}{Case \Alph{enumi}\arabic{enumii}\alph{enumiii}.}{\setlength{\leftmargin}{\leftmarginiii}\usecounter{enumiii}}
\item $|\lambda - \xi^3| = M$.

We bound \pref{7.12'} by $$C\ic\int\blim_{\substack{\xi = \xi_1 + \xi_2 \\ \lambda = \lambda_1 + \lambda_2}} d(\xi, \lambda)\frac{\hat{g}_1(\xi_1, \lambda_1)}{(1 + |\lambda_1 - \xi_1^3|)^b}\frac{|\xi_2|^{1 - 2b}\hat{g}_2(\xi_2, \lambda_2)}{(1 + |\lambda_2 - \xi_2^3|)^b}.$$  This may be rewritten as $\int D.G_1.D_x^{1 - 2b}G_2$, which we estimate as \linebreak $||D||_{L^2}||G_1||_{L^4}||D_x^{1 - 2b}G_2||_{L^4}$.  Then \pref{7.9'} proves \pref{7.8'} in this case, provided $1 - 2b \leq \frac{1}{8}$ and $\frac{3}{8} < b$; that is, provided $b \geq \frac{7}{16}$.  We require henceforth that $b \geq \frac{7}{16}$.
\item $|\lambda_1 - \xi_1^3| = M$.

We bound \pref{7.12'} by $$C\ic\int\blim_{\substack{\xi = \xi_1 + \xi_2 \\ \lambda = \lambda_1 + \lambda_2}} \frac{|\xi|^{1 - 2b}d(\xi, \lambda)}{(1 + |\lambda - \xi^3|)^b}\hat{g}_1(\xi_1, \lambda_1)\frac{\hat{g}_2(\xi_2, \lambda_2)}{(1 + |\lambda_2 - \xi_2^3|)^b},$$ and this may also be estimated using \pref{7.8'}.
\item $|\lambda_2 - \xi_2^3| = M$.

This case is similar.
\end{list}
\item $|\xi_1| \leq 1$.

We write $|\xi| = |\xi|^{1/2}|\xi|^{1/2} \leq C|\xi|^{1/2}|\xi_2|^{1/2}$ and bound \pref{7.12'} by $$C\ic\int\blim_{\substack{\xi = \xi_1 + \xi_2 \\ \lambda = \lambda_1 + \lambda_2}} \frac{|\xi|^{1/2}d(\xi, \lambda)}{(1 + |\lambda - \xi^3|)^b}\frac{\hat{g}_1(\xi_1, \lambda_1)\chi(\xi_1)}{1 + \chi(\xi_1)|\lambda_1|^\alpha}\frac{|\xi_2|^{1/2}\hat{g}_2(\xi_2, \lambda_2)}{(1 + |\lambda_2 - \xi_2^3|)^b}.$$  This may be reexpressed as $\int\tbic\int D_x^{1/2}D.H_\alpha.D_x^{1/2}G_2$, where $\hat{H}_\alpha(\xi, \lambda) = \frac{\hat{g}_1(\xi, \lambda)}{(1 + |\lambda|)^\alpha}$.  Applying H\"older, we bound by $C||D_x^{1/2}D||_{L_x^4L_t^2}||H_\alpha||_{L_x^2L_t^\infty}||D_x^{1/2}G_2||_{L_x^4L_t^2}$, which is controlled using \pref{7.10'} and \pref{7.11'}, since $\alpha > \frac{1}{2}$ and $b > \frac{1}{4}$.
\end{list}
\end{list}
All that remains is the third piece of the $Y_b$ norm appearing in \pref{5.2}.  Observe that $$\int \biggl(\int \frac{\hat{f}(\xi, \lambda)}{1 + |\lambda - \xi^3|}d\lambda\biggr)^2d\xi = \int \biggl(\int \frac{\hat{f}(\xi, \lambda)}{(1 + |\lambda - \xi^3|)^b}\frac{1}{(1 + |\lambda - \xi^3|)^{1 - b}}d\lambda\biggr)^2d\xi.$$  By Cauchy-Schwarz, this is bounded by $$\int \biggl(\int \frac{|\hat{f}(\xi, \lambda)|^2}{(1 + |\lambda - \xi^3|)^{2b}}d\lambda\biggr)\biggl(\int \frac{1}{(1 + |\lambda - \xi^3|)^{2(1 - b)}}d\lambda\biggr)d\xi.$$ Since $2(1 - b) > 1$, the second $\lambda$-integral is bounded and we reduce matters to controlling the first term in \pref{5.2}.
\end{proof}

The next result establishes a local well-posedness result for the $k = 1$ case of \pref{1.1} when the data $(f, \phi) \in H_t^{1/3} \times L_x^2$ are sufficiently small.

\begin{theorem}
\label{Th7.5}
There exists a $\delta_0 > 0$ such that, if $(f, \phi) \in H_t^{1/3} \times L_x^2$ satisfies
\begin{equation}
\label{7.14'}
||(f, \phi)||_{H_t^{1/3} \times L_x^2} \leq \delta_0,
\end{equation}
then the ``integral equation''
\begin{equation}
\label{7.15'}
w = \HS_1(f, \phi) + \IHS_1\bigl(\partial_x(w^2/2)\bigr)
\end{equation}
has a unique fixed point in the space
\begin{equation}
\label{7.16'}
C\bigl((-\infty, \infty); H_x^s) \cap C\bigl((-\infty, \infty); H_t^{(s + 1)/3}) \cap X_b\comment{In the text it is $\{w : ||w||_{X_b} < \infty\}$; but isn't this what $X_b$ is?}
\end{equation}
(for suitable $b < \frac{1}{2}$ and $\alpha > \frac{1}{2}$).  The resulting $w$ solves \pref{7.1} in ${\mathbb R}^+ \times (0, 1)$.
\end{theorem}

\begin{proof}
Theorem \ref{Th6.1} and Lemmas \ref{Lm5.2} and \ref{Lm5.3} show that $\HS_1(\phi, f)$ satisfies \pref{7.16'}.  Lemmas \ref{Lm5.4} and \ref{Lm7.1} allow us to prove the contraction estimate under the condition \pref{7.14'}.
\end{proof}

Note that the above Theorem extends, with almost identical proof, to the case of data $(f, \phi)$ in $H^{(s + 1)/3}({\mathbb R}^+) \times H^s({\mathbb R}^+)$, $0 < s \leq 1$, $s \leq \frac{1}{2}$ (with the compatibility condition $f(0) = \phi(0)$ when $s > \frac{1}{2}$), and $||(f, \phi)||_{H^{(s + 1)/3}({\mathbb R}^+) \times H^s({\mathbb R}^+)} \leq \delta_s$.  We just need to use the $X_{s, b}$ spaces in Remark \ref{Rk5.1}, and show the corresponding estimates.

\begin{remark}
\label{Rk7.4}
There are two possible approaches to eliminating the restriction on the size of $(f, \phi)$.  The first one, as in \cite{5}, leads to considering functions on $[-T, T]$ and seeing that, as $T$ is small, a power of $T$ is gained in the above estimates.  The other approach is simply to scale things down:  $u$ solves $$\partial_tu + \partial_x^3u + u\partial_xu = 0$$ if and only if $u_\lambda(x, t) = \lambda^2u(\lambda x, \lambda^3t)$ does.  One can then choose a small $\lambda$ so that the datum $(f_\lambda, \phi_\lambda)$ corresponding to $u_\lambda$ has norm smaller than $\delta_s$.  (The resulting $\lambda$ depends only on $s$ and $||(f, \phi)||_{H^{(s + 1)/3} \times H^s}$.)  Once this is done, one uses the fact that $u_\lambda$ is defined for $0 < t < 1$, and hence that $u$ is defined for $0 < t < 1/\lambda^3$, to obtain a solution to \pref{7.1} in ${\mathbb R}^+ \times (0, T)$, with $$T = T(||(f, \phi)||_{H^{(s + 1)/3} \times H^s}).$$
\end{remark}

\nw{\begin{remark}
\label{Rk7.5}
The results just explained easily extend to equations with drift terms $c\partial_x$ by simply noting that the estimate in Lemma \ref{Lm7.1} also holds in the spaces $\tilde{X}_b$ introduced in Remark \ref{Rk6.4}.  This is because the key estimate \pref{7.13'} remains unchanged if we replace $\xi^3$, $\xi_1^3$ and $\xi_2^3$ by $\xi^3 + c\xi$, $\xi_1^3 + c\xi_1$ and $\xi_2^3 + c\xi_2$, respectively, under the convolution constraint $\xi = \xi_1 + \xi_2$.
\end{remark}}

\subsection{Global well-posedness for $1 \leq k \leq 4$}

When $f \equiv 0$, it is easy to see from Theorem \ref{Th7.5} and Remark \ref{Rk7.4} (using the argument in the proof of \pref{6.14} to obtain the {\it a priori} bound $$||u(-, t)||_{L^2({\mathbb R}^+)} \leq ||u(-, 0)||_{L^2({\mathbb R}^+)}$$ for sufficiently smooth solutions $u$, and an approximation argument using the solutions constructed in Theorem \ref{Th7.5}) that, for each $0 < T < \infty$, one can construct solutions as in Theorem \ref{Th7.5} (for $s = 0$) for \pref{7.1} (for $k = 1$) with $\phi \in L^2({\mathbb R})$.  In general, we have:

\begin{theorem}
\label{Th7.6}
Given $\phi \in L^2({\mathbb R}^+)$, $f \in H^{7/12}({\mathbb R}^+)$ and $T > 0$, the solutions constructed in Remark \ref{Rk7.4} (for $s = 0$) can be extended to the interval $(0, T)$.
\end{theorem}

To obtain the theorem, it suffices to establish the {\it a priori} estimate
\begin{equation}
\label{7.11}
\sup_{0 \leq t \leq T} ||u(-, t)||_{L^2({\mathbb R}^+)} \leq C_T(||(f, \phi)||_{H^{7/12}({\mathbb R}^+) \times L^2({\mathbb R}^+)}).
\end{equation}
To obtain this estimate, we use a device used in Proposition 5.4 of \cite{1}.  In fact, let $v$ be the solution to the linear problem
\begin{equation}
\label{7.12}
\begin{cases}
\partial_tv + \partial_x^3v = 0 & \text{in }{\mathbb R}^+ \times (0, T) \\
v|_{x = 0} = f              &                                       \\
v|_{t = 0} = \phi_0.        &
\end{cases}
\end{equation}
where $\phi_0$ is chosen so that $\phi_0%%Should this be here?
(0) = f(0)$ and $\anm{\phi_0}{3/4}{\mathbb R} \leq C\nm{f}{7/12}$, and $v$ is constructed in Theorem \ref{Th6.1} (for $s = 3/4$).  We will presently show that $v$ satisfies, in addition to the estimates in Theorem \ref{Th6.1}, the inequality
\begin{equation}
\label{7.13}
||\partial_xv||_{L_T^4L_x^\infty} \leq C||(f, \phi_0)||_{H^{7/12} \times H^{3/4}}
\end{equation}
(here we see the reason for the exponent $\frac{7}{12} = \frac{3/4 + 1}{3}$).  Let $w = u - v$, so that $w$ satisfies
\begin{equation}
\label{7.14}
\begin{cases}
\partial_tw + \partial_x^3w + \partial_x(w^2/2) + \partial_x(wv) + \partial_x(v^2/2) = 0 \\
w|_{x = 0} = 0 \\
w|_{t = 0} = \phi - \phi_0.
\end{cases}
\end{equation}
It clearly suffices, in view of Theorem \ref{Th6.1}, to establish the analogue of \pref{7.11} for $w$.  Multiplying the %%`The'?
equation in \pref{7.14} by $2w$ and integrating by parts in a manner similar to the proof of \pref{6.14}, and using the boundary conditions in \pref{7.14}, we obtain:
\begin{multline}
\label{7.15}
\int\blim_{x > 0} w^2(x, t)dx + \int_0^t \partial_xw(0, t')^2dt' \\
\leq \int\blim_{x \geq 0} (\phi - \phi_0)^2dx - 2\ic\int_0^t \int\blim_{x > 0} \partial_x(wv)w - \int_0^t \int\blim_{x > 0} v\partial_xv.w.
\end{multline}
Since $2\tic\int_0^t\int\tblim_{x > 0} \partial_x(wv)w = \int_0^t \int\tblim_{x > 0} w^2.\partial_xv$, it is easy to see that the right-hand side of \pref{7.15} is bounded (in view of Theorem \ref{Th6.1} and \pref{7.13}) by
\begin{multline*}
\frac{1}{2}\sup_{0 < t < T} ||w(-, t)||_{L^2({\mathbb R}^+)}^2 + ||\phi - \phi_0||_{L^2}^2 + C_T(||(f, \phi_0)||_{H^{7/12} \times H^{3/4}}) \\
+ C_T(||(f, \phi_0)||_{H^{7/12} \times H^{3/4}})\Bigl(\int_0^t ||w(-, t')||_{L^2({\mathbb R}_x^+)}^{8/3}dt'\Bigr)^{3/4},
\end{multline*}
which (upon raising both sides of \pref{7.15} to the power $4/3$), combined with Gronwall's inequality, gives the desired bound.  To finish the proof, we sketch the argument leading to \pref{7.13}.  In view of Theorem \ref{Th6.1} and its proof, and Theorem 2.1 in \cite{12}, it suffices to show that, if $h \in C_0^\infty\bigl((0, 1)\bigr)$ and $w(x, t) = \int_0^t S(t - t')(\delta_0)(x)h(t')dt'$, then we have
\begin{equation}
\label{7.16}
||\partial_xw||_{L_t^4L_x^\infty} \leq C\dnm{h}{-1/12}.
\end{equation}
We establish \pref{7.16} by means of the well-known variant of the three-lines theorem in Lemma 4.2 of Chapter V of \cite{17}, from the estimates
\begin{gather}
\label{7.17}||{\cal H}D_x^{i\gamma}D_x^{1 + 1/4%%Or (1 + 1)/4?
}w||_{L_t^4L_x^\infty} \leq Ce^{C|\gamma|}||h||_{L^2} \\
\intertext{and}
\label{7.18}||{\cal H}D_x^{i\gamma}w||_{L_t^4L_x^\infty} \leq Ce^{C|\gamma|}\dnm{h}{-1/12%%Or is it just 1/2?
- 1/3}.
\end{gather}
\pref{7.17} follows from Proposition 3.5(2) in \cite{7}.

To establish \pref{7.18}, we need, in view of \pref{4.16}, to establish the estimate for ${\cal H}D_x^{i\gamma}D_t^{1/12 + 1/3}A(\delta_0 \otimes g)$ and for $${\cal H}D_x^{i\gamma}D_t^{1/12 + 1/3}\dic\int_{-\infty}^{+\infty} S(t - t')(\delta_0)(x)g(t')dt'$$ with $g \in L^2$.  For the second estimate, it suffices to establish it for $${\cal H}D_x^{i\gamma}D_x^{1/4}D_x^1\ic\int_{-\infty}^{+\infty} S(t - t')(\delta_0)(x)g(t')dt'.$$  This, in turn, using the argument in the proof of \pref{4.31}, reduces to the `kernel estimate'
\begin{equation}
\label{7.19}
|H_{x_0}(x, y, t, s)| \leq \frac{C}{|t - s|^{1/2}},
\end{equation}
where $$H_{x_0}(x, y, t, s) = \int D_x^{i\gamma}{\cal H}D_x^{1 + 1/4}A(x_0 - x, t - t')D_y^{-i\gamma}D_y^{1 + 1/4}A(y - x_0, t' - s)dt'.$$

Note that a calculation like the one preceding \pref{4.33} gives the formula $$H_{x_0}(x, y, t, s) = \frac{1}{3}\ic\int e^{-i(x - y)\xi}e^{i(t - s)\xi^3}|\xi|^{1/2}d\xi,$$ and the estimate now follows from Lemma 2.7 in \cite{12}.  In order to estimate \linebreak ${\cal H}D_x^{i\gamma}D_t^{1/12 + 1/3}A(\delta_0 \otimes g)$, we first estimate $$D_t^{1/12 + 1/3}A(\delta_0 \otimes g) = \int\bic\int e^{ix\xi}e^{it\tau}\frac{|\tau|^{1/12 + 1/3}}{\tau - \xi^3}\hat{g}(\tau)d\tau\,d\xi.$$  As in the proof of \pref{4.38}, we first calculate the integral over $\xi$, and we are reduced to three terms, the first of which is
\begin{multline*}
C_1\ic\int e^{it\tau}|\tau|^{1/12 + 1/3}\sgn(x\tau^{1/3})e^{-ix\tau^{1/3}}\hat{g}(\tau)\frac{d\tau}{\tau^{2/3}} \\
= C_1(\sgn x)\ic\int e^{it\tau}e^{-ix\tau^{1/3}}\sgn(\tau^{1/3})\hat{g}(\tau)\frac{d\tau}{|\tau|^{1/4}}.
\end{multline*}
The desired estimate for this term follows from Theorem 2.5 in \cite{12}.  The remaining two terms can be handled in a similar manner, in the spirit of the above argument and the proof of \pref{4.38}.  For the proof of the estimate for ${\cal H}D_x^{i\gamma}D_t^{1/12 + 1/3}A(\delta_0 \otimes g)$, we again introduce a family of operators $$T_\eta(g)(x, t) = \int_{-\infty}^{+\infty} e^{it\tau}e^{ix\tau^{1/3}\eta}\hat{g}(\tau)\frac{d\tau}{|\tau|^{1/4}},$$ which we easily check satisfies
\begin{equation}
\label{7.20}
||T_\eta(g)||_{L_t^4L_x^\infty} \leq C||g||_{L^2},
\end{equation}
as a consequence of \pref{7.19}.  We then split our operator, as in the proof of \pref{4.39}, into the sum of three operators $L_i$, $i = 1, 2, 3$.  The bounds for $L_1$ and $L_3$ easily follow from \pref{7.20}.  For $L_2$, we use the argument in \pref{4.39}, together with \pref{7.20}, to reduce matters to estimating
\begin{multline*}
\int e^{it\tau}\frac{|\tau|^{i\gamma/3}\sgn(\tau^{1/3})\hat{g}(\tau)}{|\tau|^{1/4}}\biggl(\int e^{ix\tau^{1/3}\eta}\frac{\varphi_2(\eta)}{1 - \eta}d\eta\biggr)d\tau \\
= \int e^{it\tau}|\tau|^{i\gamma/3}\frac{\sgn(\tau^{1/3})}{|\tau|^{1/4}}e^{-ix\tau^{1/3}}\biggl(\int e^{ix\tau^{1/3}\mu}\frac{\varphi_2(1 - \mu)}{\mu}d\mu\biggr)d\tau
\end{multline*}
-- but the integral in parentheses equals $\int \sgn(x\tau^{1/3} - z)\theta(z)dz$ for some $\theta \in {\cal S}({\mathbb R})$, and, if we now use Corollary 2.9 in \cite{12} to estimate $$T_2(g)%%What *is* this thing?  That's what it seems to be in the text, but this makes no sense.
(x, t) = \int e^{it\tau}e^{-ix\tau^{1/3}}\sgn(x\tau^{1/3} - z)|\tau|^{i\gamma/3}\sgn(\tau^{1/3})\hat{g}(\tau)\frac{d\tau}{|\tau|^{1/4}}$$ in $L_t^4L_x^\infty$, uniformly in $z$, and Minkowski's integral inequality, then the proof is finished.

For $k > 1$, it is well-known (see \cite{3}) that global well-posedness for $\phi$ in $H^1$ does not follow readily from Theorem \ref{Th7.1} when $f \not\equiv 0$.  (When $f \equiv 0$, the results are completely analogous to the ones obtained in \cite{10} for the whole upper half-plane, namely, for $k = 2, 3$ there is global well-posedness, and also for $k \geq 4$, provided the $L^2$ norms of the initial data are small enough.)  When $f \not\equiv 0$, we can establish:

\begin{theorem}
\label{Th7.7}
Given $T > 0$ and $\phi \in H^1({\mathbb R}^+)$, we have:
\begin{enumerate}
\item\label{Th7.7(i)} If $f \in H^{11/12}({\mathbb R}^+)$ and $||f||_{L^2({\mathbb R}^+)}$ is sufficiently small, the solution given in Theorem \ref{Th7.1} (for $k = 2$ and $s = 1$) extends to the interval $(0, T)$.
\item\label{Th7.7(ii)} If $f \in H^{5/4}({\mathbb R}^+)$ and $||f||_{L^2({\mathbb R}^+)}$ is sufficiently small, the solution given in Theorem \ref{Th7.1} (for $k = 3$ and $s = 1$) extends to the interval $(0, T)$.
\item\label{Th7.7(iii)} If $f \in H^{11/12}({\mathbb R}^+)$ and $||f||_{L^2({\mathbb R}^+)}$ and $||\phi||_{L^2({\mathbb R}^+)}$ are sufficiently small, the solution given in Theorem \ref{Th7.1} (for $k \geq 4$ and $s = 1$) extends to the interval $(0, T)$.
\end{enumerate}
\end{theorem}

We will limit ourselves to a sketch of the proof of this result.  We start from the following two identities, readily obtained by integration by parts for sufficiently nice solutions (here $w = u - v$, where $v$ solves \pref{7.12}):
\begin{multline}
\label{7.21}
\int\blim_{x > 0} u^2(x, t)dx + \int_0^t \partial_xu(0, t')^2dt' + 2\ic\int_0^t \partial_x^2u(0, t')f(t')dt' \\
- \frac{2}{k + 2}\int_0^t f(t')^{k + 2}dt' = \int\blim_{x > 0} \phi^2(x)dx
\end{multline}
and
\begin{multline}
\label{7.22}
\int\blim_{x > 0} \bigl(\partial_xw(x, t)\bigr)^2dx + \int_0^t \bigl(\partial_x^2w(0, t')\bigr)^2dt' + \int_0^t \frac{f^{2k + 2}(t')}{(k + 1)^2}dt' \\
+ \frac{2}{k + 1}\int_0^t \partial_x^2w(0, t')\!\cdot\!f^{k + 1}(t')dt' \\
- 2\ic\int_0^t\int\blim_{x > 0} \partial_x^3v\frac{u^{k + 1}}{k + 1} - \frac{2}{(k + 1)(k + 2)}\ic\int\blim_{x > 0} u^{k + 2}(x, t)dx \\
+ \frac{2}{(k + 1)(k + 2)}\ic\int\blim_{x > 0} \phi^{k + 2}(x)dx = \int\blim_{x > 0} (\partial_x\phi)^2(x)dx.
\end{multline}
From \pref{7.21} and \pref{7.22}, together with estimates on $v$, one obtains, under the conditions of Theorem \ref{Th7.7}, the {\it a priori} bound
\begin{equation}
\label{7.23}
\sup_{0 < t < T} \anm{u(-, t)}{1}{{\mathbb R}^+} \leq C_T(f, \phi),
\end{equation}
from which Theorem \ref{Th7.7} easily follows.  In fact, to obtain \pref{7.23} in case \pref{Th7.7(i)}, one notes that $$2\ic\int_0^t \int\blim_{x > 0} \partial_x^3v\frac{u^3}{3} = -2\ic\int_0^t \partial_x^2v(0, t')\frac{f^3(t')}{3}dt' - 2\ic\int_0^t \int\blim_{x > 0} \partial_x^2v.\partial_xu.u^2,$$ and that an extension of the argument used to prove \pref{7.13} gives (with $\phi_0$ chosen so that $\nm{\phi_0}{7/4} \leq C\nm{f}{11/12}$) the estimate $$||\partial_x^2v||_{L_T^4L_x^\infty} \leq C||(f, \phi_0)||_{H^{11/12} \times H^{7/4}}.$$  Applying now simple manipulations to these facts and \pref{7.21} and \pref{7.22} yields the estimate
\begin{multline}
\label{7.24}
\sup_{0 < t < T} ||u(-, t)||_{L^2({\mathbb R}^+)}^2 + \sup_{0 < t < T} ||\partial_xu(-, t)||_{L^2({\mathbb R}^+)}^2 + \int_0^T \bigl(\partial_x^2u(0, t)\bigr)^2dt \\
\leq C + C\!\sup_{0 < t < T} \int\blim_{x > 0} u^4(x, t)dx + C\ic\int_0^T |\partial_x^2u(0, t)|\,|f(t)|dt,
\end{multline}
where $C = C\bigl(T, (f, \phi)\bigr)$.  Recall now the elementary estimate $$||w||_{L^\infty({\mathbb R}^+)} \leq ||w||_{L^2({\mathbb R}^+)}^{1/2}||\partial_xw||_{L^2({\mathbb R}^+)}^{1/2}$$ ($w(0) = 0$) and use it, together with the estimates on $v$, to conclude that
\begin{multline*}
\sup_{0 < t < T} ||u(-, t)||_{L^2({\mathbb R}^+)}^2 + \sup_{0 < t < T} ||\partial_xu(-, t)||_{L^2({\mathbb R}^+)}^2 + \int_0^T \bigl(\partial_x^2u(0, t)\bigr)^2dt \\
\leq C + C\Bigl(\int_0^T |\partial_x^2u(0, t)|^2\Bigr)^{1/2}||f||_{L^2} + C\Bigl(\sup_{0 < t < T} \int\blim_{x > 0} u^2\Bigr)^3.
\end{multline*}
Let now $h(t) = \int_0^T \bigl(\partial_x^2u(0, t)\bigr)^2dt$ and insert \pref{7.21} into the right-hand side of it.  We then obtain that $$h(T) \leq C + Ch(T)^{1/2}||f||_{L^2} + Ch(T)^{3/2}||f||_{L^2}^3.$$  Since $h(0) = 0$ and $h$ is continuous, if $||f||_{L^2}^3$ is small enough, then we obtain an {\it a priori} bound for $h(T)$.  This in turn yields a bound for $\sup\inlim_{0 < t < T} ||u(-, t)||_{L^2({\mathbb R}^+)}^2$, and from this our {\it a priori} bound follows.

The cases dealt with in \pref{Th7.7(ii)} and \pref{Th7.7(iii)} are treated similarly.  The key difference is that, when the non-linearity is of a higher power, the integration by parts after \pref{7.23} is not useful, since then $\partial_xu$ appears and forces too high a power on $u$, and hence on $\partial_xu$.  This is avoided by using the fact that, when $f \in H^{5/4}({\mathbb R}^+)$, we have the estimate $||\partial_x^3v||_{L_T^4L_x^\infty} < \infty$, which allows the previous argument to be carried out without this integration by parts.  This is actually only of interest in case \pref{Th7.7(ii)}, since, in case \pref{Th7.7(iii)}, we are assuming that $||\phi||_{L^2({\mathbb R}^+)}$ is small, which allows us to revert to the method used in \pref{Th7.7(i)}.  The details are omitted.

\newpage

\bibliographystyle{amsplain}

% %\nocite{*}
 \bibliography{KdVonRplus}

\providecommand{\bysame}{\leavevmode\hbox to3em{\hrulefill}\thinspace}
\providecommand{\MR}{\relax\ifhmode\unskip\space\fi MR }
\providecommand{\MRhref}[2]{%
  \href{http://www.ams.org/mathscinet-getitem?mr=#1}{#2}
}
\providecommand{\href}[2]{#2}

\noindent University of Toronto

\noindent University of Chicago

\end{document}